\documentclass[12pt]{amsart}
\usepackage{amsmath,amssymb}

\newcommand\inv{^{-1}}

\newcommand\fh{\mathcal H}
\newcommand\fb{\mathcal B}

\newcommand\fL{\mathfrak L}

\newcommand\ff{\mathfrak F}

\newcommand\G{\Gamma}\newcommand\g{\gamma}

\newcommand\Ta{\mathbb T}

\newcommand\Ra{\mathbb R}

\newcommand\Za{\mathbb Z}

\newcommand\cu{\mathcal{U}}
\newcommand\Na{\mathbb N}
\newcommand\so{_{\omega}}
\newcommand\dk{\disp_K}

\DeclareMathOperator{\disp}{disp} 
\DeclareMathOperator{\Id}{Id} 

\DeclareMathOperator{\Diff}{Diff} 
 \DeclareMathOperator{\Vect}{Vect}
\DeclareMathOperator{\Sect}{Sect}

 \DeclareMathOperator{\Exp}{Exp}
\DeclareMathOperator{\Hom}{Hom} 
 \DeclareMathOperator{\supp}{supp}
 
\DeclareMathOperator{\Isom}{Isom} 
\DeclareMathOperator{\CAT}{CAT}

\newtheorem{theorem}{Theorem}[section]
\newtheorem{proposition}[theorem]{Proposition}
\newtheorem{lemma}[theorem]{Lemma}
\newtheorem{corollary}[theorem]{Corollary}
\newtheorem{defn}[theorem]{Definition}

\begin{document}

\title[Almost isometric actions]{Almost isometric actions, property $(T)$, and local rigidity}
\author[D. Fisher and G.Margulis]{David Fisher and Gregory Margulis}
\thanks{First author partially supported by NSF grants DMS-9902411
and DMS-0226121 and a PSC-CUNY grant. Second author partially
supported by NSF grants DMS-9800607 and DMS-0244406. The authors
would also like to thank the FIM at ETHZ for hospitality and
support.}

\begin{abstract}
Let $\Gamma$ be a discrete group with property $(T)$ of Kazhdan.
We prove that any Riemannian isometric action of $\Gamma$ on a
compact manifold $X$ is locally rigid.  We also prove a more
general foliated version of this result. The foliated result is
used in our proof of local rigidity for standard actions of higher
rank semisimple Lie groups and their lattices in \cite{FM2}.

One definition of property $(T)$ is that a group $\Gamma$ has
property $(T)$ if every isometric $\Gamma$ action on a Hilbert
space has a fixed point.   We prove a variety of strengthenings of
this fixed point properties for groups with property $(T)$. Some
of these are used in the proofs of our local rigidity theorems.
\end{abstract}

\maketitle

\section{{\bf Introduction}}

\label{sec:intro}

One of the main results of this paper is the following.

\begin{theorem}
\label{theorem:isomrigid} Let $\Gamma$ be a discrete group with
property $(T)$.  Let $X$ be a compact smooth manifold, and let
$\rho$ be a smooth action of $\Gamma$ on $X$ by Riemannian
isometries. Then the action is $C^{\infty,\infty}$ locally rigid
and $C^{k,k-\kappa}$ locally rigid for every $\kappa>0$ for
$k{>}1$.
\end{theorem}

We recall the definition of local rigidity.

\begin{defn}
\label{defn:locallyrigid} Given a locally compact group $\Gamma$
and a $\G$ action ${\rho}:\Gamma{\times}X{\rightarrow}X$ by $C^k$
diffeomorphisms on a manifold $X$, we say that the action is {\em
$C^{k,r}$ locally rigid}, where $r{\leq}k$, if any action $\rho'$
by $C^k$ diffeomorphisms, that is sufficiently $C^k$ close to
$\rho$ is conjugate to $\rho$ by a small $C^r$ diffeomorphism. We
say an action is {\em $C^{\infty,\infty}$ locally rigid} if any
action by $C^{\infty}$ diffeomorphisms which is sufficiently
$C^{\infty}$ close to $\rho$ is conjugate to $\rho$ by a small
$C^{\infty}$ diffeomorphism.
\end{defn}

\noindent {\bf Remark:} Throughout this paper, we assume that $X$,
the $\rho(\G)$ invariant metric $g$, and therefore the action
$\rho$, are much smoother than any perturbation we consider.  This
assumption is in some sense redundant: given a compact $C^k$
manifold $X$ and a $C^{k-1}$ metric $g$ on $X$, one can show that
there is a $C^{\infty}$ structure on $X$ and a $C^{\infty}$ metric
$g'$ on $X$ invariant under $\Isom(X,g)$.

We topologize the space of $C^k$ actions of $\Gamma$ by taking the
compact open topology on the space $\Hom(\Gamma, \Diff^k(X))$. The
special case of $C^{k,k}$ local rigidity is exactly local rigidity
of the homomorphism $\rho:\G{\rightarrow}\Diff^k(X)$. Since the
$C^{\infty}$ topology is defined as the inverse limit of the $C^k$
topologies, two $C^{\infty}$ diffeomorphisms are $C^{\infty}$
close if they are $C^k$ close for some large $k$.  Our proof shows
explicitly that, for any $\kappa>0$,  a $C^{\infty}$ perturbation
$\rho'$ of $\rho$ which is sufficiently $C^k$ close to $\rho$ is
conjugate to $\rho$ by a $C^{\infty}$ diffeomorphism which is
$C^{k-\kappa}$ close to the identity.  Many local rigidity results
have been proven for actions of higher rank semisimple groups and
their lattices. See the introduction to \cite{FM2} for a more
detailed historical discussion.

In fact, Theorem \ref{theorem:isomrigid} follows (though with
lower regularity) from a more general foliated version, whose
somewhat complicated statement we defer to the next section. We
also give two self-contained proofs of Theorem
\ref{theorem:isomrigid}, since many of the ideas are clearer in
that special case, and since one proof gives better regularity in
that case. Our foliated result is a principal ingredient in our
proof of local rigidity for quasi-affine actions of higher rank
semisimple Lie groups and their lattices \cite{FM1,FM2}.

Some prior results about local rigidity of isometric actions are
known. The question was first investigated for lattices $\Gamma$
in groups $G$, where $G$ is a semisimple Lie group with all simple
factors of real rank at least 2.  In \cite{Z1}, Zimmer proved that
any ergodic, volume preserving perturbation of an ergodic,
isometric actions of such $\Gamma$ preserves a $C^0$ Riemannian
metric $g$. In \cite{Z1.5}, he showed that $g$ was actually
smooth.  In \cite{Z3}, Zimmer extended his result to cover all
groups with property $(T)$ of Kazhdan, but still required that the
perturbation be ergodic and volume preserving and only constructed
an invariant metric rather than a conjugacy. In \cite{Be},
Benveniste proves $C^{\infty,\infty}$ local rigidity for isometric
actions of cocompact lattices in semisimple groups $G$ as above.
As a direct generalization of Zimmer's result, we have the
following.

\begin{theorem}
\label{theorem:invariantmetric} Let $\Gamma$ be a discrete group
with property $(T)$.  Let $(X,g)$ be a compact Riemannian manifold
and let $\rho$ be a smooth action of $\Gamma$ on $X$ preserving
$g$. For any $\kappa>0$, any $\G$ action $\rho'$ which is
sufficiently $C^{k+1}$ close to $\rho$ preserves a $C^{k-\kappa}$
Riemannian. Furthermore, if $\rho'$ is a $C^{\infty}$ action
$C^{\infty}$ close to $\rho$, then $\rho'$ preserves a
$C^{\infty}$ metric.
\end{theorem}

\noindent Though this theorem is a corollary of Theorem
\ref{theorem:isomrigid}, we give a simple direct proof of the
finite regularity case of Theorem \ref{theorem:invariantmetric} in
section \ref{section:parametrizing}.

A locally compact, $\sigma$-compact group $\Gamma$ has {\em
property $(T)$} if any continuous isometric action of $\Gamma$ on
a Hilbert space has a fixed point. In this paper we generalize
this standard fixed point property to a wider class of actions.
One can view our results as showing that this fixed point property
persists for actions which are perturbations of isometric actions.
In fact the fixed point property holds quite generally, even for
actions which are only partially defined. Note that it is a
theorem of Kazhdan that any discrete group with property $(T)$ is
finitely generated and any locally compact, $\sigma$-compact group
with property $(T)$ is compactly generated  \cite{K}.

The definition above is not Kazhdan's original definition of
property $(T)$, but is equivalent by work of Delorme and
Guichardet \cite{De,Gu}. Kazhdan defined a group $\Gamma$ to have
property $(T)$ if the trivial representation of $\Gamma$ is
isolated in the Fell topology on the unitary dual of $\Gamma$. For
detailed introductions to property $(T)$ see \cite{HV} or
\cite[Chapter III]{M}.  A key step in our proofs is to strengthen
standard fixed point properties for groups with property $(T)$.
For our foliated local rigidity theorems, we also require an
effective method for finding fixed points. One corollary of our
general method is a simpler proof of Shalom's result that any
finitely generated group with property $(T)$ is a quotient of a
finitely presented group with property $(T)$ \cite{S}. See also
\cite{Zk} for related results. We also prove a similar result for
compactly generated groups with property $(T)$, see Theorem
\ref{theorem:compactlypresented} below.



We now state a special case of our general fixed point property,
that suffices for the proof of $C^{k,k-\frac{1}{2}\dim(X)}$ local
rigidity.

\begin{defn}
\label{definition:epsilonisometry} Let $\varepsilon{\geq}0$ and
$Z$ and $Y$ be metric spaces.  Then a map $h:Z{\rightarrow}Y$ is
an {\em $\varepsilon$-almost isometry} if
$$(1-{\varepsilon})d_Z(x,y){\leq}{d_Y(h(x),h(y))}{\leq}(1+{\varepsilon}){d_Z(x,y)}$$
for all $x,y{\in}Z$.
\end{defn}

\noindent The reader should note that an $\varepsilon$-almost
isometry is a bilipschitz map.  We prefer this notation and
vocabulary since it emphasizes the relationship to isometries.

\begin{defn}
\label{definition:displacement} Given a group $\Gamma$ acting on a
metric space $X$, a compact subset $K$ of $\Gamma$ and a point
$x{\in}X$.  The number $\sup_{k{\in}K}d(x,k{\cdot}x)$ is called
the $K$-displacement of $x$ and is denoted $\dk(x)$.
\end{defn}


\begin{theorem}
\label{theorem:fixedpointsimplegen} Let $\Gamma$ be a locally
compact, $\sigma$-compact group with property $(T)$ and $K$ a
compact generating set.  There exist positive constants
$\varepsilon$ and $D$, depending only on $\Gamma$ and $K$, such
that for any continuous action of $\Gamma$ on a Hilbert space
$\fh$ where $K$ acts by $\varepsilon$-almost isometries there is a
fixed point $x$; furthermore for any $y$ in $X$, the distance from
$y$ to the fixed set is not more than $D\dk(y)$.
\end{theorem}

We note that in most of our applications, the $\varepsilon$-almost
isometric action to which we apply Theorem
\ref{theorem:fixedpointsimplegen} and its generalizations are
linear, and therefore automatically has a fixed point, the $0$
vector. The importance of the final claim in Theorem
\ref{theorem:fixedpointsimplegen} then becomes clear: we have a
linear relationship between the distance from a point to the fixed
set and the $K$-displacement of the point.   In our applications
this is used to find non-zero fixed vectors in certain linear
actions.  That the fixed vector is close to a particular vector
with particular prescribed properties is also central to the
proof.

Preliminary forms of Theorems \ref{theorem:isomrigid} and
\ref{theorem:fixedpointsimplegen} were announced by the second
author at a talk in Jerusalem in 1997.

Theorem \ref{theorem:fixedpointsimplegen} is proved by
contradiction. We assume the existence of a sequence of
$\varepsilon$-almost isometric actions not satisfying the
conclusion of that theorem, with $\varepsilon$ going to zero. One
then constructs a limit action which is isometric and therefore
must have fixed points. One then uses a quantitative strengthening
of the fixed point property to show that actions ``close enough"
to the limit action must have fixed points as well.  In this
article we use ultra-filters and ultra-limits to produce the limit
action which considerably simplifies earlier versions of the
argument. The argument is further simplified by our use of a
stronger quantitative strengthening of the fixed point property
which is, in fact, an iterative method for producing fixed points.
For $\G$ not discrete, additional difficulties arise from the fact
that the limit action is not a priori continuous.

Though the approaches and applications are different, our
strengthenings of property $(T)$ are related to the strengthenings
discussed by Gromov in \cite{Gr2}.  In particular, in section
3.13B, Gromov outlines a proof of Theorem
\ref{theorem:fixedpointsimplegen}, though only for a certain class
of ``random" infinite, discrete groups with property $(T)$ and
only for affine $\varepsilon$-almost isometric actions. See
Appendix \ref{appendix:gromov} for further discussion.

Our original approach to proving Theorem \ref{theorem:isomrigid}
remains incomplete, though the idea is instructive. Given an
isometric action $\rho$ of $\Gamma$ on a compact manifold $X$ and
a perturbation $\rho'$ of $\rho$, a conjugacy is a diffeomorphism
$f:X{\rightarrow}X$ such that
$\rho(\gamma){\circ}f=f{\circ}\rho'(\gamma)$ for all $\gamma$ in
$\Gamma$. Rearranging, the conjugacy is a fixed point for the
$\Gamma$ action on the group $\Diff^k(X)$ of diffeomorphisms of
$X$ defined by
$f{\rightarrow}{\rho(\gamma)}{\circ}f{\circ}\rho'(\gamma){\inv}$.
Ideally we would parameterize diffeomorphisms  of $X$ locally as a
Hilbert space and then use Theorem \ref{theorem:fixedpointpartial}
below, a generalization of Theorem
\ref{theorem:fixedpointsimplegen} for partially defined actions,
to find a fixed point or conjugacy.  This approach does not work,
see Appendix \ref{appendix:kam} for further discussion.

Our two proofs of Theorem \ref{theorem:isomrigid} have distinct
advantages.  We outline here the simpler one which allows us to
prove our general foliated result.  We discuss here only the
result that uses Theorem \ref{theorem:fixedpointsimplegen} and
only indicate a proof of $C^{k,k-\dim(X)-1}$ local rigidity. Even
combined with arguments below which improve regularity, this proof
requires the loss of $1+\kappa$ derivatives. The precise
regularity of Theorem \ref{theorem:isomrigid} requires a different
argument and requires stronger assumptions on the action in the
foliated case. In subsection \ref{subsection:isomrigid} we include
the other proof of Theorem \ref{theorem:isomrigid} but only
briefly indicate how, and when, it can be foliated.

Given a compact Riemannian manifold $X$, there is a canonical
construction of a Sobolev inner product on $C^k(X)$ such that the
Sobolev inner product is invariant under isometries of the
Riemannian metric,  see section \ref{section:parametrizing} below.
We call the completion of $C^k(X)$ with respect to the metric
induced by the Sobolev structure $L^{2,k}(X)$. Given an isometric
$\Gamma$ action $\rho$ on a manifold $M$ there may be no
non-constant $\Gamma$ invariant functions in $L^{2,k}(X)$.
However, if we pass to the diagonal $\Gamma$ action on
$X{\times}X$, then any function of the distance to the diagonal is
$\Gamma$ invariant and, if $C^k$, is in $L^{2,k}(X{\times}X)$.

We choose a smooth function $f$ of the distance to the diagonal in
$X{\times}X$ which has a unique global minimum at $x$ on
$\{x\}{\times}X$ for each $x$, and such that any function $C^2$
close to $f$ also has a unique minimum on each $\{x\}{\times}X$.
This is guaranteed by a condition on the Hessian and the function
is obtained from $d(x,y)^2$ by renormalizing and smoothing the
function away from the diagonal. This implies $f$
 is invariant under the diagonal $\Gamma$ action defined by $\rho$.
 Let $\rho'$ be another action
$C^k$ close to $\rho$. We define a $\Gamma$ action on $X{\times}X$
by acting on the first factor by $\rho$ and on the second factor
by $\rho'$. For the resulting action $\bar \rho'$ of $\Gamma$ on
$L^{2,k}(X{\times}X)$ and every $k{\in}K$, we show that $\bar
\rho'(k)$ is an $\varepsilon$-almost isometry and that the
$K$-displacement of $f$ is a small number $\delta$, where both
$\varepsilon$ and $\delta$ can be made arbitrarily small by
choosing $\rho'$ close enough to $\rho$. Theorem
\ref{theorem:fixedpointsimplegen} produces a $\bar \rho'$
invariant function $f'$ close to $f$ in the $L^{2,k}$ topology.
Then $f'$ is $C^{k-\dim(X)}$ close to $f$ by the Sobolev embedding
theorems and if $k-\dim(X){\geq}2$, then $f$ has a unique minimum
on each fiber $\{x\}{\times}X$ which is close to $(x,x)$. We
verify that the set of minima is a $C^{k-\dim(X)-1}$ submanifold
and, in fact, the graph of a conjugacy between the $\Gamma$
actions on $X$ defined by $\rho$ and $\rho'$.

In the context of Theorem \ref{theorem:fixedpointsimplegen}, we
can prove that given any vector $v$, one can produce any invariant
vector $v_0$ by an iterative method of ``averaging over balls in
$\G$".  The proofs of the $C^{\infty}$ cases of a Theorems
\ref{theorem:isomrigid} and \ref{theorem:invariantmetric} rely on
this iterative method and additional estimates. If our
perturbation $\rho'$ is $C^k$ close to $\rho$, using this
iterative method, convexity estimates on derivatives and estimates
on compositions we produce a sequence of $C^{\infty}$
diffeomorphisms which converge to conjugacy in the $C^l$ topology
for some $l>k$.  We then apply an additional iterative argument
loosely inspired by the KAM method, to produce the actual
$C^{\infty}$ conjugacy. For a discussion of the relation between
our work and the KAM method, see Appendix \ref{appendix:kam}. The
proof of $C^{k,k-\kappa}$ local rigidity for any $\kappa>0$ and of
the lower loss of regularity in Theorem
\ref{theorem:invariantmetric} follow from a somewhat technical
result which allows us to show that the iterative procedure
defined by ``averaging over balls" also converges in $L^p$ type
Sobolev spaces where $p>2$.   We defer statements of these results
to subsection \ref{subsection:TresultsBanach}.  Once one replaces
standard consequences of property $(T)$ with an observation of
Bader and Gelander \cite{BG}, the proof of this result is similar
to the proofs of our results concerning actions on Hilbert spaces.

The proof of the foliated generalization of Theorem
\ref{theorem:isomrigid} follows a similar outline, but is more
difficult at several steps. The choice of initial invariant
function is slightly more complicated since leaves of the
foliation are generally non-compact.   The absence of a natural
topology on the set of pairs of points on the same leaf forces us
to work on the holonomy groupoid of the foliation. Since we need
to work in a Sobolev space defined by only taking derivatives
along the leaves of the foliation, having small norm in this
topology on functions only gives a good $C^k$ estimate on the
conjugacy on a set $S$ of large measure. To guarantee that the
orbit of $S$ covers all of $X$ , we use our effective method of
producing $\tilde f$ from $f$ by ``averaging over balls" in
$\Gamma$.

\noindent {\bf Plan of the paper:} In section \ref{section:defs}
we make the necessary definitions and state our general results.
First, in subsection \ref{subsection:Tresults} we discuss various
generalizations of Theorem \ref{theorem:fixedpointsimplegen} for
actions on Hilbert spaces. Second in subsection
\ref{subsection:TresultsBanach} we discuss various generalizations
of Theorem \ref{theorem:fixedpointsimplegen} for actions on more
general Banach spaces. Then in subsection
\ref{subsection:foliatedresults} we describe our foliated
generalization of Theorem \ref{theorem:isomrigid}.  Subsection
\ref{subsection:fpstandard} and \ref{subsection:limitactions}
contain preliminaries on, respectively, groups with property $(T)$
and limits of actions.  We then proceed to prove the results from
subsections \ref{subsection:Tresults} and
\ref{subsection:TresultsBanach} in subsections
\ref{subsection:Tproofs} and \ref{subsection:TproofsBanach}
respectively. Section \ref{section:parametrizing} gives an
explicit construction of various Sobolev metrics on various spaces
of tensors on Riemannian manifolds and more general spaces with
Riemannian foliations. Section \ref{section:parametrizing} also
contains a proof of Theorem \ref{theorem:invariantmetric}. Section
\ref{section:conjugacy} contains two proofs of Theorem
\ref{theorem:isomrigid}.  In section
\ref{section:convexityofderivatives}, we prove the $C^{\infty}$
case of Theorem \ref{theorem:isomrigid}. Section
\ref{section:foliatedproofs} contains some additional background
on foliations, a discussion of the holonomy groupoid of a
foliation, and a proof of the foliated generalization of Theorem
\ref{theorem:isomrigid}.  On first reading the paper, the reader
may wish to skip subsections \ref{subsection:TresultsBanach},
\ref{subsection:foliatedresults}, \ref{subsection:TproofsBanach},
read \ref{section:parametrizing} assuming $p=2$ everywhere and
assuming that the foliation is by a single leaf and then read
subsection \ref{subsection:isomrigid}. This allows the reader to
read the proof of Theorem \ref{theorem:isomrigid} for the
$C^{k,k-\frac{\dim(X)}{2}}$ case, before beginning to study the
techniques for improving regularity and/or the, significantly more
complicated, formulation and proof of the foliated version.

{\noindent}{\it Acknowledgements:}  The authors would like to
thank:  Dmitry Dolgopyat, Bassam Fayad and Raphael Krikorian for
useful conversations concerning the KAM method and convexity of
derivatives; Uri Bader and Tsachik Gelander for sharing their
observations on property $(T)$ and $L^p$ spaces and for useful
conversations concerning Banach spaces and positive definite
functions; and Tim Riley for useful conversations concerning
ultrafilters, ultralimits and ultraproducts of metric spaces. We
also thank the referees for copious helpful comments that
considerably improved the exposition.

\section{\bf Definitions and statements of main results}
\label{section:defs}

In this section we give the necessary definitions and state our
general results. The first subsection is devoted to general
results on actions and partially defined actions of groups with
property $(T)$ on Hilbert spaces.  The second subsection concerns
generalizations of some of these results to more general Banach
spaces. The third subsection concerns the foliated version of
Theorem \ref{theorem:isomrigid}.

\medskip
\noindent{\bf On Constants:} Throughout this paper, we use a
convention to simplify the specification of which constants depend
on which other choices. When introducing a constant $C$, we will
use the notation $C=C(\alpha,\beta,S)$ to specify that $C$ depends
on choices of $\alpha,\beta$ and $S$.  We make one exception to
this rule: as most constants in this paper depend on a choice of a
group $\G$ and a generating $K$, we will always leave this
dependence implicit.  The few cases where constants do not
actually depend on an ambient choice of $\G$ and $K$ are clear
from context as they appear in statements where $\G$ and $K$ are
irrelevant.

\subsection{\bf Fixed points for actions of groups with property $(T)$ on Hilbert spaces}
\label{subsection:Tresults}

Throughout this subsection $\Gamma$ will be locally compact,
$\sigma$-compact, group generated by a fixed compact subset $K$,
which contains a neighborhood of the identity. It follows from
work of Kazhdan \cite{K} that any locally compact,
$\sigma$-compact $\Gamma$ with property $(T)$ is compactly
generated.  Given any compact generating set $C$, a simple Baire
category argument shows that $C^s$ contains a neighborhood of the
identity for some positive integer $s$. (Given a subset $K$ of a
group $\Gamma$, we write $K^s$ for the set of all elements of
$\Gamma$ that can be written as a product of $s$ elements of $K$.)






Theorem \ref{theorem:fixedpointsimplegen} suffices to prove
$C^{k,k-\frac{\dim(X)}{2}}$ local rigidity in Theorem
\ref{theorem:isomrigid}.  To obtain better finite regularity, a
$C^{\infty,\infty}$ local rigidity result, and to prove our more
general results, we will need more precise control over how one
obtains an invariant vector from an almost invariant vector.  As
noted above, most of the applications of our results are to the
case where the $\varepsilon$ almost isometric actions are actually
linear representations.  Since the statements of our results do
not simplify in any useful way in that setting, we leave it to the
interested reader to state the special cases.

Fix a (left) Haar measure $\mu_{\Gamma}$ on $\Gamma$.  We let
$\mathcal{U}(\Gamma)$ denote the set of continuous non-negative
functions $h$ with compact support on $\Gamma$ with
$\int_{\Gamma}hd\mu_{\Gamma}=1$. Given $h{\in}\mathcal{U}(\G)$ and
an action $\rho$ of $\Gamma$ on a Hilbert space $\fh$, we can
define an operator $\rho(h)$ on $\fh$.  Let
$\rho(h)v=\int_{\G}\rho(\gamma)(h(\gamma)v)d\mu_{\G}$. It is
straightforward to see that $\int_{\G}hd\mu_{\Gamma}=1$ implies
that this definition does not depend on the choice of basepoint in
$\fh$. When the action $\rho$ is not affine, $\rho(h)$ is not
necessarily an affine transformation. We let $\mathcal{U}_2(\G)$
be the subset of functions $h{\in}\mathcal{U}(\G)$ such that $h>0$
on $K^2$. We denote by $f*g$ the convolution of integrable
functions $f$ and $g$. Note that if $f,g{\in}\mathcal{U}(\G)$ then
so is $f*g$. Given a positive integer $d$, we denote by $f^{*d}$
the $d$-fold convolution of $f$ with itself.  More generally, if
$\mathcal{P}(\G)$ is the set of probability measures on $\G$ and
$\nu{\in}\mathcal{P}(\G)$, we can also define
$\rho(\nu)v=\int_{\G}\rho(\gamma)vd\nu$.  Note that
$\mathcal{U}(\G){\subset}\mathcal{P}(\G)$ and that this definition
generalizes the one above.

We now state a theorem which implies Theorem
\ref{theorem:fixedpointsimplegen}.  This theorem implies that
iterates of certain averaging operators converge to a bounded
projection onto the set of fixed points for the action.

\begin{theorem}
\label{theorem:decreasingdisplacement}  If $\G$ has property $(T)$
and $f{\in}{\mathcal{U}_2}(\G)$ and $0<C<1$, there exists
${\varepsilon}>0$, and positive integers $m=m(C,f)$ and
$M=M(C,f)$, such that, letting $h=f^{*m}$, for any Hilbert space
$\fh$, any continuous action of $\Gamma$ on $\fh$ such that $K$
acts by $\varepsilon$-almost isometries, and any $x{\in}{\fh}$ we
have
\begin{enumerate}
\item $d_{\fh}(x,\rho(h)(x)){\leq}M\dk(x)$

\item $\dk(\rho(h)(x)){\leq}C\disp_K(x)$.
\end{enumerate}
\end{theorem}

\noindent {\bf Remark:} In this theorem choosing smaller values of
$C$ increases the value of $m$. The number $M$ is the least
integer with $\supp(h)=\supp(f^{*m}){\subset}K^M$.

\begin{proof}[Proof of Theorem \ref{theorem:fixedpointsimplegen}
from Theorem \ref{theorem:decreasingdisplacement}] The hypotheses
of the theorems are almost identical. Since the $\G$ action in
Theorem \ref{theorem:fixedpointsimplegen} is continuous, it
follows that every point $x{\in}\fh$ has finite $K$-displacement.
Given a point $x{\in}{\fh}$ with $K$-displacement $\delta$, we
look at the sequence $y_n=\rho(h)^n(x)$. Theorem
\ref{theorem:decreasingdisplacement} implies that the
$\dk(y_n){\leq}C^n{\delta}$ and that
$d(y_n,y_{n+1}){\leq}MC^n{\delta}$. This implies that $y_n$ is a
Cauchy sequence and $y=\lim_{n{\rightarrow}{\infty}}y_n$ clearly
has $K$ displacement zero.  Letting
$D=\sum_{i=1}^{\infty}MC^n=\frac{MC}{1-C}$, then
$d_{\fh}(x,y)<D{\delta}$ which completes the proof.
\end{proof}





We now make precise our notion of a partially defined action.
  By $B(x, r)$ we denote the ball around $x$
of radius $r$.

\begin{defn}
\label{definition:actsonball} Let $X$ a metric space and fix a
point $x{\in}X$. Given $r,s,{\varepsilon},{\delta}>0$, we call a
map $\rho:K^s{\times}B(x,r){\rightarrow}X$ an {\em $(r, s,
{\varepsilon}, {\delta},K)$-almost action of $\G$ on $X$ at $x$}
if the following conditions hold.

\begin{enumerate}
\item For each $d{\in}K^s$, the map
$\rho(d,{\cdot}):B(x,r){\rightarrow}X$ is an $\varepsilon$-almost
isometry.

\item $\dk(x)<\delta$.

\item With the notation $\rho(d,z)=\rho(d)z$, if $ab,a$ and $b$
are in $K^s$ then $\rho(a)(\rho(b)y)=\rho(ab)y$ whenever
$\rho(b)y$ is in $B(x, r)$.

\end{enumerate}
\end{defn}

\noindent  When $K$ is fixed, we sometimes abbreviate the above
notation by calling an $(r,s,\varepsilon,\delta,K)$-almost action
an $(r,s,\varepsilon,\delta)$-almost action.  We denote by a
$({\infty}, s, {\varepsilon}, \delta,K)$-almost action the case
when $B(x,r)$ in the definition above can be replaced by $X$. The
following theorem now produces fixed points for partially defined
actions on Hilbert spaces that are "close enough" to isometric
ones.

\begin{theorem}
\label{theorem:fixedpointpartial} If $\G$ has property $(T)$ and
$\delta_0>0$ there exist ${\varepsilon}>0, D>0$, a positive
integer $s$, and $r=r(\delta_0)>0$ such that for any Hilbert space
$X$, any $\delta{\leq}\delta_0$ and any $x{\in}X$, any continuous
$(r, s, {\varepsilon}, {\delta},K)$-action of $\G$ on $X$ at $x$
has a fixed point. Furthermore, the distance from the fixed point
to $x$ is not more than $D{\delta}$.
\end{theorem}

Fixing $\delta_0$ is only necessary as a normalization. If we
compose a given action with a homothety, we may always assume
$\delta_0$ is $1$.  The constants $s$ and ${\varepsilon}$ remain
unchanged by this process, but $r$ becomes ${\delta_0}r$. The
utility of considering partially defined actions is illustrated by
our proof of the observation of Shalom stated in the introduction.
In fact, we prove the following generalization, which is used in
section \ref{section:foliatedproofs}.  For background on the
notion of a compact presentation see \cite{Ab}.

\begin{theorem}
\label{theorem:compactlypresented} Let $\Gamma$ be a locally
compact, $\sigma$-compact group with property $(T)$.  Then
$\Gamma$ is a quotient of a compactly presented locally compact,
$\sigma$-compact group with property $(T)$.
\end{theorem}

\begin{proof}
As remarked above, by work of Kazhdan, $\G$ is compactly
generated, and we fix a compact generating set $K$. Possibly after
replacing $K$ with a power of $K$, we can assume that $K$ contains
a neighborhood of the identity. The group $\Gamma$ is the quotient
of the group $\Gamma'$ generated by $K$ satisfying all relations
of $\G$ of the form $xy=z$ where $x,y,z{\in}K^s$. Since $\G'$
satisfies all the relations contained in $K$, we can topologize
$\G'$ so that the projection $\G'{\rightarrow}\G$ is a
homeomorphism in a neighborhood of the identity and therefore
$\G'$ is locally compact and $\sigma$-compact. We believe that
this fact is known, but state it as Proposition \ref{lemma:lc} in
Appendix \ref{appendix:lc} where we also sketch a proof, as we did
not find a reference in the literature. Since a continuous
isometric $\Gamma'$ action is a continuous
$(\infty,s,0,\delta)$-action of $\Gamma$ at $x$ where $\delta$ is
the $K$-displacement of $x$, Theorem
\ref{theorem:fixedpointpartial} implies that, if we choose $s$
large enough, $\Gamma'$ has property $(T)$.  It is clear that
$\Gamma'$ is compactly presented.
\end{proof}

{\noindent}{\bf Remarks:} \begin{enumerate} \item Theorem
\ref{theorem:compactlypresented} is used in the proof of Theorem
\ref{theorem:almostconjugacygen}, the foliated generalization of
Theorem \ref{theorem:isomrigid}. It is used to show that an action
of a locally compact group with property $(T)$ on a compact
foliated space lifts to an action on the holonomy groupoid of the
foliation. \item It is also possible to prove Theorem
\ref{theorem:almostconjugacygen} directly from Theorem
\ref{theorem:decreasingdisplacementgen}  and Corollary
\ref{corollary:contractingoperatorbanach} below.
\end{enumerate}

We now state a generalization of Theorem
\ref{theorem:decreasingdisplacement} which implies Theorem
\ref{theorem:fixedpointpartial}.  We note that the operator
$\rho(h)$ is well defined for a
$(r,s,\varepsilon,\delta,K)$-action $\rho$, provided the support
of $h$ is contained in $K^s$.

\begin{theorem}
\label{theorem:decreasingdisplacementgen} If $\G$ has property
$(T)$ and  $f{\in}\mathcal{U}_2(\G), 0<C<1$ and $\delta_0>0$ there
exist $r=r(\delta_0,f,C)>0$ and $\varepsilon>0$ and positive
integers $m=m(f,C),s=s(f,C)$ and $M=M(f,C)$ such that, letting
$h=f^{*m}$, for any Hilbert space $X$, any $\delta{\leq}\delta_0$,
any $x{\in}X$, and any continuous $(r, s, {\varepsilon},
{\delta},K)$-action $\rho$ of $K$ on $X$ at $x$ we have
\begin{enumerate}
\item $d_{\fh}(x,\rho(h)(x)){\leq}M\dk(x)$;

\item $\dk(\rho(h)(x)){\leq}C\dk(x)$.
\end{enumerate}
\end{theorem}

\begin{proof}[Proof of Theorem \ref{theorem:fixedpointpartial}
from Theorem \ref{theorem:decreasingdisplacementgen}] The proof is
almost identical to the proof of Theorem
\ref{theorem:fixedpointsimplegen} from Theorem
\ref{theorem:decreasingdisplacement}.  One point requires
additional care: if $r_0$ and $C$ are the constants given by
Theorem \ref{theorem:decreasingdisplacementgen}, we need to take
$r$ in Theorem \ref{theorem:fixedpointpartial} to be at least
$r_0+M\sum_{i=1}^{\infty}C_1^i$.  This insures that we can apply
Theorem \ref{theorem:decreasingdisplacementgen} to each
$\rho(h)^i(x)$ successively, since it implies that $\rho$ defines
an $(r,s,\varepsilon,C^i\delta)$-action on $\fh$ at
$\rho(h)^i(x)$.
\end{proof}

\noindent{\bf Remarks:}\begin{enumerate} \item That Theorem
\ref{theorem:decreasingdisplacementgen} implies Theorem
\ref{theorem:decreasingdisplacement} is clear from the
definitions.  Section \ref{section:Tproofs} is devoted to the
proof of Theorem \ref{theorem:decreasingdisplacementgen}.

\item For most of our dynamical applications Theorem
\ref{theorem:decreasingdisplacement} suffices.  However, as
remarked above, we need Theorem \ref{theorem:compactlypresented},
and therefore Theorem \ref{theorem:fixedpointpartial}, for the
proof of Theorem \ref{theorem:almostconjugacygen}.  As remarked
above, one can also prove Theorem \ref{theorem:almostconjugacygen}
using Theorem \ref{theorem:decreasingdisplacementgen} in place of
the combination of Theorem \ref{theorem:decreasingdisplacement}
and Theorem \ref{theorem:compactlypresented}.
\end{enumerate}

\subsection{\bf Property $(T)$ and uniformly convex Banach spaces}
\label{subsection:TresultsBanach}

In this subsection we describe some generalizations of the results
in the previous subsection to non-Hilbertian Banach spaces.
Throughout this subsection $\G$ and $K$ will be as in the previous
subsection. For $1<p{\leq}2$, we will call a Banach space $\fb$ a
{\em generalized $L^p$ space}, if the function $\|x\|^p$ is
negative definite on $\fb$ or equivalently if $\exp(-t\|x\|^p)$ is
positive definite for all $t>0$.  A theorem of Bretagnolle,
Dacunba-Castelle and Krivine implies that any generalized $L^p$
space is a closed subspace of an $L^p$ space, see \cite[Theorem
8.9]{BL}. For $q>2$, we will called a Banach space $\fb$ a {\em
generalized $L^q$ space} if the dual of $\fb$ is a generalized
$L^p$ space where $\frac{1}{p}+\frac{1}{q}=1$.   Given a finite
dimensional Euclidean space $V$ with Euclidean norm
$\|{\cdot}\|_V$ and a measure space $(S,\mu)$ we define a norm on
measurable maps $f:S{\rightarrow}V$ by
$\|f\|^p=\int_S\|f(s)\|_V^pd\mu$ and let $L^p(S,\mu,V)$ be the
space of equivalence classes of maps $f$ with finite norm.  If
$\dim(V)=n$ and $1<p<{\infty}$, we will call $L^p(S,\mu,V)$ a {\em
Banach space of type $L^p_n$.}  It is easy to verify that if
$\frac{1}{p}+\frac{1}{q}=1$ then the dual of a Banach space of
type $L^p_n$ is a Banach space of type $L^q_n$. It is also easy to
verify, for $1{<}p<\infty$,  that a Banach space of type $L^p_n$
is a generalized $L^p$ space.  For $p{\leq}2$ this is shown by
embedding $L^p(S,\mu,V)$ into
$L^p(S{\times}S_1(V),\mu{\times}\nu)$ where $S_1(V)$ is the unit
sphere in $V$ and $\nu$ is (normalized) Haar measure.  For $p>2$
it is immediate from the definitions.

We now state a variant of Theorem
\ref{theorem:decreasingdisplacement} for affine actions on Banach
spaces.

\begin{theorem}
\label{theorem:decreasingdisplacementBanach} If $\G$ has property
$(T)$ and $f{\in}\mathcal{U}_2(\G),\delta_0>0,0<C<1$, and
$\eta>0$, there exist $\varepsilon=\varepsilon(\eta)>0$ positive
integers $m=m(f,C,\eta), s=s(f,C,\eta)$ and numbers
$r=r(\delta_0,\eta,f,C)>0$ and $M=M(f,C,\eta)$ such that, letting
$h=f^{*m}$, for any generalized $L^p$ space $\fb$ where
$1+\eta<p{\leq}2$, any $\delta{\leq}\delta_0$, any $x{\in}\fb$,
and any continuous affine $(r, s, {\varepsilon},
{\delta},K)$-action $\rho$ of $K$ on $\fb$ at $x$ we have
\begin{enumerate}
\item $d_{\fh}(x,\rho(h)(x)){\leq}M\dk(x)$;

\item $\dk\rho(h)(x){\leq}C\dk(x)$.
\end{enumerate}
\end{theorem}




\noindent Though they are only concerned with finding fixed points
and do not discuss the iterative method, the special case of
Theorem \ref{theorem:decreasingdisplacementBanach} for (globally
defined) isometric actions is essentially contained in \cite{BG}.
Modulo that fact, the proof of this theorem is quite similar to
the proof of Theorem \ref{theorem:decreasingdisplacementgen}.  In
\cite{BG}, it is also proven that a version of Theorem
\ref{theorem:decreasingdisplacement}  holds for unitary
representations in $L^p$ spaces with $2<p<\infty$ using a simple
duality argument. (Once again they only find fixed points, and do
not describe the iterative method for finding them.) It would be
interesting to know if this is true for representations which are
only ``almost unitary" and $p>2$, but we only need a weaker
statement for our applications, which we now deduce.  We first
define the relevant notion of an ``almost unitary" representation.

\begin{defn}
\label{definition:almostunitary} \begin{enumerate} \item Let
$\sigma$ be a continuous linear representation of $\Gamma$ on a
Banach space $\fb$. Given $\varepsilon>0$, we say that $\sigma$ is
{\em $(K,\varepsilon)$-almost unitary} if for any $k$ in $K$, the
map $\sigma(k)$ is an $\varepsilon$-almost isometry. \item If
$\sigma$ is an $(\infty, s, \varepsilon, 0,0)$-almost action of
$\G$ on a Banach space $\fb$, we call $\sigma$ a {\em
$(K,\varepsilon,s)$-almost unitary representation}.
\end{enumerate}
\end{defn}

\noindent{\bf Remark:} When a fixed choice of $K$ has been made,
we frequently refer to a $(K,\varepsilon)$-almost unitary
representation as an $\varepsilon$-almost unitary representation.

We begin by noting some consequences of Theorem
\ref{theorem:decreasingdisplacementBanach} for a $(K,\varepsilon,s
)$-almost unitary representation $\sigma$ of $\G$ on a generalized
$L^p$ space $\fb$ where $1<p<2$ where $K,\varepsilon,s$ are chosen
to satisfy the conclusions of Theorem
\ref{theorem:decreasingdisplacementBanach} for some values of
$M,C,h$. It is immediate that
$\|\sigma(h)\|{\leq}(1+\varepsilon)^M$ and that
$\|\sigma(h)^n\|{\leq}(1+\varepsilon)^{nM}$.  We can define an
operator $P$ by letting
$Pv=\lim_{n{\rightarrow}{\infty}}\sigma(h)^n(v)$.  It is easy to
see that $\|\sigma(h)^{n+1}-\sigma(h)^n\|{\leq}C^nM$ and therefore
that $\|\sigma(h)^n-P\|{\leq}C^{n-2}M$.  One can then deduce that
$\|P\|<1+\alpha$ where $\alpha$ depends only on $\varepsilon$ and
$\alpha{\rightarrow}0$ as $\varepsilon{\rightarrow}0$.

If we have an $\varepsilon$-almost unitary representation
$\sigma^*$ of $\G$ on a generalized $L^q$ space $\fb$ with
$2<q<\infty$, then the adjoint representation $\sigma$ of
$\sigma^*$ on $\fb^*$ is an $\varepsilon$-almost unitary
representation of $\G$ and $\fb^*$ is a generalized $L^p$ space
for $1<p<2$.  Assuming $f(\g)=f(\g\inv)$ and therefore
$h(\g)=h(\g{\inv})$, it follows that $\sigma(h)^*=\sigma^*(h)$.
Since $\|A\|=\|A^*\|$ for any bounded operator $A$, so the
estimates above carry over for $\sigma^*(h)$, and
$\|\sigma^*(h)\|{\leq}(1+\varepsilon)^M$ and
$\|\sigma^*(h)^n\|{\leq}(1+\varepsilon)^{nM}$.  Furthermore, the
operator $P^*$ defined by letting
$P^*v=\lim_{n{\rightarrow}{\infty}}\sigma^*(h)^n(v)$ is the
adjoint of $P$ and so bounded and a projection. Ideally, $P^*$
would project on $\G$ invariant vectors.  This is easy to verify
if $\varepsilon=0$, but unclear in general.  It is also immediate
that $\|\sigma^*(h)^{n+1}-\sigma^*(h)^n\|{\leq}C^nM$ and that
$\|\sigma^*(h)^n-P^*\|{\leq}C^{n-2}M$.  We summarize this
discussion as follows:

\begin{corollary}
\label{corollary:contractingoperatorbanach} If $\G$ has property
$(T)$ and $f{\in}\mathcal{U}_2(\G)$ satisfies $f(\g)=f(\g{\inv})$
and $0<C<1$ and $1<p_0<\infty$, there exist positive integers
$M=M(f,p_0), s=s(f,p_0)$ and $m=m(f,p_0)$ and
${\varepsilon}=\varepsilon(p_0)>0$ such that, letting $h=f^{*m}$,
for any $p<p_0$ and any $(K,\varepsilon,s)$-almost unitary
representation ${\sigma}$ of $\Gamma$ on a generalized $L^p$ space
$\fb$ and any vector $v$, we have
$d_{\fb}(\sigma(h)^{n+1}v,\sigma(h)^n(v))<MC^n{\dk(v)}$.
Furthermore $Pv=\lim_{n{\rightarrow}{\infty}}\sigma(h)^nv$ is a
bounded linear operator such that
$d_{\fb}(v,Pv){\leq}\frac{MC}{1-C}\dk(v)$.
\end{corollary}

\noindent {\bf Remarks:} \begin{enumerate} \item We emphasize
again that we do not know if $Pv$ is necessarily $\G$ invariant
unless $\sigma$ is unitary. For applications, we will be dealing
with Banach spaces $\fb$ which are $L^p$ type function spaces and
so subspaces of a Hilbert space $\fh$ which is a function space of
type $L^2$. The operator $\sigma(h)$ will be defined on $\fh$ and
we will know, by Theorem \ref{theorem:decreasingdisplacement},
that $\sigma(h)^nv$ converges to a $\G$ invariant vector $v'$.
Corollary \ref{corollary:contractingoperatorbanach} will be used
in conjunction with the Sobolev embedding theorems to obtain
stronger estimates on the regularity of $v'$.  For this argument
to work, it is important to know that we can choose $h$ satisfying
both Corollary \ref{corollary:contractingoperatorbanach} and
Theorem \ref{theorem:decreasingdisplacement} at the same time.  It
is for this reason that we emphasize throughout that $h$ can be
any large enough convolution power of any
$f{\in}\mathcal{U}_2(\G)$. \item We only explicitly use below the
variant of this corollary for $(K,\varepsilon)$-almost unitary
representations.  As remarked above, the version for partially
defined representations can be used in conjunction with Theorem
\ref{theorem:decreasingdisplacementgen} to give a proof of Theorem
\ref{theorem:almostconjugacygen} that does not use Theorem
\ref{theorem:compactlypresented}.\end{enumerate}

\subsection{\bf Foliating Theorem \ref{theorem:isomrigid}}
\label{subsection:foliatedresults}

We now discuss the necessary notions to state our generalization
of Theorem \ref{theorem:isomrigid}.  Though our applications are
to smooth foliations of smooth manifolds, here we work in a
broader setting.

To motivate the results in this section, we state one corollary of
the results of \cite{FM2}, for which Theorem
\ref{theorem:almostconjugacygen} is a key ingredient in the proof.
We call an action $\rho$ of a group $\G$ on $\Ta^n$ {\em linear}
if it is defined by a homomorphism from $\G$ to $GL(n,\Za)$, the
full  group of linear automorphisms of $\Ta^n$.

\begin{corollary}[\cite{FM2}]
\label{corollary:fm2} Let $G$ be a semisimple Lie group with all
simple factors of real rank at least $2$ and let $\G<G$ be a
lattice.  Then any linear action of $\G$ on $\Ta^n$ is
$C^{\infty,\infty}$ locally rigid and there exists a positive
integer $k_0$ depending on the action, such that the action is
$C^{k,k-\frac{n}{2}-2}$ locally rigid for all $k{\geq}k_0$.
\end{corollary}

\noindent This result follows from a more general local rigidity
theorem in \cite{FM2} whose proof uses both Theorem
\ref{theorem:almostconjugacygen} and our results from \cite{FM1}.

Throughout this section $X$ will be a locally compact, second
countable metric space and $\ff$ will be a foliation of $X$ by $n$
dimensional manifolds.  For background on foliated spaces, their
tangent bundles, and transverse invariant measures, the reader is
referred to \cite{CC} or \cite{MS}.  Recall that $\ff$ is a
partition of $X$, satisfying certain additional conditions, into
smooth manifolds called {\em leaves} of the foliation. We will
often refer to the leaf containing $x$ as $\fL_x$.

We let $\Diff^k(X,\ff)$ be the group of homeomorphisms of $X$
which preserve $\ff$ and restrict to $C^k$ diffeomorphisms on each
leaf with derivatives depending continuously on $x$ in $X$.  For
$1{\leq}k{\leq}{\infty}$, there is a natural $C^k$ topology on
$\Diff^k(X,\ff)$. The definition of this topology is
straightforward and we sketch it briefly.  As is usual, the
topology on $\Diff^{\infty}(X,\ff)$ is the inverse limit of the
topologies on $\Diff^k(X,\ff)$ so we now restrict to the case of
$k$ finite. If $X$ is compact, we fix a finite cover of $X$ by
charts $\tilde U_i$ which are products, such that there are proper
subsets $U_i{\subset}\tilde U_i$ which are also products and which
cover $X$.  Without loss of generality, we can identify each $U_i$
as $B(0,r){\times}V_i$ where $B(0,r)$ is standard Euclidean ball
and $V_i$ is an open set in the transversal and identify $\tilde
U_i$ as $B(0,2r){\times}\tilde V_i$ where $\tilde V_i$ is an open
set in the transversal such that $V_i{\varsubsetneq}{\tilde V_i}$.
(See Proposition \ref{proposition:regularatlas} below for a
precise description of such charts.)  A neighborhood of the
identity in $\Diff^k(X,\ff)$ will consist of homeomorphisms $\phi$
which map each $U_i$ inside $\tilde U_i$ and which are uniformly
$C^k$ small as maps from each $B(0,r){\times}\{v\}$ to
$B(0,2r){\times}\{v'\}$, where $v'$ is the point in $\tilde V_i$
such that
$\phi(B(0,r){\times}\{v\}){\subset}B(0,2r){\times}\{v'\}$. When
$X$ is non-compact, there are two possible topologies on
$\Diff^k(X,\ff)$.  The {\em weak topology} is given by taking the
inverse limit of the topologies described above for an increasing
union of compact subsets of $X$.  To define the {\em strong
topology}, we cover $X$ by a countable collection of neighborhoods
$U_i{\subset}\tilde U_i$ as described above, and take the same
topology. When $X$ is not compact, we will always consider the
strong topology.  Though non-compact foliated spaces arise in the
proofs, for the remainder of this subsection, we consider only
compact $X$.

We now define the type of perturbations of actions that we will
consider.

\begin{defn}
\label{definition:foliatedperturbation} Let $\G$ be a compactly
generated topological group and $\rho$ an action of $\G$ on $X$
defined by a homomorphism from $\G$ to $\Diff^{\infty}(X,\ff)$.
Let $\rho'$ be another action of $\G$ on $X$ defined by a
homomorphism form $\G$ to $\Diff^k(X,\ff)$. Let $U$ be a (small)
neighborhood of the identity in $\Diff^k(X,\ff)$ and $K$ be a
compact generating set for $\G$. We call $\rho'$ a {\em
$(U,C^k)$-foliated perturbation} of $\rho$ if:
\begin{enumerate}
\item  for every leaf $\fL$ of $\ff$ and every $\g{\in}\G$, we
have $\rho(\g)\fL=\rho'(\g)\fL$ and,

\item $\rho'(\g)\rho(\g){\inv}$ is in $U$ for every $\g$ in $K$.

\end{enumerate}
\end{defn}

We fix a continuous, leafwise smooth Riemannian metric $g_{\ff}$
on $T{\ff}$, the tangent bundle to the foliation and note that
$g_{\ff}$ defines a volume form and corresponding measure on each
leaf $\fL$ of $\ff$, both of which we denote by $\nu_{\ff}$.
(Metrics $g_{\ff}$ exist by a standard partition of unity
argument.)  Let $\G$ be a group and $\rho$ an action of $\G$ on
$X$ defined by a homomorphism from $\G$ to $\Diff^k(X,\ff)$. We
say the action is {\em leafwise isometric} if $g_{\ff}$ is
invariant under the action.  When $\G=\Za$ and $\Za=<f>$, we will
call $f$ a {\em leafwise isometry}.

For the remainder of the paper, we will assume that the foliation
has a transverse invariant measure $\nu$. By integrating the
transverse invariant measure $\nu$ against the Riemannian measure
on the leaves of $\ff$, we obtain a measure $\mu$ on $X$ which is
finite when $X$ is compact. In this case, we normalize $g_{\ff}$
so that $\mu(X)=1$. We will write $(X,\ff,g_{\ff},\mu)$ for our
space equipped with the above data, sometime leaving one or more
of $\ff, g_{\ff}$ and $\mu$ implicit.  We will refer to the
subgroup of $\Diff^k(X,\ff)$ which preserves $\nu$ as
$\Diff_{\nu}^k(X,\ff)$.  Note that if $\rho$ is an action of $\G$
on $X$ defined by a homomorphism into $\Diff_{\nu}^k(X,\ff)$ and
$\rho$ is leafwise isometric, then $\rho$ preserves $\mu$.
Furthermore if $\rho$ is an action of $\G$ on $X$ defined by a
homomorphism into $\Diff_{\nu}^k(X,\ff)$ and $\rho'$ is a
$(U,C^k)$-leafwise perturbation of $\rho$, then it follows easily
from the definition that $\rho'$ is defined by a homomorphism into
$\Diff^k_{\nu}(X,\ff)$ since the induced map on transversals is
the same.

The following foliated version of Theorem \ref{theorem:isomrigid}
is one of the key steps in the proof of the main results in
\cite{FM2}. We denote by $B_{\ff}(x,r)$ the ball in $\fL_x$ about
$x$ of radius $r$.  For a sufficiently small value of $r>0$, we
can canonically identify each $B_{\ff}(x,2r)$ with the ball of
radius $2r$ in Euclidean space via the exponential map from
$T\ff_x$ to $\fL_x$.  To state our results, we will need a
quantitative measure of the size of the $k$-jet of $C^k$ maps. We
first consider the case when $k$ is an integer, where we can give
a pointwise measure of size.  Recall that a $C^k$ self map of a
manifold $Z$ acts on $k$-jets of $C^k$ functions on $Z$. Any
metric on $TZ$ defines a pointwise norm on each fiber of the
bundle of $J^k(Z)$ of $k$-jets of functions on $Z$. For any $C^k$
diffeomorphism $f$ we can define $\|j^k(f)(z)\|$ as the operator
norm of the map induced by $f$ from $J^k(Z)_z$ to $J^k(Z)_{f(z)}$.
For a more detailed discussion on jets and an explicit
construction of the norm on $J^k(Z)_z$, see section
\ref{section:parametrizing}.  We say that a map $f$ has $C^k$ size
less than $\delta$ on a set $U$ if $\|j^k(f)(z)\|<\delta$ for all
$z$ in $U$. If $k$ is not an integer, we say that $f$ has $C^k$
size less than $\delta$ on $U$ if $f$ has $C^{k'}$ size less than
$\delta$ on $U$ where $k'$ is the greatest integer less than $k$
and $j^{k'}(f)$ satisfies a (local) H\"older estimate on $U$. See
section \ref{section:parametrizing} for a more detailed discussion
of H\"older estimates.

\noindent{\bf Remark:} This notion of $C^k$ size is not very
sharp.  The size of the identity map will be $1$, as will be the
size of any leafwise isometry.  We only use this notion of size to
control estimates on a map at points where the map is known to be
``fairly large" and where we only want bounds to show it is ``not
too large".

 For the following theorem, we assume
that the holonomy groupoid of $(X,\ff)$ is Hausdorff. This is a
standard technical assumption that allows us to define certain
function spaces on ``pairs of points on the same leaf of
$(X,\ff)$".  See subsection \ref{subsection:functionsonpairs},
\cite{CC} and \cite{MS} for further discussion.  All the
foliations considered in \cite{FM2} are covered by fiber bundles,
 in which case it is easy to show that the holonomy groupoid is
Hausdorff.


\begin{theorem}
\label{theorem:almostconjugacygen} Let $\Gamma$ be a locally
compact, $\sigma$-compact group with property $(T)$.  Let $\rho$
be a continuous leafwise isometric action of $\Gamma$ on $X$
defined by a homomorphism from $\G$ to
$\Diff_{\nu}^{\infty}(X,\ff)$. Then for any $k{\geq}3,\kappa>0$
and any $\varsigma>0$ there exists a neighborhood $U$ of the
identity in $\Diff^k(X,\ff)$ such that for any continuous
$(U,C^k)$-foliated perturbation $\rho'$ of $\rho$ there exists a
measurable $\Gamma$-equivariant map $\phi:X{\rightarrow}X$ such
that:

\begin{enumerate}
\item $\phi{\circ}\rho(\g)=\rho'(\g){\circ}\phi$ for all
$\g{\in}\G$,

\item $\phi$ maps each leaf of $\ff$ into itself,


\item there is a subset $S{\subset}X$ with $\mu(S)=1-\varsigma$
and $\G{\cdot}S$ has full measure in $X$, and a constant
$r{\in}\Ra^+$, depending only on $X,\ff$ and $g_{\ff}$, such that,
for every $x{\in}S$, the map $\phi:B_{\ff}(x,r){\rightarrow}\fL_x$
is $C^{k-1-\kappa}$-close to the identity; more precisely, with
our chosen identification of $B_{\ff}(x,2r)$ with the ball or
radius $2r$ in Euclidean space,
$\phi-\Id:B_{\ff}(x,r){\rightarrow}B_{\ff}(x,2r)$ has
$C^{k-1-\kappa}$ norm less than $\varsigma$ for every $x{\in}S$,
and

\item there exists $0<t<1$ depending only on $\G$ and $K$ such
that the set of $x{\in}X$ where the $C^{k-1-\kappa}$ size of
$\phi$ on $B_{\ff}(x,r)$ is not less than $(1+\varsigma)^{l+1}$
has measure less than $t^l\varsigma$ for any positive integer $l$.

\end{enumerate}
Furthermore, for any $l{\geq}k$, if $\rho'$ is a $C^{2l-k+1}$
action, then by choosing $U$ small enough, we can choose $\phi$ to
be $C^{l}$ on $B_{\ff}(x,r)$ for almost every $x$ in $X$.  In
particular, if $\rho'$ is $C^{\infty}$ then for any $l{\geq}k$, by
choosing $U$ small enough, we can choose $\phi$ to be $C^l$ on
$B_{\ff}(x,r)$ for almost every $x$ in $X$.
\end{theorem}

\noindent {\bf Remarks:}\begin{enumerate} \item The map $\phi$
constructed in the theorem is not even $C^0$ close to the identity
on $X$. However, the proof of the theorem shows that for every
$1{\leq}q{<}{\infty}$, possibly after changing $U$ depending on
$q$, we have $\int_X(d(x,\phi(x))^qd\mu{\leq}\varsigma$.



\item This theorem implies a version of Theorem
\ref{theorem:isomrigid}, but with lower regularity. \item In some
special cases it is possible to slightly improve the regularity of
$\phi$.  It is possible to show that $\phi$ is $C^l$ for some
given choice of $l$ even if $\rho'$ is only  $C^{l+1}$ provided
$U$ is small enough, see section
\ref{section:convexityofderivatives} for more discussion.  Unlike
in Theorem \ref{theorem:isomrigid}, it does not seem possible to
show that $\phi$ is $C^{\infty}$ without some assumption on the
action transverse to $\ff$.  Again see section
\ref{section:convexityofderivatives} for more details.
\end{enumerate}

In the case when $X$ is a direct product, we can prove slightly
greater regularity.

\begin{theorem}
\label{theorem:almostconjugacyproduct} If $X=Y{\times}Z$ and the
foliation $\ff$ has leaves of the form
$\{\{y\}{\times}Z|y{\in}Y\}$, then $\phi$ in Theorem
\ref{theorem:almostconjugacygen} is $C^{k-\kappa}$ and all
estimates in that theorem for the $C^{k-1-\kappa}$ topology can be
replaced by analogous estimates in the $C^{k-\kappa}$ topology.
\end{theorem}

\noindent {\bf Remark:} We do not give a proof of Theorem
\ref{theorem:almostconjugacyproduct} here.  The proof of Theorem
\ref{theorem:isomrigid} given in subsection
\ref{subsection:isomrigid} can be combined with the techniques of
section \ref{section:foliatedproofs} to give such a proof, which
we leave as an exercise for the interested reader.

\section{\bf Proof of Theorem
\ref{theorem:decreasingdisplacementgen} and variants}
\label{section:Tproofs}

In this section we prove Theorem
\ref{theorem:decreasingdisplacementgen}.   In the first
subsection, we give a proof of the analogue of Theorem
\ref{theorem:decreasingdisplacementgen} for isometric actions of
groups with property $(T)$ on Hilbert spaces.  In the second
subsection we develop a general method of constructing limit
actions from sequences of actions. In the third subsection, we
prove Theorem \ref{theorem:decreasingdisplacementgen} modulo some
observations contained in the appendix to this paper, which are
required only when the action is not affine and $\G$ is not
discrete. In the final subsection we recall some results from
\cite{BG} and some facts about Banach spaces of type $L^p_n$ and
indicate the modifications to prior arguments needed to prove the
results in subsection \ref{subsection:TresultsBanach}.

\subsection{\bf Finding fixed points for isometric actions of groups
with Property $(T)$} \label{subsection:fpstandard}

Theorem \ref{theorem:decreasingdisplacementgen} is a
generalization of the following consequence of property $(T)$.
Though this fact is a variant of well-known consequences of
property $(T)$, we did not find a prior reference for this precise
statement and so give a detailed proof.

We first fix some notation.  As in subsection
\ref{subsection:Tresults} we fix a locally compact,
$\sigma$-compact group $\Gamma$ with a compact generating set
$K{\subset}\Gamma$ containing a neighborhood of the identity, and
a (left) Haar measure $\mu$ on $\G$. Given a function
$h{\in}C_c(\Gamma)$ and $\gamma_0{\in}\Gamma$ we write
$\g_0{\cdot}h$ for the function $\g{\rightarrow}h(\g_0{\inv}\g)$.
The subsets $\cu(\G)$ and $\cu_2(\G)$ of $C_c(\G)$ are as in
subsection \ref{subsection:Tresults}

\begin{proposition}
\label{proposition:decreasingdisplacementstandard} If $\G$ has
property $(T)$ and $f{\in}\mathcal{U}_2(\G)$ and $0<C_0<1$ there
exist positive integers $M=M(f,C_0)$ and $m=m(f,C_0)$, such that,
letting $h=f^{*m}$,  for any Hilbert space $\fh$, any continuous
isometric action $\rho$ of $\Gamma$ on $\fh$, and any
$x{\in}{\fh}$ we have
\begin{enumerate}
\item $d_{\fh}(x,\rho(h)(x)){\leq}M\dk(x)$

\item $\dk(\rho(h)(x)){\leq}C_0\dk(x)$.
\end{enumerate}
\end{proposition}

Given a unitary representation $\sigma$ of $\G$ on $\fh$, we let
$\fh_{\sigma}$ be the $\sigma$ invariant vectors and
$\fh_{\sigma}^{\perp}$ it's orthogonal complement. We recall a
fact about groups with property $(T)$.

\begin{lemma}
\label{lemma:contraction} Let $\fh$ be a Hilbert space and
$\sigma$ a continuous unitary representation of $\Gamma$ on $\fh$.
Then for any  $f{\in}\cu_2(\G)$, we have
$\sigma(f)|_{\fh_{\sigma}^{\perp}}$ is a contraction. More
precisely, there exists a constant $0<D<1$ such that
$\|\sigma(f)(x)\|<D\|x\|$ for any $x{\in}{\fh_{\sigma}^{\perp}}$.
\end{lemma}

\noindent This lemma is an immediate consequence of Kazhdan's
definition of property $(T)$ and the characterization of the Fell
topology in Lemma III.1.1 of \cite{M}. Though explicitly stated
there only for some $f$, the proof is valid for any
$f{\in}\mathcal{U}_2(\G)$. For a proof of a more general fact see
Lemma \ref{lemma:almostinvariantequalscontracting} below.  The
following lemma is elementary from the fact that isometries of
Hilbert spaces are affine \cite{MU}.

\begin{lemma}
\label{lemma:convolutionisometrichilbert} Let $\G$ be a group,
$\fh$ a Hilbert space and $\rho$ an isometric $\G$ action on
$\fh$. Then for any measures $\mu,\lambda{\in}\mathcal{P}(\G)$, we
have $\rho(\mu)\rho(\lambda)=\rho(\mu*\lambda)$.
\end{lemma}

\begin{lemma}
\label{proposition:decreasingdisplacementrep}  If $\G$ has
property $(T)$ and $f{\in}\mathcal{U}(\G)$ and $0<C_0<1$,there
exist positive integers $M=M(f,C_0)$ and $m=m(f,C_0)$ such that,
letting $h=f^{*m}$, for any Hilbert space $\fh$, any continuous
unitary representation $\sigma$ of $\Gamma$ on $\fh$, and any
$v{\in}{\fh}$ we have
\begin{enumerate}
\item $d_{\fh}(v,\sigma(h)(v)){\leq}M\dk(v)$

\item $\dk(\sigma(h)(v)){\leq}C_0\dk(v)$.
\end{enumerate}
\end{lemma}

\begin{proof}
Since if $v=(v_1,v_2)$, where $v_1{\in}\fh_{\sigma}$ and
$v_2{\in}\fh_{\sigma}^{\perp}$,  $\dk(v)=\dk(v_2)$, it suffices to
assume $\fh=\fh_{\sigma}^{\perp}$. Since $\Gamma$ has property
$(T)$ this implies that there exists $\varepsilon$ such that there
are no $(K,\varepsilon)$-invariant vectors in $\fh$, i.e.
$$\varepsilon\|v\|<\|\sigma(k)v-v\|{\leq}2\|v\|$$
for any $k{\in}K$ and any $v{\in}\fh$.   Let $D$ be the
contraction factor from \ref{lemma:contraction} and choose $m$
such that ${2}{D^m}{\leq}C_0\varepsilon$ and let $h=f^{*m}$. Note
that Lemma \ref{lemma:convolutionisometrichilbert} implies that
$\sigma(h)=\sigma(f^{*m})=\sigma(f)^m$.  Let $v$ be a vector with
$K$-displacement $\delta$. It follows from the equation above that
$\delta>{\varepsilon}\|v\|$. Direct computation shows that
$\|\sigma(k)\sigma(h)(v)-\sigma(h)(v)\|<2\|\sigma(h)v\|<
2D^m\|v\|{\leq}C_0\varepsilon\|v\|<C_0\delta$ which is the second
conclusion of the lemma. Letting $M$ be the smallest value such
that $\supp(h){\subset}K^M$, the first conclusion follows as well.
\end{proof}

\begin{proof}[Proof of Proposition
\ref{proposition:decreasingdisplacementstandard}]

Fix the function $h{\in}\cu(\G)$ and the constant $0<C_0<1$ from
the conclusion of Lemma
\ref{proposition:decreasingdisplacementrep}. As any continuous
affine, isometric action of a group with property $(T)$ on a
Hilbert space has a fixed point, $\G$ fixes some point $x$ in
$\fh$ \cite{De}. Viewing $\fh$ as a vector space with $x$ as
origin allows us to view our action $\rho$ as a unitary
representation, and the proposition is now an immediate
consequence of Proposition
\ref{proposition:decreasingdisplacementrep} with the same $h,M$
and $C_0$.

\end{proof}

\subsection{\bf Limits of sequences of actions}
\label{subsection:limitactions} In this subsection we give a very
general process for constructing a limit action  from a sequence
of actions, or partially defined actions.  The reader primarily
interested in actions of discrete groups may compare this with the
discussion of scaling limits in \cite{Gr2} and the references
cited there. Throughout this subsection $\G$ is a group and $K$ is
a generating set for $\G$.

Let $X_n$ be a sequence of complete metric spaces, with
distinguished points $x_n{\in}X_n$, and let $\rho_n$ be
$(r_n,s_n,\varepsilon_n,\delta_n,K)$-almost actions of a group
$\G$ on $X_n$ at $x_n$.  We construct our limit space $X$ as a
quotient of a certain subspace ${\tilde X}$ in $\prod{X_n}$.  We
will use ultrafilters and ultralimits to define ${\tilde X}$, and
a pseudo-metric on $\tilde X$, and then let $X$ be $\tilde X$
modulo relation of being at distance zero in the pseudo-metric.

\begin{defn}
\label{definition:ultrafilter} A {\em non-principal ultrafilter}
is a finitely additive probability measure $\omega$ defined on all
subsets of $\mathbb N$ such that
\begin{enumerate}
\item $\omega(S)=0$ or $1$, \item $\omega(S)=0$ if $S$ is finite.
\end{enumerate}
\end{defn}

This definition is the one given by Gromov in \cite{Gr1}, at the
beginning of section $2.A.$ on page $36$. It is not clear with
this definition that non-principal ultrafilters exist.  To show
existence, one defines an {\em ultrafilter} as a maximal {\em
filter}, and shows that maximal objects exist using Zorn's lemma.
For a more traditional definition of ultrafilters, see
\cite[I.6.4]{BTG}.  In the context of group theory ultrafilters
were first used to construct limits of sequences of metric spaces
in \cite{vDW}, though their use in the study of Banach space
theory is much older than that, see \cite{BL} and \cite{He} for
more history. In what follows, we fix a non-principal ultrafilter
$\omega$.

Let $\{y_n\}$ be sequence in $\mathbb R$, the {\em $\omega$-limit}
of $\{y_n\}$ is $\omega$-$\lim{y_n}=y$ if for every $\epsilon>0$
it follows that $\omega\{n|d(y_n,y)<\epsilon\}=1$.  The following
well-known proposition can be proven easily by mimicking the proof
that bounded sequences have limit points.

\begin{proposition}
\label{proposition:ultralimit} Any bounded sequence of real
numbers has a unique $\omega$-limit.
\end{proposition}

More generally, if $X$ is a Hausdorff topological space, and
$\{y_n\}$ is a sequence of points in $X$, the {\em $\omega$-limit}
of $\{y_n\}$ is $\omega$-$\lim{y_n}=y$ if every neighborhood of
$y$ has measure $1$ with respect to the pushforward of $\omega$
under the map $n{\rightarrow}y_n$.  The following, almost
tautological proposition, is from \cite[I.10.1]{BTG}:

\begin{proposition}
\label{proposition:ultralimitgeneral} The space $X$ is compact if
and only if, for every ultrafilter $\omega$, every sequence
$\{y_n\}$ has a unique $\omega$-limit.
\end{proposition}

 We
let $\tilde X= \{y{\in}\prod{X_n}| \omega$-$\lim
d_n(y_n,x_n)<\infty\}$. We put a metric on $\tilde X$ by letting
${\tilde d}(\{v_n\},\{w_n\})=\omega$-$\lim{d_{n}(v_n,w_n)}$. It is
easy to check that $\tilde d$ is a pseudo-metric on $\tilde X$. We
can define an equivalence relation on $\tilde X$ by letting
$v{\sim}w$ if ${\tilde d}(v,w)=0$. We let $X={\tilde X}/{\sim}$
with the metric $d$ defined by ${\tilde d}$. For an arbitrary
sequence $y{\in}\tilde X$, we refer to the image of $y=\{y_n\}$ in
$X$ as $y_{\omega}$ and write $y_{\omega}=\omega$-$\lim{y_n}$. The
space $X$ has a natural basepoint given by
$x_{\omega}=\omega$-$\lim{x_n}$. The space $(X,d)$ is often called
{\em the $\omega$-ultraproduct}, or simply {\em the ultraproduct},
with $\omega$ implicit, of $(X_i,d_i,x_i)$. The following
straightforward proposition is standard.

\begin{proposition}
\label{proposition:complete} The space $(X,d)$ is complete.
\end{proposition}

\begin{proof}
 Let $x_{\omega}^j$ be a Cauchy sequence in $X$, where
$x_{\omega}^j=\omega$-$\lim{x_n^j}$.  Let $N_1=\Na$. Inductively,
there is an $\omega$-full measure subset $N_j{\subseteq}N_{j-1}$
such that $n{\in}N_j$ implies that
$|d_n(x^k_n,x^l_n)-d(x^k_{\omega},x^l_{\omega})|{\leq}\frac{1}{2^n}$
for $1{\leq}k,l{\leq}j$.  For $n{\in}N_j{\setminus}N_{j-1}$,
define $y_n=x_n^j$.  Then $x^j_{\omega}$ converge to $y_{\omega}$.
\end{proof}

We record here one additional fact about limits of sequences of
Hilbert spaces that we will use in the proof of Theorem
\ref{theorem:decreasingdisplacementgen}, compare \cite{He}.

\begin{proposition}
\label{proposition:ishilbert} If the spaces $X_n$ are Hilbert
spaces with $x_n=0$ and inner product $<,>_n$, then the space
$(X,d)$ is Hilbert space with $x_{\omega}=0$ and inner product
defined by $<v_{\omega},w_{\omega}>=\omega$-$\lim<v_n,w_n>$.
\end{proposition}

\begin{proof} Since we already know $X$ is complete, we need only
check that $X$ is a vector space and that $<\cdot,\cdot>$ is a
positive definite symmetric bilinear form.  Letting
$V=\{v{\in}\tilde X|d(v,\{x_n\})=0\}$, it is immediate that $V$ is
a sub-vector space of the vector space ${\tilde X}$ and that
$\tilde X/V=\tilde X/{\sim}$.

That $<,>$ is symmetric and bilinear is immediate.   The
definition implies that $<v\so,v\so>=d(v\so,x\so)^2>0$ if
$v\so{\neq}x\so$, so the form is positive definite.
\end{proof}

We now proceed to define a $\Gamma$ action $\rho$ on $X$.  If
$\delta_n$ and $\varepsilon_n$ were bounded sequences and each
$\rho_n$ were globally defined, we could define a $\Gamma$ action
$\tilde \rho$ on $\tilde X$ simply by acting on each coordinate.
Instead we define $\tilde \rho(\gamma)(y)$ to be the sequence
whose $n$th coordinate is $\rho_n(\gamma)(y_n)$ when
$\rho_n(\gamma)(y_n)$ is defined and whose $n$th coordinate is
$x_n$ otherwise.  Though this is not an action,  we have:

\begin{proposition}
\label{proposition:limitaction} If
$\omega$-$\lim\varepsilon_n=\varepsilon<\infty,
\omega$-$\lim\delta_n=\delta<\infty$ and
$\omega$-$\lim_{n{\rightarrow}\infty}r_n=
\omega$-$\lim_{n{\rightarrow}\infty}s_n=\infty$, then for every
$\gamma$ in $\Gamma$, the map $\tilde \rho(\gamma)$ descends to a
well-defined bilipschitz map $\rho(\gamma)$ of $X$ and the map
$\rho:\Gamma{\times}X{\rightarrow}X$ is an action of $\G$ on $X$.
Furthermore, $\rho(k)$ is an $\varepsilon$-almost isometry of $X$
for every $k$ in $K$ and the $K$-displacement of $x_{\omega}$ is
at most $\delta$.
\end{proposition}

\noindent {\bf Remark:} For our applications, we will have
$\delta_n$ uniformly bounded,
$\lim_{n{\rightarrow}\infty}\varepsilon_n=0$,
$\lim_{n{\rightarrow}\infty}r_n=\infty$ and
$\lim_{n{\rightarrow}\infty}s_n=\infty$.

\begin{proof}
To verify that $\tilde \rho(\gamma)$ descends to $X$ and that it
is bilipschitz, it suffices to verify that $\tilde \rho(\gamma)$
is an almost isometry of the pseudo-metric $\tilde d$. Let
$u,v{\in}{\tilde X}$ and $\gamma{\in}\Gamma$. We fix the minimal
$s$ such that $\gamma{\in}K^s$.  By ignoring an $\omega$-measure
zero, finite set of indices, we may assume $\rho_n(\gamma)u_n$ and
$\rho_n(\gamma)v_n$ are defined. By definition

$${\tilde d}({\tilde \rho}(\gamma)u,{\tilde
\rho}(\gamma)v)=\omega{\text
-}\lim{d_n}(\rho_n(\g)u_n,\rho_n(\g)v_n).$$

\noindent Since $\rho_n(k)$ acts by $\varepsilon_n$-almost
isometries for all $k{\in}K$, we have

\noindent
$$(1-\varepsilon_n)^sd_n(u_n,v_n){\leq}d_n(\rho_n(\g)u_n,\rho_n(\g)v_n){\leq}
(1+\varepsilon_n)^sd_n(u_n,v_n).$$

\noindent By taking the $\omega$-limit of the above equation, we
have

$$(1-\varepsilon)^s{\tilde d}(u,v){\leq}{\tilde d}({\tilde \rho}(\g)u,{\tilde \rho}(\g)v){\leq}
(1+\varepsilon)^s{\tilde d}(u_n,v_n).$$

\noindent Which shows that $\tilde \rho(\g)$ preserve the
equivalence relation of being at $\tilde d$ distance zero as well
as showing that the map $\rho(\g)$ on the quotient $X$ is
bilipschitz and in fact an $\varepsilon$-almost isometry when
$\g{\in}K$.

That these maps form a $\Gamma$ action is almost obvious.  Fix
$\gamma_1,\gamma_2{\in}\Gamma$ and $v{\in}{\tilde X}$.  By
ignoring a finite set $S_{\g_1,\g_2,v}$ of $\omega$-measure zero,
we can insure that $\rho_n(\g_1\g_2)v_n, \rho_n(\g_2)v_n$ and
$\rho_n(\g_1)(\rho_n(\g_2)v_n)$ are well-defined and that
$\rho_n(\g_1\g_2)v_n= \rho_n(\g_1)(\rho_n(\g_2)v_n)$ for
$n{\notin}S_{\g_1,\g_2,v}$. This implies
$\rho(\g_1\g_2){v\so}=\rho(\g_1)(\rho(\g_2)v\so)$. Since this
verification (though not the set $S_{\g_1,\g_2,v}$) is independent
of $\g_1,\g_2$ and $v$, it follows that $\rho$ is an action. That
the $K$-displacement of $x_{\omega}$ is less than $\delta$ follows
since
$d(\rho(k)x_{\omega},x_{\omega})=\omega$-${\lim}d_n(\rho_n(k)x_n,x_n){\leq}\delta$
for all $k{\in}K$.
\end{proof}

{\noindent}{\bf Remark:} As shorthand for the construction above,
we will write $\rho=\omega$-$\lim \rho_n$.

\subsection{Proof of Theorem
\ref{theorem:decreasingdisplacementgen}}
\label{subsection:Tproofs}

In this subsection $\G$ and $K$ are as in subsection
\ref{subsection:fpstandard}. We fix the function $h{\in}\cu(\G)$
and the constant $C_0$ given by Proposition
\ref{proposition:decreasingdisplacementstandard}. We also fix an
arbitrary $C$ with $C_0<C<1$.

\begin{proof}[Proof of Theorem
\ref{theorem:decreasingdisplacementgen} for $\G$ discrete.] Fix
$\eta_0=C-C_0$.  The proof proceeds by contradiction, so we assume
the theorem is false. Let $r_n=2^n, s_n=n$ and
$\varepsilon_n=\frac{1}{n}$ and $0<\delta_n<\delta_0$. By the
assumption that Theorem \ref{theorem:decreasingdisplacementgen} is
false there exists a sequence of Hilbert spaces $\fh_n$, points
$x_n{\in}\fh_n$ and $(r_n, s_n, \varepsilon_n, \delta_n,K)$-almost
actions $\rho_n$ at $x_n$ such that $\dk(\rho_n(h)(x))>C\delta_n$.
By conjugation by a homothety at $x_n$, it suffices to consider
the case where $\delta_n=1$ for all $n$. Conjugating by this
homothety makes $\rho_n$ a
$(r_n\frac{1}{\delta_n},s_n,\varepsilon_n,1,K)$-almost action at
$x_n$ and it remains true that $r_n{\rightarrow}{\infty}$ as
$n{\rightarrow}{\infty}$.  We will denote the distance on $\fh_n$
as $d_n$ and the inner product as $<{\cdot},\cdot>_n$.  Letting
$\tilde \fh{\subset}\prod\fh_n$ be as in the paragraph following
Proposition \ref{proposition:ultralimitgeneral} and $\tilde
d=\omega$-$\lim{d_n}$ and $V$ the set of points in $\tilde \fh$
with $\tilde d(v,\{x_n\})=0$ as above, it follows from Proposition
\ref{proposition:limitaction}, the fact that
$\lim_{n{\rightarrow}{\infty}}\varepsilon_n=0$, and Proposition
\ref{proposition:ishilbert} that the action
$\rho=\omega$-$\lim{\rho_n}$ of $\G$ on $\fh={\tilde \fh}/V$ is an
isometric action on a Hilbert space.  It also follows from
Proposition \ref{proposition:limitaction} that
$\dk(x_{\omega}){\leq}1$.   By Proposition
\ref{proposition:decreasingdisplacementstandard}, this implies
that $\dk(\rho(h)(x\so)){\leq}C_0$. It is immediate from the
definitions that $\rho(h)(x\so)=\{\rho_n(h)x_n\}\so$ and that
$d(\rho(k)y\so,y\so)=\omega$-$\lim{d_n}(\rho_n(k)y,y)$ for any
$k{\in}K$ and $y{\in}\tilde \fh$.  Letting $y=\rho(h)(x\so)$, we
have a set $S_k$ of full $\omega$-measure such that
$d_n(\rho_n(k)\rho_n(h)(x\so),\rho_n(h)(x\so)){\leq}C_0+\eta_0$
for all $n$ in $S_k$. Letting $S={\cap}_{k{\in}K}S_k$ we see that
$\dk(\rho_n(h)(x_n)){\leq}C=C_0+\eta_0$ for any $n{\in}S$. Since
$K$ is finite, $\omega(S)=1$, and we have a contradiction.
\end{proof}

Before proving the theorem for more general groups, we state some
additional results needed because the limit action we construct is
not necessarily continuous.

For the remainder of this subsection, we assume that each $\rho_n$
is continuous and that $\rho_n$ satisfy the conditions of
Proposition \ref{proposition:limitaction}.  We can then define a
limit action $\rho=\omega$-$\lim\rho_n$ as in Proposition
\ref{proposition:limitaction}.  In general it is not true that
$\rho$ is continuous, but we now describe a (possibly trivial)
continuous subaction of $\rho$.

Given $y_n{\in}X_n$, we have an orbit map
$\rho_n^{y_n}:K^{s^n}{\rightarrow}X_n$ defined by
$\rho^{y_n}_n(\g)=\rho_n(\g)(y_n)$. We call a sequence $\{y_n\}$
{\em $\omega$-equicontinuous on compact sets} if for any compact
subset $D$ of $\G$, there exists a subset $S{\subset}\Na$ with
$\omega(S)=1$ such that the orbit maps $\rho_n^{y_n}$ are
equicontinuous on $D$ for $n$ in $S$. Since the collection of
actions $\rho_n$ are uniformly bilipschitz, to prove a sequence is
$\omega$-equicontinuous on compact sets, it suffices to prove that
it is $\omega$-equicontinuous at the identity in $\G$, i.e. given
$\varepsilon>0$, there is a neighborhood $U$ of the identity in
$\G$ such that $\rho_n^{y_n}(U){\subset}B(y_n,\varepsilon)$ for
$n$ in a set $S$ with $\omega(S)=1$. We denote by $\Omega$ the set
of $\omega$-equicontinuous sequences in $\tilde X$ and by $\bar
\Omega$ the image of $\Omega$ in $X$. Keeping in mind that $\rho$
is an action by bilipschitz maps, it is straightforward to verify
the following.

\begin{proposition}
\label{proposition:continuoussubaction} The set $\bar \Omega$ is
closed and $\G$ invariant.  The restriction of $\rho$ to $\bar
\Omega$ is continuous.
\end{proposition}

\noindent We state a result giving sufficient conditions for $\bar
\Omega$ to be an affine Hilbert subspace when the $X_n$ are all
Hilbert spaces.

\begin{proposition}
\label{proposition:continuoussubactionHilbert} Let $\fh_n$ be a
sequence of Hilbert spaces with basepoints $x_n$. Let $\rho_n$ be
a sequence of continuous
$(r_n,s_n,\varepsilon_n,\delta_n,K)$-almost actions of $\G$ on
$X_n$ at $x_n$ with $\omega$-$\lim\varepsilon_n=0,
\omega$-$\lim\delta_n=\delta<\infty$ and
$\omega$-$\lim_{n{\rightarrow}\infty}r_n=
\omega$-$\lim_{n{\rightarrow}\infty}s_n=\infty$.  Then the set
$\bar \Omega{\subset}\fh$ is an affine Hilbert subspace of $\fh$.
Furthermore if $f{\in}\mathcal{U}(\G)$, then for any sequence
$\{y_n\}{\in}{\tilde \fh}$, the point $\omega$-$\lim\rho_n(f)y_n$
is in $\bar \Omega$.
\end{proposition}

\noindent We will also need the following generalization of Lemma
\ref{lemma:convolutionisometrichilbert} for almost isometric
actions on Hilbert spaces.

\begin{proposition}
\label{proposition:convolutionalmostisometrichilbert} Let $\fh$ be
a Hilbert space.  Then for every $r,\eta>0$ there is a
$\varepsilon>0$ such that for any continuous
$(r,s,\varepsilon,\delta,K)$-almost action of $\G$ on $\fh$ at $x$
and any measures $\mu,\lambda{\in}\mathcal{P}(\G)$ with
$\supp(\mu),\supp(\lambda)$ and $\supp(\mu*\lambda)$ contained in
$K^s$ and $s\delta<\frac{r}{2}$, we have
$$d(\rho(\mu)\rho(\lambda)x,\rho(\mu*\lambda)x){\leq}{\eta}.$$
\end{proposition}

\noindent {\bf Remark:} For our applications to local rigidity, it
suffices to prove Theorem \ref{theorem:decreasingdisplacementgen}
for affine $\varepsilon$-almost isometric actions.  In this case,
the proof of Theorem \ref{theorem:decreasingdisplacementgen} is
almost the same, but we can assume $\rho_n$ is affine for all $n$.
We will therefore only prove those cases of Propositions
\ref{proposition:continuoussubactionHilbert} and
\ref{proposition:convolutionalmostisometrichilbert} here.
 Only readers interested in Theorem
\ref{theorem:decreasingdisplacementgen} for the case of $\rho$ not
affine and $\G$ not discrete, need refer to Appendix
\ref{appendix:limits} of this paper for proofs of the general
cases of Propositions \ref{proposition:continuoussubactionHilbert}
and \ref{proposition:convolutionalmostisometrichilbert}.

\begin{proof}[Proof of Propositions \ref{proposition:continuoussubactionHilbert} and
\ref{proposition:convolutionalmostisometrichilbert} for affine
actions] It is immediate that $\bar \Omega$ is an affine Hilbert
subspace and that
$\rho_n(\mu)\rho_n(\lambda)=\rho_n(\mu*\lambda)$.  We now prove
that $\rho_n(f)y_n$ is $\omega$ equicontinuous. To do so we use
the following estimate:
$$d(\rho_n(\g_0)\rho_n(f)y_n,\rho_n(f)y_n){\leq}$$
$$\|\rho(\g_0{\cdot}f-f)y_n\|{\leq}$$
$$\|\g_0{\cdot}f-f\|_{L^1}D_{\g_0,f}$$
where $D_{\g_0,f}=\sup_{\supp(\g_0{\cdot}f-f)}d(\rho_n(\g)x,x).$
This estimate, our assumptions on $\rho_n$, the fact that $K$
contains a neighborhood of the identity in $\G$, and continuity of
the $\G$ action on $L^1(\G)$ imply that for any $\eta>0$ there is
a neighborhood $U$ of the identity in $\G$ such that whenever
$\g_0{\in}U$, we have
$d(\rho_n(\g_0)\rho_n(f)y_n,\rho_n(f)y_n){\leq}\eta$.
\end{proof}

\noindent {\bf Remark:} The proof of the first assertion of
Proposition \ref{proposition:continuoussubactionHilbert} for
non-affine actions occurs in subsection \ref{subsection:appendix}
and the proof of the second assertion is found following the proof
of Lemma \ref{lemma:barycentervarepsilonisometry} in subsection
\ref{subsection:barycenter}.  The proof of Proposition
\ref{proposition:convolutionalmostisometrichilbert} for non-affine
actions is found at the end of subsection
\ref{subsection:convolution}.

In the discrete case, we implicitly used finiteness of $K$ to show
that the $K$-displacement of the $\omega$-limit of a sequence is
equal to the $\omega$-limit of the $K$-displacements.  This is
true more generally for sequences which are
$\omega$-equicontinuous on compact sets.

\begin{proposition}
\label{proposition:limitsanddisplacement} Let
$\{y_n\}{\in}\Omega$. Then $\dk(y\so)=\omega$-$\lim\dk(y_n)$.
\end{proposition}

\begin{proof}
We let $k_n$ be the sequence of elements in $K$ such that the
$K$-displacement of $y_n$ is $d(\rho_n(k_n)y_n,y_n)$.  By
Proposition \ref{proposition:ultralimitgeneral}, there is a unique
$k\so=\omega$-$\lim{k_n}$.  Since $y_n{\in}\Omega$, we know that
$d(\rho_n(k)y_n,y_n)$ are equicontinuous functions of $k{\in}K$.
This implies
$\omega$-$\lim{d(\rho_n(k\so)y_n,y_n)}=\omega$-$\lim{d(\rho_n(k_n)y_n,y_n)}$
which suffices to prove the proposition.
\end{proof}

Lastly, we need the following trivial lemma.

\begin{lemma}
\label{lemma:closepointsclosedisplacement} Let $X$ be a metric
space and $\rho$ an $(r,s,\varepsilon,\delta,K)$-almost action of
$\G$ on $X$ at $x$.  Then if $d(x,y)=\eta$, then
$\dk(y){\leq}\delta+2\eta+\varepsilon\eta$.
\end{lemma}

\begin{proof}[Proof of Theorem
\ref{theorem:decreasingdisplacementgen} for non-discrete $\G$] The
proof begins as in the discrete case, we assume the theorem is
false and let $\fh_n,x_n$ be a sequence of Hilbert spaces and
$\rho_n$ be as in the proof for $\G$ discrete, assuming we have
already renormalized so $\delta_n=1$.  We note that to prove the
theorem it suffices to show the existence of some
$h{\in}\mathcal{U}(\G)$, so we may assume that
$\dk(\rho(h)x_n)>C\dk(x_n)$ for every $h{\in}\mathcal{U}(\G)$.
Arguing as in the discrete case, we can produce a isometric limit
action $\rho$ on Hilbert space $\fh$ where the $K$-displacement of
$x\so=\omega$-$\lim{x_n}$ is $1$. By Proposition
\ref{proposition:continuoussubactionHilbert}, $\rho$ is continuous
on a closed affine subspace $\fh'{\subset}\fh$ and for any
sequence $\{y_n\}{\in}{\tilde \fh}$ and any
$g{\in}\mathcal{U}(\G)$, the point $\omega$-$\lim\rho_n(g)y_n$ is
in $\fh'$.  Together with Proposition
\ref{proposition:limitsanddisplacement}, the arguments for the
discrete case imply that for any $\{y_n\}$ with $y\so{\in}\fh'$,
we have $\dk(\rho(h)y_n){\leq}C\dk(y_n)$ for $\omega$-almost every
$n$. If $x\so{\in}\fh'$ this completes the proof.  Otherwise let
$x_n'=\rho(f)x_n$ where $f{\in}\cu(\G)$ with $\supp(f){\subset}K$
and $\supp(f)$ containing a neighborhood of the identity in $\G$.

In this case, we will prove the theorem with $h$ replaced by
$h^{*d}*f$ for a positive integer $d$ such that $4C^d{\leq}C$. We
know from Proposition \ref{proposition:continuoussubactionHilbert}
that the sequence $\{\rho_n(f)x_n\}$ is equicontinuous, as is
$\{\rho_n(h)^i\rho_n(f)x_n\}$ for every positive integer $i$.
Since $\rho_n(f)x_n$ is in the ball of radius $1$ about $x$, Lemma
\ref{lemma:closepointsclosedisplacement} implies that
$\dk(\rho_n(f)x_n){\leq}3+\varepsilon_n$. Therefore, we know that
$\dk(\rho_n(h)^d\rho_n(f)x_n){\leq}C^d(3+\varepsilon_n)$. Choosing
$\eta<\frac{C^d}{10d}$ by $d$ applications of Proposition
\ref{proposition:convolutionalmostisometrichilbert}, we have that
$d_n(\rho_n(h)^d\rho_n(f)x_n,\rho_n(h^{*d}*f)x_n)<\frac{C^d}{10}$
for $\omega$ almost every $n$. Then by Lemma
\ref{lemma:closepointsclosedisplacement}, we know that
$\dk(\rho_n(h^{*d}*f)x_n)<C^d(3+\varepsilon_n+\frac{2+\varepsilon_n}{10})$
for $\omega$ almost every $n$. Since for $n$ large enough,
$C^d(3+\varepsilon_n+\frac{2+\varepsilon_n}{10})<4C^d$, this
implies that, for $\omega$ almost every $n$,
$\dk(\rho_n(h^{*d}*f)x_n)<4C^d<C$, contradicting our assumptions.
\end{proof}

\noindent {\bf Remark:} In the discrete case, it is possible to
prove Theorem \ref{theorem:decreasingdisplacementgen} for the same
function $h$ as in Proposition
\ref{proposition:decreasingdisplacementstandard} and any $C>C_0$
from that proposition.  The reader should note that this is no
longer possible when $\G$ is not discrete, but that $h$ can be
replaced by $h^{*l}$ where $l$ is a constant depending only on
$\G,K$ and $h$.

\subsection{Proofs of Theorem
\ref{theorem:decreasingdisplacementBanach}}
\label{subsection:TproofsBanach}

We indicate the modifications to the proof of Theorem
\ref{theorem:decreasingdisplacementgen} needed to prove Theorem
\ref{theorem:decreasingdisplacementBanach}.

For detailed discussion and definitions about properties of Banach
spaces, the reader should refer to \cite[Chapter 8]{BL} for
positive definite functions and to \cite[Appendix A]{BL} for
uniform convexity.  We recall some consequences and definitions
here.  We will let $\fb$ be a uniformly convex Banach space. Let
$\fb_1$, respectively $\fb_1^*$, be the unit ball. Then there is a
map $j:\fb_1{\rightarrow}\fb_1^*$, called the {\em duality map}
defined by letting $j(x)$ be the unique functional such that
$\|j(x)\|=1$ and $<j(x),x>=1$ such that $j{\inv}$ is uniformly
continuous. (For a proof of uniform continuity and estimates on
the modulus of continuity in terms of the modulus of convexity see
\cite[Appendix A]{BL}.) An easy consequence of the definitions is
the existence of a strictly increasing function $\zeta$ on $[0,1]$
with $\zeta(0)=0$ such that for any pair of vectors
$v,w{\in}\fb_1$ we have $<j(w),v>{\leq}1-\zeta(\varepsilon)$ if
and only if $d(v,w){\geq}\varepsilon$.

The following lemma is used in place of Lemma
\ref{lemma:contraction} above. For any representation $\sigma$ of
$\G$ on a Banach space $\fb$, we denote by $\fb^{\sigma}$ the set
of $\sigma$ invariant vectors.

\begin{lemma}
\label{lemma:almostinvariantequalscontracting} Let $\G$ be a
locally compact, compactly generated group and $K$ a compact
generating set.  Let $\fb$ be a uniformly convex Banach space.
Then for any unitary representation $\sigma$ of $\G$ on $\fb$ the
following are equivalent:
\begin{enumerate}
\item there exists $M>0$ such that for any $\delta{\geq}0$, any
$(K,\delta)$-almost invariant vector $v$ is within $M\delta$ of
$\fb^{\sigma}$;

\item for any function $f{\in}\mathcal{U}_2(\G)$ there exists
$0<C<1$ such that for any $v{\in}\fb$ we have
$d(\sigma(f)v,\fb^{\sigma}){\leq}Cd(v,\fb^{\sigma})$.
\end{enumerate}
\end{lemma}

\begin{proof}
The proof that $(2)$ implies $(1)$ and the reverse implication in
the discrete case are straightforward and similar to the proof of
\cite[Lemma III.1.1]{M}.  Therefore we only give an argument for
$(1)$ implies $(2)$.

Fix a function $f{\in}\mathcal{U}_2(\G)$, then there exists
$\eta>0$ such that $f(\g)>\eta$ for every $\g{\in}K^2$. Fix a
vector $v$ with $d(v,\fb^{\sigma})>0$. By re-scaling and changing
basepoint, we can assume $d(v,\fb^{\sigma})=1$ and in fact that
$d(v,0)=1$ where $0$ is the origin in $\fb$. This uses the fact
that $\fb$ is uniformly convex which implies that there is a point
in $\fb^{\sigma}$ realizing $d(v,\fb^{\sigma})$. There exists
$\g_0{\in}K$ for which $d(\g\g_0v,\g v){\geq}{\frac{1}{M}}$. This
implies that
$$\mu\{\g{\in}K^2|d(\g v, w){\geq}\frac{1}{2M}\}{\geq}\mu(K)$$
for any unit vector $w$.  Applying $j(w)$ to $\sigma(f)v$ we have
that
$$\int{j(w)(f(\g)\sigma(\g)v)}{\leq}1-\zeta(\frac{1}{2M})\mu(K)\eta.$$
Since this is true for any $w$, this implies that
$d(\sigma(f)v,0){\leq}1-\zeta{\frac{1}{2M}}\mu(K)\eta$.
\end{proof}

\noindent We now state a replacement for Proposition
\ref{proposition:decreasingdisplacementstandard}.

\begin{proposition}
\label{proposition:decreasingdisplacementBanach}  If $\G$ has
property $(T)$ and $f{\in}\mathcal{U}_2(\G)$ and $)<C_0<1$, there
exist positive integers $M=M(f,C_0)$ and $m=m(f,C_0)$ such that,
letting $h=f^{*m}$, for any generalized $L^p$
 space $\fb$, any continuous isometric action $\rho$ of $\Gamma$
on $\fb$, and any $x{\in}{\fb}$ we have
\begin{enumerate}
\item $d_{\fb}(x,\rho(h)(x)){\leq}M\dk(x)$

\item $\dk(\rho(h)(x)){\leq}C_0\dk(x)$.
\end{enumerate}
\end{proposition}

\noindent As this is essentially contained in \cite{BG}, we only
provide a sketch.

\begin{proof}[Sketch of Proof]
Let $g$ be the positive definite function on $\fb$.  Then $g$
defines maps $T_t:\fb{\rightarrow}\fh_1$ where $\fh$ is a Hilbert
space and $\fh_1$ is the unit sphere.  The map $T_t$ satisfies
$<Tx,Ty>=g(t(x-y))$ where $g(x)=\exp(-\|x\|^p)$.  One can then
apply the standard proof that any affine action of a group with
property $(T)$ on a Hilbert space has a fixed point, see for
example \cite{HV} or \cite{De}, to produce a $\G$ fixed point on
$\fb$. In fact, the proof shows more.  It produces a constant $C$,
depending only on $\G$ and $K$ such that the $\G$ displacement of
any point $y$ in $\fb$ is bounded by  $C$ times the $K$
displacement.  This then implies that the distance from $y$ to a
fixed point (the barycenter of $\G{\cdot}y$) is bounded by a
constant times the $K$ displacement of $y$.  To verify these facts
one uses the fact that
$d_{\fh}(T_tx,T_ty)^2=2td_{\fb}(x,y)^p+O(t^2)$ for all $t$.  The
existence of $h$ then follows from Lemma
\ref{lemma:almostinvariantequalscontracting} and an argument as in
the proof of Proposition
\ref{proposition:decreasingdisplacementrep}.
\end{proof}

The following fact about ultra-product spaces is left to the
reader, compare \cite{He}.

\begin{proposition}
\label{proposition:ultralimitpositivedefinite} Let $p_i$ be
sequence of numbers with $1<p_i<2$ and $\fb_i$ be a sequence of
generalized $L^p$ spaces. Let $\omega$ be an ultra-filter and
$p=\omega$-$\lim p_i$. Then the function $f(x)=\exp(-\|x\|^p)$ is
positive definite on the ultra-product $\fb$ of $\fb_n$ and $\fb$
is uniformly convex with the same modulus of convexity as $L^p$.
\end{proposition}

\begin{proof}[Proof of Theorem
\ref{theorem:decreasingdisplacementBanach}] In the discrete case,
the proof is a verbatim repetition of the proof of Theorem
\ref{theorem:decreasingdisplacementgen} with Proposition
\ref{proposition:decreasingdisplacementstandard} replaced by
Proposition \ref{proposition:decreasingdisplacementBanach} and
Proposition \ref{proposition:ishilbert} replace by Proposition
\ref{proposition:ultralimitpositivedefinite}.

For non-discrete $\G$, one needs to verify that versions of
Proposition \ref{proposition:continuoussubaction} and
\ref{proposition:convolutionalmostisometrichilbert} still hold,
but modifying the proofs of those statements is straightforward if
we assume all $\G$ actions are affine.
\end{proof}

\section{\bf Inner products on tensor spaces and existence of invariant metrics}
\label{section:parametrizing}

In this subsection we define intrinsic leafwise Sobolev structures
on spaces of tensors on $T\ff$.  For our applications, we need
these structures on the space of functions and the space of
symmetric two forms, but we develop it more generally. We define a
family of norms on leafwise $C^r$ tensors and then complete with
respect to the corresponding metric.  For us the key fact about
the norms we use is that they are invariant under isometries of
the leafwise Riemannian metric.  To illustrate the utility of this
construction, we prove Theorem \ref{theorem:invariantmetric} using
this construction and results from subsection
\ref{subsection:Tresults} and \ref{subsection:TresultsBanach}.

As in subsection \ref{subsection:foliatedresults} we let $X$ be a
locally compact, $\sigma$-compact, metric space, $\ff$ a foliation
of $X$ by manifolds of dimension $n$, and $g_{\ff}$ a leafwise
Riemannian metric. We also let $\nu_{\ff}$ denote the Riemannian
volume (and corresponding measure) on leaves of $\ff$ and assume a
transverse invariant measure $\nu$. We define norms on the set of
tensors which are continuous globally and $C^r$ along leaves of
$\ff$. The definitions given below are standard when $X$ is a
single leaf. To make the norms intrinsic, we work with $k$-jets of
sections of tensor bundles. The special case of functions,
particularly important in our applications, is sections of the
trivial one dimensional vector bundle which corresponds to tensors
of the form $\otimes^0T\ff$. Here $T\ff$ is the tangent bundle to
the foliation $\ff$.  We will denote by $\xi$ an arbitrary bundle
of tensors in $T\ff$ and $\Sect^k(\xi)$ the space of sections of
$\xi$ that are globally continuous and $C^k$ along leaves of
$\ff$.  By {\em globally continuous (resp. measurable) and $C^k$
along leaves of $\ff$} we will always mean that an object is $C^k$
along leaves and varies continuously (resp. measurably) in the
$C^k$ topology transverse to leaves. Particular examples include
vector fields $\xi=T\ff$, symmetric two tensors $\xi=S^2(T\ff^*)$
or functions $\xi=X{\times}{\mathbb R}$.

\noindent {\bf Remark:} For our first proof of Theorem
\ref{theorem:isomrigid}, we allow one additional choice for $\xi$.
Let $\xi=X{\times}\Ra^n$ be a trivial bundle. Given an action
$\rho$ of $\G$ on $X$, we normally would associate the trivial
action on $\xi$.  Instead we allow the possibility of the
existence of finite dimensional unitary representation $\sigma$ of
$\G$ on $\Ra^n$, and define the action of $\rho$ on $\xi$ by
$\rho(\g)(x,v)=(\rho(\g)x,\sigma(\g)v)$. Similarly for any
perturbation $\rho'$ of $\rho$, we define the action $\rho'$ on
$\xi$ by $\rho'(\g)(x,v)=(\rho'(\g)x,\sigma(\g)v)$.

Given $g_{\ff}$, there is a canonical Levi-Civita connection on
$T\ff$ which we denote by $\nabla^T$ associated to the metric. For
any choice of $\xi$, this defines a connection on $\xi$, see for
example III.2 of \cite{KN}, which we will view as
$\nabla^{\xi}:\Sect(T\ff){\times}\Sect(\xi){\rightarrow}\Sect(\xi)$.
Note that $\nabla$ is always invariant under isometries of
$g_{\ff}$. There is also a natural metric on $\xi$ associated to
the metric on $T\ff$ (see for example section 20.8.3 in \cite{D}).
In the particular case of functions, the metric is any metric
given by identifying all fibers with $\mathbb R$, the connection
is given by $\nabla_X{f}=Xf$, and invariance of the connection is
immediate.

We will let $J^k(\xi)$ denote the bundle of leafwise $k$-jets of
sections of $\xi$.  This is a bundle where the fiber over a point
$x$ is the set of equivalence classes of continuous, leafwise
$C^k$ sections where two sections are equivalent if they agree to
order $k$ at the point $x$.  We will denote the fiber by
$J^k(\xi)_x$. There is a natural identification:
$$J^k(\xi){\simeq}{\bigoplus}_{j=0}^k(S^j(T{\ff}^*){\otimes}{\xi}).$$
As a special case of the discussion above, a metric on $T(\ff)$
defines one on $S^j(T{\ff}^*)$ for all $j$
 and together with the metric on $\xi$, this identification induces a metric on
$J^k(\xi)$.  We briefly review the identification above to show
that this metric is indeed invariant under isometries of
$g_{\ff}$. The exposition that follows draws mostly from section
9 of \cite{Palais}, where the interested reader can find more
proofs and explicit constructions.

We can view the leafwise Levi-Civita connection on $T{\ff}$ as a
map
$\nabla^T:\Sect(T\ff){\rightarrow}\Sect(T{\ff}^*){\otimes}\Sect(T\ff)$.
By the discussion above, we also have a connection
$\nabla^{T*}:\Sect(T{\ff^*}){\rightarrow}\Sect({\otimes}^2T{\ff}^*)$.
Similarly the connection on $\xi$ can be viewed as a map
$\nabla^{\xi}:\Sect(\xi){\rightarrow}\Sect(T{\ff}^*{\otimes}\xi)$
by viewing $\nabla_X{\sigma}$ as a one form with $X$ as the
variable.

Similarly, we can define a canonical covariant derivative
$$\nabla^{(i)}:{\Sect}({\otimes}^i(T{\ff}^*){\otimes}\xi){\rightarrow}\Sect({\otimes}^{i+1}(T{\ff}^*){\otimes}\xi)$$
via the formula
$$\nabla^{(i)}(V_1{\otimes}{\cdots}{\otimes}V_i{\otimes}f)=\sum_{j=1}^{i}(V_1{\otimes}{\cdots}{\otimes}
{\nabla^{T*}}V_j{\otimes}{\cdots}{\otimes}V_i{\otimes}f)+V_1{\otimes}{\cdots}{\otimes}V_i{\otimes}\nabla^{\xi}{f}$$
where $f$ is a section of $\xi$ and the $V_i$ are elements of
$T\ff^*$.  The composition
$\nabla^{(k-1)}{\cdots}\nabla^{(1)}\nabla:\Sect(\xi){\rightarrow}\Sect({\otimes}^k(T{\ff}^*){\otimes}
{\xi})$ is called a $k$th covariant derivative.

We now define the total covariant derivative.  First let $S^{(k)}$
be the natural symmetrization operator from the $k$th tensor power
of $T{\ff}$ to $S^k(T{\ff})$ the $k$th symmetric power.  We define
the $k$th total differential
$D_k=S_*^{(k)}{\nabla}^{k-1}{\cdots}\nabla^{(1)}\nabla$.  In other
words, the $k$th total differential is the symmetrization of the
$k$th covariant derivative.  The isomorphism mentioned above
$J^k(\xi){\simeq}{\bigoplus}_{m=0}^k (S^j(T{\ff}^*){\otimes}\xi)$
is given by the map $j^k(v)={\{}D_m(v){\}}_{0{\leq}m{\leq}k}$.  We
then define the metric on $J^k(\ff)$ via this isomorphism and
abuse notation by calling it $g$.  Define a family of norms on
$\Sect(J^k(\xi))$ via:
$$\|u\|^p_{p}=\int_Xg(u_x,u_x)^{\frac{p}{2}}d\nu_{\ff}d\nu,$$
for every $u{\in}\Sect(J^k(\xi))$. Since elements of
$\Sect^k(\xi)$ define elements of $\Sect(J^k(\xi))$ we can
restrict this to an inner product on $\Sect^k(\xi)$ defined by
$$\|u\|^p_{k,p}=\int_X<j^k(u),j^k(u)>^{\frac{p}{2}}d\nu_{\ff}d\nu,$$
for all $u{\in}\Sect^k(\xi)$.

\noindent{\bf Notational convention:} Throughout this paper when
$f$ is leafwise smooth homeomorphism of $(X,\ff)$  and therefore
induces a map on functions or sections of a tensor bundle $\xi$
over $X$, we abuse notations by writing $f$ for the map on
functions or sections. This remark also applies to group actions.

\begin{proposition}
\label{proposition:isometries} Let $f$ be a leafwise isometry of
$(X,\ff,g_{\ff})$ which preserves the transverse invariant measure
$\nu$. Then the action of $f$ on $\Sect^k(\xi)$ and
$\Sect(J^k(\xi))$ preserves all of the norms defined above.
\end{proposition}

\begin{proof}
This is clear from the definition of the inner product and the
fact that isometries of $g_{\ff}$ commute with all the
differential operators used in the construction and that $f$
preserves the measure in the integral above.
\end{proof}


\noindent We now have norms defined $\Sect(J^k(\xi))$ which
restrict to norms defined on $\Sect^k(\xi)$. We define distance
functions on $\Sect(J^k(\xi))$ by $d_p(u,v)={\|u-v\|_p}$ and refer
to the completion with respect to this metric as $L^p(J^k(\xi))$.
Note that $d_p$ restricts to a metric $d_{p,k}$ on $\Sect^k(\xi)$
that is exactly the metric induced by $\|{\cdot}\|_{k,p}$.
Completing $\Sect^{k}(\xi)$ with respect to $d_{p,k}$ we obtain a
standard Sobolev completion of that space, a Banach subspace of
$L^p(J^k(\xi))$, which we denote by $L^{p,k}(\xi,\ff)$. If the
foliation is the trivial foliation by a single leaf $X$, we omit
the $\ff$ and simply write $L^{p,k}(\xi)$ and $L^p(J^k(\xi))$. In
the special case of functions, we use the notation
$L^{p,k}(X,\ff), L^p(J^k(X))$ or $L^{p,k}(X)$ in place of
$L^{p,k}( X{\times}{\mathbb R},\ff), L^p(J^k(X{\times}\Ra)$ or
$L^{p,k}(X{\times}{\mathbb R})$ respectively. It is clear that if
$f$ is a homeomorphism of $X$ as in Proposition
\ref{proposition:isometries} then the action of $f$ on
$\Sect(J^k(\xi))$ and $\Sect^k(\xi)$ extend to isometric actions
on $L^p(\Sect(J^k(\xi))$ and $L^{p,k}(\xi)$. Since any
$u{\in}L^{p,k}(\ff,\xi)$ is a limit of $u_i{\in}\Sect^k(\xi)$ with
respect to the norm above, it follows that $j^k(u_i)$ converge in
$L^p(J^k(\xi))$ to a section we denote by $j^k(u)$ and call the
{\em weak $k$-jet} of $u$.

We also have the following fact about perturbations of isometric
actions which will be used heavily in the next section.

\begin{proposition}
\label{proposition:almostisometric} Let $f$ be a leafwise isometry
of $(X,\ff,g_{\ff})$ which preserves $\nu$ and let
$s{\in}\Sect^k(\xi)$ be $f$ invariant.  If $\xi$ is a trivial
bundle, we let $l=k$, if $\xi$ is non-trivial, we let $l=k+1$. For
any $p_0>1,\varepsilon>0$ and $\delta>0$ there exists a
neighborhood $U$ of the identity in $\Diff_{\nu}^l(X,\ff)$ such
that if $f'$ is an $(U,C^l)$-foliated perturbation of $f$ then

\begin{enumerate}
\item for any $p{\leq}p_0$, the action of $f'$ on
$L^{p}(J^k(\xi))$ (and therefore $L^{p,k}(\xi,\ff)$) is by
$\varepsilon$-almost isometries,

\item the $f'$ displacement of $s$ in $L^{p,k}(\xi,\ff)$ is less
than $\delta$,

\item if $V{\subset}X$ is any $f'$ invariant set of positive
measure, then the action of $f'$ on $L^{p}(J^k(\xi))|_V$ (and
therefore $L^{p,k}(\ff,\xi)|_V$) is by $\varepsilon$-almost
isometries,

\item if $V{\subset}X$ is any $f'$ invariant set of positive
measure, then the $f'$ displacement of $s|_V$ in
$L^{p,k}(\ff,\xi)|_V$ is less than $\delta\mu(V)$.
\end{enumerate}

\noindent Furthermore if $H$ is a topological group and $\rho$ is
a continuous leafwise isometric action of $H$ on $X$, then the
resulting $H$ action on $L^{p}(J^k(\xi))$ (and therefore
$L^{p,k}(\xi,\ff)$) is continuous. The same is true for any
continuous $(U,C^k)$-foliated perturbation $\rho'$ of $\rho$
\end{proposition}

\noindent {\bf Remarks:}
\begin{enumerate}
\item The choice of $l$ is required since while $C^l$
diffeomorphisms act on $C^l$ functions on $X$, they only act on
$C^{l-1}$ sections of any non-trivial tensor bundle $\xi$.  It is
easy to verify that if $f$ is a $C^l$ diffeomorphism and $s$ is a
$C^{l-1}$ section of $\xi$, then
$j^{l-1}(s{\circ}f)=j^{l-1}(s){\circ}j^l(f)$. \item  Since if $s$
is in $\Sect^{\infty}(\xi)$ then $s$ is in $\Sect^{k}(\xi)$ we can
use Proposition \ref{proposition:almostisometric} to study
translates of $C^{\infty}$ sections inside Sobolev spaces. \item
There is no better statement for the $C^{\infty}$ case, since a
$C^{\infty}$ neighborhood of  $f$ is exactly a $C^n$ neighborhood
of  $f$ for some large integer $n$. If $f'$ is close to $f$ in the
$C^n$ topology but not the $C^{n+1}$ topology, then even if $f'$
is $C^{\infty}$, $f'$ is not be an $\varepsilon$-almost isometry
on any space defined using more than $n$ derivatives.  Therefore
we can only obtain an estimates for the $f'$ action on spaces
whose norms depend on no more than $n$ derivatives.
\end{enumerate}

\begin{proof}
First given $\varepsilon$, we find $U$ such that
$(1-\varepsilon)\|s\|_{p}\leq\|s{\circ}f'\|_{p}\leq(1-\varepsilon)\|s\|_{p}$
for any $s{\in}L^p(\Sect(\xi))$ and for any $f'$ which is
$(U,C^l)$-close to $f$. Since continuous sections are dense in
$L^{p}(\Sect(\xi))$ (this follows from the fact that continuous
functions are dense in $L^p$), we can assume that $s$ is
continuous. We can write $s{\circ}f'$ as
$s{\circ}(f{\circ}f{\inv}){\circ}f'=({s}{\circ}f){\circ}(f{\inv}{\circ}f')$.
Since $f$ is an isometry of $L^p(\Sect(\xi))$ it suffices to show
that
$(1-\varepsilon)\|s\|_{p}\leq\|s{\circ}(f{\inv}{\circ}f'){\cdot}\|_{p}\leq(1-\varepsilon)\|s\|_{p}$
for leafwise smooth $s$.  For any $\eta>0$, we can choose $U$, an
open set in $\Diff^l(X,\ff)$ containing the identity, such that
$1-\eta<\|j^l(f{\inv}{\circ}f')(x)\|<1+\eta$ for all $x$. Then the
chain rule implies the pointwise bound
$$(1-\eta)\|j^k(s)(x)\|{\leq}
\|j^k(s{\circ}(f{\inv}{\circ}f'))(f{\inv}(f'(x))\|{\leq}(1+\eta)\|j^k(s)(x)\|.$$
We further restrict $U$ so that the Jacobian of $f'$ along $\ff$
is bounded between $1+\eta$ and $1-\eta$, and then the result
follows from the fact that $f'$ preserves the transverse measure
$\nu$ provided $\varepsilon<(1+\eta)^{p+1}-1$.  This argument also
verifies that $f'$ acts by $\varepsilon$-almost isometries on
$L^{p}(\Sect(\xi))|_V$ for any $V\subset{X}$ of positive measure.

The remaining conclusions follow from the fact that
$\Diff^k(X,\ff)$ acts continuously on $\Sect(J^k(\xi))$ and
$\Sect^k(\xi)$ and therefore on $L^p(J^k(\xi))$ and
$L^{p,k}(\xi)$.
\end{proof}

In order to obtain optimal results, we need to make precise some
notions of H\"older regularity in order to have a norm on
$\Sect^k(\xi)$ where $k$ is not integral.  For the remainder of
this section, we allow the possibility that $k$ is not integral
and let $k'$ to denote the greatest integer less than $k$. Given
$x{\in}X,y{\in}\fL_x$ and a piecewise $C^1$ curve $c$ in $\fL_x$
joining $x$ to $y$, for any natural vector bundle $V$ over $X$, we
denote the parallel translation of $v{\in}V_y$ to $V_x$ by
$P_y^xv$ and by $l(c)$ the length of $c$. We then define
$$\|s\|_k= \|s\|_{k'}+\sup
\frac{\|P^x_yj^{k'}(s)(y)-j^{k'}(s)(x)\|}{l(c)^{k-k'}}$$ where
$k'$ is the least integer not greater than $k$ and the supremum is
taken over $x{\in}X$, $y{\in}\fL_x$ and piecewise $C^1$ curves $c$
in $\fL_x$ joining $x$ to $y$.  It is easy to verify that this
definition agrees with the usual one in the Euclidean case.

We now also make precise the $C^k$ size of a $C^k$ map
$f:Z{\rightarrow}Z$ where $Z$ is a Riemannian manifold and $k$ is
not an integer. This notion is needed to make precise the
conclusion $(4)$ of Theorem \ref{theorem:almostconjugacygen}.  We
already have a notion of pointwise $C^{k'}$ size, defined in
subsection \ref{subsection:generalfoliated}, which we denote by
$\|j^{k'}(f)(x)\|$.  Recall that
$j^{k'}(f)(x):J^{k'}(Z,\Ra)_x{\rightarrow}J^{k'}(Z,\Ra)_{f(x)}$ is
a linear map between vector spaces.  Given a curve $c$ in $Z$, we
can compose $j^{k'}(f)(x)$ with parallel translation $P_{f(x)}^x$
along $c$ to obtain a self-map $P_{f(x)}^x{\circ}j^{k'}(f)(x)$ of
$J^{k'}(Z,\Ra)_x$.  We define the $C^k$ size of $f$ to be
$$\|f\|_k=\sup_x\|j^{k'}(f)(x)\|+\sup\frac{\|P_{f(x)}^x{\circ}j^{k'}(f)(x)\|}{l(c)^{k-k'}}$$
where the supremum is taken over $x{\in}X$, $y{\in}Z$ and
piecewise $C^1$ curves $c$ in $Z$ joining $x$ to $y$.  We can also
measure the $C^{k'}$ size of $f$ on any subset $U$ of $Z$ by
restricting the above supremum to $x{\in}U$.

\begin{proposition}
\label{proposition:sobolev} Let $(X,\ff)$ be a compact foliated
space and $g_{\ff}$ a continuous, leafwise smooth metric on
$(X,\ff)$. Then for any $\xi$ and any $p>1$, there are uniformly
bounded inclusions $L^{p,k}(\tilde
\fL_x,\xi){\subset}\Sect^{k-\frac{d}{p}}(\xi|_{\tilde \fL_x})$ for
all $x$, where $d=\dim(Z)$ and $\tilde \fL_x$ is any covering
space of the leaf through $x$.
\end{proposition}

\begin{proof}
The standard Sobolev embedding theorems provide an bounded
inclusion of $L^{p,k}(\Ra^{d})$ in $C^{k-\frac{d}{p}}(\Ra^{d})$
which easily implies a bounded embedding of
$L^{p,k}(\Ra^{d},\Ra^n)$ in $C^{k-\frac{d}{p}}(\Ra^{d},\Ra^n)$.
Compactness of $X$ and the fact that $g_{\ff}$ is continuous and
leafwise smooth, imply that we can cover $X$ with a finite
collection of charts $(U_i,\phi_i)$ with
$\phi_i(U_i)=V_i{\times}B(0,c)$ and such that there is a uniform
bound on the resulting inclusions
$$L^{p,k}(\ff,\xi|_{B_{\ff}(v_i,c)}){\subset}L^{p,k}(B(0,c),\Ra^n)$$
and
$$C^{k-\frac{p}{d}}(B(0,c),\Ra^n){\subset}\Sect^{k-\frac{d}{p}}(\xi|_{B_{\ff}(v_i,c)})$$
for every $U_i$ and every $v_i{\in}V_i$, where
$B_{\ff}(v_i,c)=\phi{\inv}(v_i{\times}B(0,c))$. So we have
uniformly bounded embeddings
$$L^{p,k}(\ff,\xi|_{B_{\ff}(v_i,c)}){\subset}L^{p,k}(B(0,c),\Ra^n)
{\subset}C^{k-\frac{d}{p}}(B(0,c),\Ra^n){\subset}\Sect^{k-\frac{d}{p}}(\xi|_{B_{\ff}(v_i,c)})$$
for every $U_i$ and every $v_i{\in}V_i$ which suffices to complete
the proof.  It is easy to see that the same bound holds for any
cover $\tilde \fL_x{\rightarrow}\fL_x$.  If $k-\frac{d}{p}$ is not
integral, this does not immediately yield the desired result,
since we only have a H\"older bound at small scales.  However,
since we have a global bound on the $C^0$ norm, it is easy to
convert this small scale H\"older bound to a worse H\"older bound
on all scales.  More precise estimates can be obtained by
following the standard proofs of H\"older regularity in the
Sobolev embedding theorems.
\end{proof}

If we are studying perturbations $\rho'$ of an action $\rho$, in
order to obtain optimal regularity in all proofs, we will need to
know that a certain section $s'$ invariant under $\rho'$ is close
in $L^{p,k}$ type Sobolev spaces to certain $\rho$ invariant
section $s$. The difficulty here is to show that $s'$ is both
invariant under $\rho'$ and close in $L^{p,k}$ for $p>2$
simultaneously. To show this, we will require the following
elementary fact.

\begin{lemma}
\label{lemma:convergingpointwiseae} Let $(X,\mu)$ be a measure
space and $V$ a finite dimensional vector space.  Assume that
$f_n{\in}L^p(X,\mu,V)$ converge in $L^p$ to a function $f$.
Further assume that $\|f_n-f_{n+1}\|_p{\leq}C^n$ where $0<C<1$.
Then $f_n$ converges pointwise almost everywhere to $f$.
\end{lemma}

\begin{proof}
Let $X_n=\{x| |f_n(x)-f_{n+1}(x)|>C^{n/2p}\}$, since
$\|f_n-f_{n+1}\|_p{\leq}C^n$, it follows that $\mu(X_n)<C^{n/2}$.
Then $\{f_n\}$  converges pointwise on the complement of
$X_{\infty}=\cap_{n=1}^{\infty}\cup_{k=n}^{\infty}X_n$. The lemma
follows from the Borel-Cantelli lemma, since $\sum_n\mu(X_n)$
converges, so $\mu(X_{\infty})=0$.
\end{proof}

Many of our uses of this fact could, with slight rewording, be
deduced from the fact that if a sequence of functions $\{f_n\}$
converges to a function $f^p$  $L^p$ and converges to a function
$f^q$ in $L^q$ then $f^p=f^q$ almost everywhere. However, the full
strength of Lemma \ref{lemma:convergingpointwiseae} is required in
the proof of Theorem \ref{theorem:almostconjugacygen}.

\begin{lemma}
\label{lemma:multiplespaces} Let $\G$ be a locally compact,
$\sigma$-compact group with property $(T)$ generated by a compact
set $K$, and let $\rho$ be a leafwise isometric action of $\G$ on
$(X,\ff,g_{\ff})$ and $s$ a $\rho$ invariant section in
$\Sect^k(\xi)$.  Let $l$ be as in Proposition
\ref{proposition:almostisometric}.  For any $p{\geq}2,\eta>0,F>0$
and $0<C<1$, there exists a neighborhood $U$ of the identity in
$\Diff_{\nu}^l(X,\ff)$ and a function $h=h(p)$ in
$\mathcal{U}(\G)$ such that if $\rho'$ is a $(U,C^l)$-foliated
perturbation of $\rho$,
\begin{enumerate}
\item $\rho'(h)^ns$ converge pointwise almost everywhere to a
$\rho'$ invariant section $s'$ in $L^{p,k}(\xi)$, \item
$\|\rho'(h)^ns-s'\|_{p,k}{\leq}\eta$ for all $n{\geq}0$ and, \item
$\|\rho(h)^{n+1}s-\rho(h)^ns\|_{p,k}<C^nF$ for all $n{\geq}0$.
\end{enumerate}
\end{lemma}

{\noindent}{\bf Remarks:}\begin{enumerate} \item For many
applications we only need conclusion $(1)$ and the case of $(2)$
where $n=0$, i.e. that $\|s'-s\|_{p,k}{\leq}\eta$. \item The
reason we do not obtain these estimates in $L^{p,k}$ for all $p,k$
when $\rho'$ is $C^{\infty}$ close to $\rho$ is explained
following Proposition \ref{proposition:almostisometric}.
\end{enumerate}

\begin{proof}
Given $p{\geq}2,\varepsilon>0$ and
$\delta=\min(\frac{F}{M},\frac{\eta(1-C)}{MC})$ for $M$ to be
specified below, by Proposition \ref{proposition:almostisometric}
there is a neighborhood $U$ in $\Diff^{l}(X)$ such that for any
$(U,C^{l})$-perturbation $\rho'$ of $\rho$, it follows that
$\rho'(k)$ is an $\varepsilon$-almost isometry of
$L^p(\Sect^k(\xi,\ff))$ (and therefore of $L^{p,k}(\xi,\ff)$) for
any $p<p_0$ and the $\dk(s)<\delta'$ in any of these spaces.

We choose $h{\in}\mathcal{U}(\G)$ satisfying both  Theorem
\ref{theorem:decreasingdisplacement} and Corollary
\ref{corollary:contractingoperatorbanach}. Then Theorem
\ref{theorem:decreasingdisplacement} shows that $\rho'(h)^ns$
converges exponentially to $s'$ in $L^{2,k}(\xi,\ff)$, which by
Lemma \ref{lemma:convergingpointwiseae}, implies that
$\rho'(h)^ns$ converges pointwise almost everywhere to $s'$. Then
applying Corollary \ref{corollary:contractingoperatorbanach} there
is a constant $M=M(h,C,p)$, such that
$\|\rho(h)^{n+1}s-\rho(h)^ns\|<C^nM\delta$ and therefore
$s'{\in}L^{p,k}(\xi,\ff)$ and
$\|\rho(h)^ns-s'\|_{p,k}{\leq}\frac{MC}{1-C}\delta$ for all
$n{\geq}0$. By our choice of $\delta$ we have
$\|\rho(h)^ns-s'\|{\leq}\eta$ and
$\|\rho(h)^{n+1}s-\rho(h)^ns\|<C^nF$ for all $n{\geq}0$ as
desired.
\end{proof}

To illustrate the application of the results in section
\ref{subsection:Tresults} to perturbations of isometric and
leafwise isometric actions, we now prove Theorem
\ref{theorem:invariantmetric} from the introduction.

\begin{proof}[Proof of Theorem \ref{theorem:invariantmetric}]
We have an action $\rho$ of a group $\G$ with property $(T)$ on
compact manifold $X$ preserving a Riemannian metric $g$.  We view
$g$ as a section of the (positive cone in) the bundle of symmetric
two tensors $S^2(TX)$.  Fix a generating set $K$ of $\G$ and a
choice of $\eta>0$ to be specified below.   Given $\varepsilon>0$
satisfying the hypotheses of Theorem
\ref{theorem:fixedpointsimplegen} and $\delta>0$ to be specified
below, by Proposition \ref{proposition:almostisometric} there is a
neighborhood $U$ in $\Diff^{k+1}(X)$ such that for any
$(U,C^{k+1})$-perturbation $\rho'$ of $\rho$, it follows that
$\rho'(k)$ is an $\varepsilon$-almost isometry of
$L^{2,k}(S^2(TX))$  and the $K$ displacement of $g$ is less than
$\delta$ in this space.  Theorem \ref{theorem:fixedpointsimplegen}
then implies that there is a number $C>0$ depending only on $\G$
and $K$ and a $\rho'(\G)$ invariant section
$g'{\in}L^{2,k}(S^2(TX))$ with $\|g-g'\|_{2,k}{\leq}\eta$ where
$\eta=C\delta$ is specified below. To obtain optimal regularity,
we choose $p>\frac{d}{\kappa}$ and let $U$ satisfy Lemma
\ref{lemma:multiplespaces} for $p$ and $\eta$ specified below and
then Lemma \ref{lemma:multiplespaces} implies that there is a
$\rho'$ invariant section $g'$ such that $\|g-g'\|_{p,k}<\eta$.

By Proposition \ref{proposition:sobolev}, this implies that
$\|g-g'\|_{k-\frac{d}{p}}{\leq}C'\eta$ where $d={\dim(X)}$ and
$C'$ depends only on $X$ and $g$. Since the cone of positive
definite metrics is open in $S^2(TX)$, we can choose $\eta$
depending only on $p$ and $g$, so $g'$ is a $C^{k-\frac{d}{p}}$
Riemannian metric on $X$, invariant under $\rho'(\G)$ and
$C^{k-\frac{d}{p}}$ close to $g$.
\end{proof}

\section{\bf Property $(T)$ and conjugacy}
\label{section:conjugacy}
\noindent

In this section we prove Theorem \ref{theorem:isomrigid} using
Theorem \ref{theorem:fixedpointsimplegen} and Lemma
\ref{lemma:multiplespaces}. In this section, we only consider
$C^k$ perturbations.  The additional arguments required for the
$C^{\infty}$ case are in section
\ref{section:convexityofderivatives}.

\subsection{\bf A proof of Theorem \ref{theorem:isomrigid}}
\label{subsection:isomrigid}

We begin by noting a classical fact about isometric actions.

\begin{proposition}
\label{proposition:fdinvsubspaces} Let $0{\leq}k{\leq}{\infty}$,
let $X$ be a compact $C^k$ manifold and $\rho$ a $C^k$ action of
$\Gamma$ on $X$ such that the image of $\G$ in $\Diff^k(X)$ is
pre-compact. Then there is a positive integer $n$, a homomorphism
$\sigma:\G{\rightarrow}O(n)$ and a $\G$ equivariant $C^k$
embedding $s:X{\rightarrow}\Ra^n$.
\end{proposition}

\noindent {\bf Remark:} For our applications, the fact that $\G$
is precompact in $\Diff^k(X)$ follows from the fact that the
isometry group of a compact Riemannian manifold is compact.

\begin{proof}
Let $C$ be the closure of $\G$ in $\Diff^k(X)$. For $C$ this is
the Mostow-Palais theorem \cite{Mo,P2}. More precisely, for all
$k$, Mostow has proven that, for some $n$, there is a map
$C{\rightarrow}O(n)$ and a $C^k$ equivariant embedding of $X$ into
the Euclidean space $\Ra^n$.  For $k=0$ this is the main result of
\cite{Mo}, for $k>0$, it is proven in section $7.4$ of that paper
by a different method. For $k=\infty$, the same result is proven
in \cite{P2} using the fact that $C$ preserves a $C^{\infty}$
Riemannian metric.  (Mostow's proofs do not explicitly use the
existence of an invariant metric. In the $C^k$ case, Palais'
method produces an equivariant embedding of lower regularity.)
\end{proof}

\noindent {\bf Remark:}  If $k{\geq}2$, and $\G$ preserves a
$C^{k,\alpha}$ Riemannian metric $g$, then one can prove the $C^k$
version of the above theorem by approximating any embedding of
$\G$ in $\Ra^n$ with an embedding defined by eigenfunctions of the
Laplacian.

Given $\sigma$ and $n$ as in Proposition
\ref{proposition:fdinvsubspaces}, we define a trivial bundle
$\xi=X{\times}\Ra^n$ with $\G$ action
$\rho(\g)(x,v)=(\rho(\g)x,\sigma(\g)v)$.  The conclusion of
Proposition \ref{proposition:fdinvsubspaces} is then equivalent to
the existence of a $\G$ invariant section $s:X{\rightarrow}\xi$.
We will show that the perturbed action preserves a section $s'$
close to $s$ and then use the following lemma to produce the
conjugacy.  Given a compact manifold $Y{\subset}\Ra^n$, there is a
neighborhood $U$ of $Y$ in the normal bundle of $Y$ in $\Ra^n$
such that the exponential map $\exp:U{\rightarrow}\Ra^n$ defined
by $\exp(x,v)=x+v$ is a diffeomorphism.  The closest point
projection $\phi$ from $\exp(U)$ to $Y$ is then $C^{\infty}$
(resp. $C^{n-1}$) when $Y$ is $C^{\infty}$ (resp. $C^n$).  This
yields the following:

\begin{lemma}
\label{lemma:sectionstodiffeos} Let $s:X{\rightarrow}\Ra^n$ be a
$C^{\infty}$ embedding, then there exists $\eta$ such that for any
integers $l{\geq}k{\geq}1$ and any $s':X{\rightarrow}\Ra^n$ a
$C^l$ map with $\|s'-s\|_k{\leq}\eta$, the map
$\psi=s{\inv}{\circ}\phi{\circ}s'$ is a $C^k$ small, $C^l$
diffeomorphism of $X$. Furthermore as $\eta{\rightarrow}0$, the
map $\psi$ tends to the identity map in the $C^k$ topology.
\end{lemma}

We now prove Theorem \ref{theorem:isomrigid}. The reader who is
only interested in a result, and not a result with optimal
regularity, may ignore the last sentence of each paragraph and
read the second paragraph assuming $p=2$. For any perturbation
$\rho'$ of $\rho$, we define an action $\rho'$ on $\xi$ by
$\rho'(\g)(x,v)=(\rho'(\g)x,\sigma(\g)v)$.

\begin{proof}[Proof of Theorem \ref{theorem:isomrigid}]
Fix a generating set $K$ for $\G$ and a constant $\eta>0$ to be
specified below. By Proposition \ref{proposition:almostisometric},
given $\varepsilon>0$ satisfying the hypotheses of Theorem
\ref{theorem:fixedpointsimplegen} and $\delta>0$ specified below,
we can choose a neighborhood of the identity
 $U{\subset}\Diff^k(X)$ such that for any $\rho'$ that is $(U,C^k)$ close to
 $\rho$, the map
$\rho'(\g)$ is a $\varepsilon$-almost isometry of $L^{2,k}(\xi)$
for any $\g{\in}K$ and such that $\dk(s)<\delta$ in
$L^{2,k}(\xi)$. Theorem \ref{theorem:fixedpointsimplegen} then
implies that there is a $\rho'(\G)$ invariant section
$s'{\in}L^{2,k}(\xi)$ with $\|s-s'\|_{2,k}{\leq}\eta$ where
$\eta=C\delta$ and $C>0$ depends only on $\G$ and $K$.  To obtain
a $C^{k-\kappa}$ conjugacy, we choose $p<\frac{d}{\kappa}$ and
choose $U$ to satisfy Lemma \ref{lemma:multiplespaces} for our
choices of $\eta$ and $p$ and then Lemma
\ref{lemma:multiplespaces} implies that there is a $\rho'$
invariant section $s'$ with $\|s-s'\|_{p,k}<\eta$.

Proposition \ref{proposition:sobolev} then implies that $s'$ is
$C^{k-\frac{d}{p}}$ and that $\|s-s'\|_{k-\frac{d}{p}}<C_0\eta$
where $C_0$ is an absolute constant depending only on $X,p$ and
$g$. We can view $s':X{\rightarrow}\Ra^n$ as a $\G$ equivariant
$C^{k-\frac{d}{p}}$ map from $X$ to $\Ra^n$ where the action on
$X$ is given by $\rho'$ and the action on $\Ra^n$ is given by
$\sigma$.  By choosing $\eta$ (and therefore $U$) sufficiently
small, the map $s':X{\rightarrow}\Ra^n$ is $C^{k-\frac{d}{p}}$
close to $s$. Then by Lemma \ref{lemma:sectionstodiffeos}, the map
$\psi=s{\inv}{\circ}{\phi}{\circ}s'$ is a $C^{k-\frac{d}{p}}$
small diffeomorphism of $X$.  Since $s,s'$ and $\phi$ are all
equivariant, $\psi$ is a conjugacy between the $\rho'(\G)$ action
on $X$ and the $\rho(\G)$ action on $X$.
\end{proof}

\subsection{\bf Another Proof of Theorem \ref{theorem:isomrigid}}
\label{subsection:isomrigid2}

In this section we give another proof of Theorem
\ref{theorem:isomrigid} which gives somewhat lower regularity, but
which generalizes to prove Theorem
\ref{theorem:almostconjugacygen}.

We will denote points in $X{\times}X$ by $(x_1,x_2)$ and denote
the diagonal in $X{\times}X$ by $\Delta(X)$. Given any group $\G$
acting on a manifold $X$, we will denote by $\bar \rho$ the
diagonal action of $\G$ on $X{\times}X$ given by $\bar
\rho(\g)(x_1,x_2)=(\rho(\g)x_1,\rho(\g)x_2)$.

We begin with two elementary facts. Recall that a normal
neighborhood of $x$ is the image under the Riemannian exponential
map of an open ball $B$ in $T_xX$, such that $\exp_x|_B$ is a
diffeomorphism and $d_X(x,\exp(v))=d_{T_xX}(0,v)$. It is immediate
that $d(x,\cdot)^2$ is a smooth function on any normal
neighborhood of $x$. Let $N(x)$ be the maximal radius of a normal
neighborhood of $x$ in $X$.  On $X{\times}X$ we have a Riemannian
metric on $g{\times}g$ and the induced distance function. By
$B(x,\varepsilon)$ we denote the ball of radius $\varepsilon$
around $x$ in $X$, by $B((x,y),\varepsilon)$ the ball of radius
$\varepsilon$ around $(x,y)$ in $X{\times}X$.  Since
$\{x\}{\times}X$ is totally geodesic in $X{\times}X$, we have
$B((x,x),\varepsilon ){\cap }\{x\}{\times}X=B(x,\varepsilon)$.

\begin{proposition}
\label{proposition:invariantfunctions} Let $\rho$ be an isometric
action of any group $\G$ on any Riemannian manifold $X$. If we
further assume the function $N(x)>d$ for some $d>0$ and all
$x{\in}X$, then for any $0<\varepsilon<\frac{d}{2}$, there exists
an invariant smooth function $f$ on $X{\times}X$ such that:
\begin{enumerate}
\item $f$ takes the value $0$ on $\Delta(X)$,

\item $f{\geq}0$ and $f(x,y)>0$ if $x{\neq}y$,

\item for any $x{\in}X$, the restriction of $f$ to
$\{x\}{\times}X$ satisfies $f{\geq}1$ outside of $B((x,x),
\varepsilon)$,

\item  the Hessian of $f$ restricted to $\{x\}{\times}X$ is
positive definite on the closure of $B((x,x),\varepsilon)$

\end{enumerate}
\end{proposition}

\begin{proof}  The
action $\bar \rho$ leaves invariant any function on $X{\times}X$
which is a function of $d(x_1,x_2)$ and we define $f$ as such a
function. To define $f$ we first define
$f_{x_0}:X{\rightarrow}\Ra$ by
$\frac{1}{\varepsilon^2}d(x,x_0)^2$. Our assumptions on
$\varepsilon$ imply that  $B(x_0,2\varepsilon)$ is contained in a
normal neighborhood of $x_0$ and so $d(x,x_0)^2$ is a smooth
function of $x$ and $x_0$ inside $B(x_0,2\varepsilon)$, see for
example \cite[IV.3.6]{KN}. It is clear that
$f(x_1,x_2)=f_{x_1}(x_2)$ satisfies all the requirements except
smoothness on points at distance greater than $2\varepsilon$ from
$\Delta(X)$. We merely need to change $f_{x_0}$ outside
$B(x_0,\varepsilon)$ to produce a smooth $f_{x_0}$ while keeping
$f_0{\geq}1$ outside $B(x_0,\varepsilon)$. This is easily done by
 choosing
any smooth function $g:\Ra{\rightarrow}\Ra$ such that $g$ agrees
with $\frac{1}{\varepsilon^2}x^2$ to all orders for all
$x{\leq}\varepsilon$ with $g{\geq}1$ for all $x>\varepsilon$ and
$g=1$ for all $x{\geq}2\varepsilon$. We then let
$f(x_1,x_2)=g(d(x_1,x_2))$.
\end{proof}

\begin{proposition}
\label{proposition:perturbedinvariantfunctions} Let $X$ be a
Riemannian manifold and $f$ a function on  $X{\times}X$ such that:
\begin{enumerate}
\item $f$ takes the value $0$ on $\Delta(X)$,

\item $f{\geq}0$,

\item for any $x{\in}X$, the restriction of $f$ to
$\{x\}{\times}X$ satisfies $f{\geq}1$ outside
$B((x,x),\varepsilon)$,

\item  the Hessian of $f$ restricted to $\{x\}{\times}X$ is
positive definite on the closure of $B((x,x),\varepsilon)$.
\end{enumerate}
Let $f'$ be a function which is $C^k$ close to $f$ where
$k{\geq}2$. Then for every $x$, the restriction of $f'$ to
$\{x\}{\times}X$ has a unique global minimum at a point $(x,x')$
which is close to the point $(x,x)$. Furthermore, if we let
$X'=\{(x,x')| f'(x') \text{ is the global minimum of } f' \text{
on } \{x\}{\times}X\}$ then $X'$ is a $C^{k-1}$ embedded copy of
$X$ which is $C^{k-1}$ close to $\Delta(X)$.
\end{proposition}

\noindent {\bf Remark:} The last statement of the Proposition
means that $X'$ is close to $\Delta(X)$ in the $C^{k-1}$ topology
on $C^{k-1}$ submanifolds of $X{\times}X$.  This actually suffices
to imply that $X'$ is diffeomorphic to $X$ by a normal projection
argument like the one used to prove Lemma
\ref{lemma:sectionstodiffeos}.

\begin{proof}
 Let
$B=\{x\}{\times}B(x,\varepsilon){\subset}\{x\}{\times}X$. We first
verify the existence of $(x,x')$ in $B$. Since $f'$ is $C^k$ close
to $f$ we have that $f'{\geq}\frac{1}{2}$ outside $B$ and $f'$ is
close to zero near $(x,x)$. We look at all local minima of $f'$ on
$\bar B$, the closure of $B$. Since $f'$ is close to $f$, at least
one such minimum occurs in $B$. Since $k{\geq}2$, if $f'$ is
sufficiently $C^k$ close to $f$, the Hessian of $f'$ is positive
definite on $B$, which implies there is exactly one local minimum
on $B$, say at $(x,x')$.  Since $f'$ is $C^k$ close to $f$, it is
easy to see that $f'(x,x')$ must be close to zero and that
$f'(x,y)$ must be close to one if $f'|_{\{x\}{\times}X}$ has a
local minimum at $(x,y)$ and $y$ is not in $\bar B$. Therefore
$x'$ is the unique global minimum of $f'$ on $\{x\}{\times}X$.

Given a function $g:X{\times}X{\rightarrow}\Ra$ we denote by
$D_2g$ the derivative with respect to the second variable, which
is naturally a map from $X{\times}TX$ to $\Ra$.   To see that $X'$
is a smooth submanifold $C^{k-1}$ we note that $X'$ is the set of
zeros of $D_2f':X{\times}TX{\rightarrow}\Ra$ in a neighborhood
$N_{\varepsilon}(\Delta(X)){\subset}X{\times}TX$. Our assumption
on the Hessian implies that these are regular values so $X'$ is
$C^{k-1}$ submanifold since $D_2f'$ is $C^{k-1}$. That $X'$ is
diffeomorphic to $X$ follows from the fact that $X'$ is $C^{k-1}$
close to $\Delta(X)$. This is immediate since $D_2f'$ is $C^{k-1}$
close to $D_2f$.
\end{proof}

For the remainder of this section $X$ will be a compact Riemannian
manifold, $\G$ will be a locally compact group with property $(T)$
and $K$ will be a fixed compact generating set,  $\rho$ will be an
isometric action of $\G$ on $X$ and $\rho'$ will be a $C^k$
perturbation of $\rho$, where $k>2$.  We will denote by $\bar
\rho'$ the $\G$ action on $X{\times}X$ given by perturbing in the
second factor: $\bar
\rho'(\g)(x_1,x_2)=(\rho(\g)x_1,\rho'(\g)x_2)$.  This induces
actions on various spaces of functions which we also denote by
$\bar \rho'$.

We now prove Theorem \ref{theorem:isomrigid}. We only give a proof
with small loss of derivatives. The reader interested in lower
regularity results depending only on Hilbert space techniques can
produce a proof by combining this one with the proof in subsection
\ref{subsection:isomrigid}.

\begin{proof}[Proof of Theorem \ref{theorem:isomrigid}]
We first choose a function $f$ invariant under $\bar \rho$ as in
Proposition \ref{proposition:invariantfunctions}.  Given
$\kappa>0$, we choose $p$ with $\kappa<\frac{d}{p}$ where
$d=\dim(X{\times}X)=2\dim(X)$.  We make a choice of $\eta>0$,
depending on $p$, to be specified below.  We choose $U$ satisfying
Lemma \ref{lemma:multiplespaces} for our choices of $p$ and $\eta$
and then Lemma \ref{lemma:multiplespaces} implies that there is a
$\bar \rho'$ invariant function $f'$ with  $\|f-f'\|_{p,k}<\eta$.

 Proposition \ref{proposition:sobolev}
implies that $\|f-f'\|_{k-\frac{d}{p}}<C_0\eta$ where $C_0$
depends only on $g$ and $p$.  Choosing $\eta$ small enough and
applying Proposition
\ref{proposition:perturbedinvariantfunctions}, we see that we have
a submanifold $X'{\subset}X{\times}X$ which is diffeomorphic to
$X$, $C^{k-\frac{d}{p}-1}$ close to $\Delta(X)$ and $\bar
\rho'(\G)$-invariant.  The first two claims are contained in that
proposition, the last follows from the definition $X'=\{(x,x')|
f'(x') \text{ is the global minimum of } f' \text{ on }
\{x\}{\times}X\}$. We let $p_i:X'{\rightarrow}X$ be the
restriction to $X'$ of the projection
$\pi_i:X{\times}X{\rightarrow}X$ on the $i$th factor where
$i=1,2$. Note that each $\pi_i$ and therefore each $p_i$ is an
equivariant map, where we view the first projection as to $X$
equipped with the action $\rho$ and the second as to $X$ equipped
with the action $\rho'$.  Since $X'$ is $C^{k-\frac{d}{p}-1}$
close to the diagonal, each $p_i$ is a $C^{k-\frac{d}{p}-1}$
diffeomorphism, and the map $p_1{\inv}{\circ}p_2$ is a
$C^{k-\frac{d}{p}-1}$ small diffeomorphism. Therefore
$p_1{\inv}{\circ}p_2$ is a $C^{k-\frac{d}{p}-1}$ conjugacy between
$\rho$ and $\rho'$.
\end{proof}

\section{Additional estimates and $C^{\infty,\infty}$ Local Rigidity} \label{section:convexityofderivatives}

In this section, we prove a key lemmas on regularity  in the
context of isometric actions and their perturbations. From this,
we deduce $C^{\infty,\infty}$ local rigidity in Theorem
\ref{theorem:isomrigid}.    In all that follows $\G$ will be a
locally compact group with property $(T)$ of Kazhdan and $K$ will
be a fixed, compact generating set for $\G$, containing a
neighborhood of the identity in $\G$.   Furthermore, for
simplicity of exposition, the letters $k$ and $l$ below always
denote integers.

The strategy of the proof of the $C^{\infty,\infty}$ version of
Theorem \ref{theorem:isomrigid} is motivated by analogy with the
iterative methods of $KAM$ theory but does not follow a $KAM$
algorithm, see Appendix \ref{appendix:kam} for discussion.

\smallskip
{\noindent}{\bf Remark:} In order to prove optimal regularity, we
make use here of Corollary
\ref{corollary:contractingoperatorbanach}  and the resulting
estimates in $L^p$ type Sobolev spaces. This allows us to give
proofs that imply that, in the context of Theorem
\ref{theorem:isomrigid}, a $C^{\infty}$ action $\rho'$ that is
sufficiently $C^2$ close to an isometric action $\rho$ is
conjugate back to $\rho$ by a $C^{\infty}$ map which is
$C^{2-\kappa}$ small, for $\kappa$ depending on the $C^2$ size of
the perturbation.  The reader only interested in obtaining a
$C^{\infty}$ conjugacy under some circumstances, rather than
optimal circumstances, can easily modify the proofs to use only
$L^2$ type Sobolev spaces and Theorem
\ref{theorem:decreasingdisplacement} instead.

We first state a proposition and lemma for isometric actions.  We
prove Proposition \ref{proposition:betterregularityunfoliated}
from Lemma \ref{lemma:smootherinvariantsections} and some results
in subsection \ref{subsection:isomrigid}, and then use Proposition
\ref{proposition:betterregularityunfoliated} to prove
$C^{\infty,\infty}$ local rigidity for isometric actions. Lemma
\ref{lemma:smootherinvariantsections} will be proven later in this
section. For notational convenience in the statement of this
proposition and the proof of the $C^{\infty}$ case of Theorem
\ref{theorem:isomrigid}, it is convenient to fix right invariant
metrics $d_l$ on the connected components of $\Diff^l(X)$ with the
additional property that if $\varphi$ is in the connected
component of $\Diff^{\infty}(X)$, then $d_l(\varphi,
\Id){\leq}d_{l+1}(\varphi,\Id)$.  To fix $d_l$, it suffices to
define inner products $<,>_l$ on $\Vect^l(X)$ which satisfy
$<V,V>_l{\leq}<V,V>_{l+1}$ for $V{\in}\Vect^{\infty}(X)$. Fixing a
Riemannian metric $g$ on $X$, it is straightforward to introduce
such metrics using the methods of section
\ref{section:parametrizing}.

\begin{proposition}
\label{proposition:betterregularityunfoliated} Let $X$ be a
compact Riemannian manifold and $\rho$ be an isometric action of
$\G$ on $X$. Then for every integer $k{\geq}2$ and every integer
$l{\geq}k$ and every $\varsigma>0$ there is a neighborhood $U$ of
the identity in $\Diff^k(X)$ such that if $\rho'$ is a
$C^{\infty}$ action of $\G$ on $X$ with
$\rho(\g)\inv\rho'(\g){\in}U$ for all $\g{\in}K$ then there exist
a sequence $\psi_n{\in}\Diff^{\infty}(X)$ such that $\psi_n$
converge to a diffeomorphism $\psi$ in $\Diff^l(X)$ and
$\psi_n\circ{\rho'}{\circ}\psi_n{\inv}$ converges to $\rho$ in
$C^l$ and $d_{k-1}(\psi_n,\Id)<\varsigma$ for all $n$.
\end{proposition}

This proposition is a consequence of the following lemma
concerning regularity of invariant sections.  We will only use
this lemma for the trivial bundle $\xi=X{\times}\Ra^n$ equipped
with the $\G$ actions $\rho(\g)(x,v)=(\rho(\g)x,\sigma(\g)v)$ and
$\rho'(\g)(x,v)=(\rho'(\g)x,\sigma(\g)v)$ where
$\sigma:\G{\rightarrow}O(n)$ is fixed, so we do not consider more
general tensor bundles $\xi$.  A similar statement is true in the
general context, though one needs to replace $\Diff^k(X)$ in the
statement with $\Diff^{k+1}(X)$.

\begin{lemma}
\label{lemma:smootherinvariantsections} Let $\G,X,\rho$ be as in
Proposition \ref{proposition:betterregularityunfoliated}, let
$\xi=X{\times}\Ra^n$ and let $s$ be a $\rho(\G)$ invariant section
of $\xi$. Then for every integer $k{\geq}2$, every integer
$l{\geq}k$ and every $\eta>0$ there is a neighborhood $U$ of the
identity in $\Diff^k(X)$ such that if $\rho'$ is a $C^{\infty}$
action of $\G$ on $X$ with $\rho(\g)\inv\rho'(\g){\in}U$ for all
$\g{\in}K$ then the sequence $s_n=\rho'(h)^ns$ satisfies:
\begin{enumerate}
\item $\|s-s_n\|_{k-1}<\eta $ for all $n$ and, \item $s_n$
converges in $\Sect^l(X)$ to a $\rho'$ invariant section $s'$.
\end{enumerate}
\end{lemma}

\noindent {\bf Remarks:} We defer the proof of Lemma
\ref{lemma:smootherinvariantsections} until later in this section.
 Given a positive integer $l>k$, the proof of the lemma only
requires that $\rho'$ is $C^{2l-k+1}$ rather than $C^{\infty}$. By
shrinking $U$, it is possible to show the same result when $\rho'$
is $C^{l+1}$.

\begin{proof}[Proof of Proposition
\ref{proposition:betterregularityunfoliated}] The proof is very
similar to the argument in subsection \ref{subsection:isomrigid}.
We apply Proposition \ref{proposition:fdinvsubspaces} to the
action $\rho$, which produces a representation
$\sigma:\G{\rightarrow}\Ra^n$ and an equivariant embedding
$s:X{\rightarrow}\Ra^n$. We let $\xi=X{\times}\Ra^n$ and define a
action of $\G$ on $\xi$ as specified before Lemma
\ref{lemma:smootherinvariantsections}. Then $s$ is $\rho(\g)$
invariant for every $\g{\in}\G$.  Given $\eta>0$, Lemma
\ref{lemma:smootherinvariantsections} implies that there is a
neighborhood $U$ of the identity in $\Diff^k(X)$ such that for any
action $\rho'$ with $\rho(\g){\inv}\rho'(\g)$ in $U$ for all
$\g{\in}K$ and the action $\rho'$ on $\xi$ defined before the
statement of Lemma \ref{lemma:smootherinvariantsections}, we have
that $s_n=\rho'(h)^ns$ satisfy $\|s-s_n\|_{k-1}<\eta$ and $s_n$
converge in $\Sect^l(X)$ to a $\rho'$ invariant section $s'$. It
is clear that each $s_n$ is $C^{\infty}$. Choosing $\eta$ small
enough and applying Lemma \ref{lemma:sectionstodiffeos}, we see
that the maps $\psi_n=s{\inv}{\circ}{\phi}{\circ}s_n$ are
$C^{k-1}$ small, $C^{\infty}$ diffeomorphisms of $X$, where $\phi$
is the normal projection from a neighborhood of $s(X)$ in $\Ra^n$
to $s(X)$. Letting $\psi=s{\inv}{\circ}{\phi}{\circ}s'$, it is
clear that $\psi_n$ converge to $\psi$ in $\Diff^l(X)$ since $s_n$
converge to $s$ in $\Sect^l(\xi)$.  That $\psi$ is a conjugacy
between the actions $\rho'$ and $\rho$ follows as in the proof of
Theorem \ref{theorem:isomrigid} in subsection
\ref{subsection:isomrigid}.
\end{proof}

\begin{proof}[Proof of $C^{\infty,\infty}$ local rigidity in
Theorem \ref{theorem:isomrigid}] If $\rho'$ is a $C^{\infty}$
perturbation of $\rho$, then there exists some $k>1$, such that
$\rho'$ is $C^k$ close to $\rho$.  We fix a sequence of positive
integers $k=l_0<l_1<l_2<\cdots<l_i<\ldots$ and will construct a
sequence of $C^{\infty}$ diffeomorphisms $\phi_i$ such that the
sequence $\{\phi_n{\circ}{\ldots}{\circ}\phi_1\}_{n{\in}\Na}$
converges in the $C^{\infty}$ topology to a conjugacy between
$\rho$ and $\rho'$.

We let $\phi^i=\phi_i{\circ}{\ldots}{\circ}\phi_1$ and
$\rho_i=\phi^i{\circ}\rho'{\circ}(\phi^i){\inv}$ and construct
$\phi_i$ inductively such that
\begin{enumerate}
\item $\rho_i$ is sufficiently $C^{l_i}$ close to $\rho$ to apply
Proposition \ref{proposition:betterregularityunfoliated} to
$\rho_i$ and $\rho$ with $l=l_{i+1}$ and
$\varsigma=\frac{1}{2^{i+1}}$, \item
$d_{l_i}(\phi_i,\Id)<\frac{1}{2^i}$ and, \item
$d_{l_{i-1}}(\rho_i(\g){\circ}\rho(\g){\inv},\Id)<\frac{1}{2^i}$
for every $\g{\in}K$.
\end{enumerate}
Given $\phi^i$ and therefore $\rho_i$, we construct $\phi_{i+1}$.
We have assumed that $\rho_i$ is close enough to $\rho$ in the
$C^{l_i}$ topology to apply Proposition
\ref{proposition:betterregularityunfoliated} with $l=l_{i+1}$ and
$\varsigma=\frac{1}{2^{i+1}}$. Then we have a sequence of
diffeomorphisms $\psi_n{\in}\Diff^{\infty}(X)$ such that
$\psi_n{\circ}\rho_i{\circ}\psi_n{\inv}$ converges to $\rho$ in
the $C^{l_{i+1}}$ topology and
$d_{k-1}(\psi_n,\Id)<\frac{1}{2^{i+1}}$. We choose $n_i$
sufficiently large so that
$\rho_{i+1}=\psi_{n_i}{\circ}\rho_{i-1}{\circ}\psi_{n_i}{\inv}$ is
close enough to $\rho$ in the $C^{l_{i+1}}$ topology to apply
Proposition \ref{proposition:betterregularityunfoliated} with
$l=l_{i+2}$ and $\varsigma=\frac{1}{2^{i+2}}$ and so that
$d_{l_i}(\rho_i(\g)\rho(\g){\inv},\Id){\leq}\frac{1}{2^{i+1}}$ and
then let $\phi_{i+1}=\psi_{n_i}$.

To start the induction it suffices that $\rho'$ is sufficiently
$C^k$ close to $\rho$ to apply Proposition
\ref{proposition:betterregularityunfoliated} with $l=l_1$ and
$\varsigma=\frac{1}{2}$.

It remains to show that the sequence
$\{\phi_n{\circ}{\ldots}{\circ}\phi_1\}_{n{\in}\Na}$ converges in
the $C^{\infty}$ topology to a conjugacy between $\rho$ and
$\rho'$.  Combining condition $(2)$ with the fact that
$d_{l_i}(\phi_m,\Id){\leq}d_{j}(\phi_m,\Id)$ for all $j{\geq}l_i$,
and the fact that $d_{l_i}$ is right invariant implies that
$d_{l_{i-1}}(\phi_m,\Id)=d_{l_{i+1}}(\phi^m,\phi^{m-1}){\leq}\frac{1}{2^m}$
for all $m{\geq}i$.  This implies that $\{\phi^m\}$ is a Cauchy
sequence in $\Diff^{l_i}(X)$ for all $i$, and therefore $\phi^m$
converge in $\Diff^{\infty}(X)$.  Similarly, condition $(3)$
implies $\rho_m$ converges to $\rho$ in the $C^{\infty}$ topology.
\end{proof}

\noindent{\bf Remark:} The proof above can be made to work in a
more general setting.  Given an action $\rho$ such that for any
large enough $k$ and any $l$ larger than $k$ and any action
$\rho'$ which is sufficiently $C^k$ close to $\rho$, we can find a
conjugacy between $\rho$ and $\rho'$ which is $C^l$ and $C^{k-n}$
small for a number $n$ which does not depend on $l$ or $k$, then
we can use the method above to produce a $C^{\infty}$ conjugacy.
More precisely, we need a bound on the $C^{k-n}$ size of the
conjugacy that depends only on the $C^k$ size of the perturbation.
To apply the argument in this setting, one produces a $C^l$
conjugacy $\varphi$ and then approximates it in the $C^l$ topology
by a $C^{\infty}$ map $\tilde \varphi$ which will play the role of
$\psi_n$ in the argument above.  We use this argument to prove
$C^{\infty,\infty}$ local rigidity in \cite{FM2}.

Before we  proceed to prove Lemma
\ref{lemma:smootherinvariantsections}, we need two additional
estimates. Similar estimates are used in KAM theory. The first is
a convexity estimate on derivatives, which is also used in the
proof of Hamilton's $C^{\infty}$ implicit function theorem, and
which we take from \cite{Ho}.  To be able to prove a foliated
variant of Lemma \ref{lemma:smootherinvariantsections} below, we
state these estimates in the context of foliated spaces. For the
next two lemmas, let $(X,\ff, g_{\ff})$ be a foliated space
equipped with a leafwise Riemannian metric as described in section
\ref{section:parametrizing}. For our applications in subsection
\ref{subsection:smootherinvariantfunctionsonpairs}, it is
important that $X$ need not be compact in either of the following
lemmas.

\begin{lemma}
\label{lemma:convexityofderivatives} Let $a,b,c$ be integers and
$0<\lambda<1$ such that $c=a(1-\lambda)+b\lambda$ and let
$f{\in}\Sect^k(\xi,\ff)$. Then there is a constant $B$ depending
only on $X,\ff$ and $g_{\ff}$ and $b$ such that:

$$\|f\|_c{\leq}B\|f\|_a^{1-\lambda}\|f\|_b^{\lambda}$$.
\end{lemma}

\noindent For $a,b,c$ not necessarily integral, this lemma is
proven for functions on $\Ra^n$ in appendix A of \cite{Ho}.  This
implies the proposition as stated by standard manipulations as in
the proof of Proposition \ref{proposition:sobolev}.

Given a collection elements
$\phi_1{\ldots},\phi_n{\in}\Diff^{\infty}(X,\ff)$ we require a
certain type of bound on the norm of the composition
$\phi_1{\circ}{\cdots}{\circ}\phi_n$ as an operator on $k$-jets of
tensors.  Recall that we have a pointwise norm $\|j^k(\phi)(x)\|$
defined to be the operator norm of
$$j^k(\phi)(x):J^k(X,\ff)_x{\rightarrow}J^k(X,{\ff})_{\phi(x)}.$$
Then we can define the $k$ norm of $\phi$ by
$\|\phi\|_k=\sup_{X}\|j^k(\phi)(x)\|$.  Though we did not find the
following precise estimate in the literature, this type of
estimate is typical in KAM theory. In Appendix
\ref{appendix:compositions}, we give a proof that may be new, at
least in that it makes no reference to coordinates.

\begin{lemma}
\label{lemma:estimateoncompostion} Let
$\phi_1,{\ldots},\phi_n{\in}\Diff^{k}(X,\ff)$. Let
$N_k=\max_{1{\leq}i{\leq}n}\|\phi_i\|_k$ and
$N_1=\max_{1{\leq}i{\leq}n}\|\phi_i\|_1$. Then there exists a
polynomial $Q$ depending only on the dimension of the leaves of
the foliation and $k$ such that:
$$\|\phi_1{\circ}{\ldots}{\circ}\phi_n\|_k{\leq}N_1^{kn}Q(nN_k)$$
for every $n{\in}\Na$.
\end{lemma}

\noindent This lemma has immediate consequences for the operator
norms of $\rho'(h)$ on $C^k(X,\ff)$ which we denote by
$\|\rho'(h)\|_k$.

\begin{corollary}
\label{corollary:estimateoncompostion} Under the assumptions of
Lemma \ref{lemma:smootherinvariantsections}, for any
$h{\in}\mathcal{U}(\G)$, we have the following estimates:
$$\|\rho'(h)^n\|_k{\leq}N_1^{kn}Q(nN_k)$$
where $Q$ is the same polynomial as in Lemma
\ref{lemma:estimateoncompostion} above and
$N_i=\max_{\supp(h)}\|\rho'(\g)\|_i$.
\end{corollary}

\noindent{\bf Remark:} We require this estimate to be able to
estimate the size of $\rho'(h)^ns$ in the $C^l$ topology, even
when the group action $\rho'$ is only $C^k$ close to $\rho$ for
some $k<l$. We do not know of another way to obtain such an
estimate.

\begin{proof}
It follows from the definition that
$$\rho'(h)^n=\int_{\G}h^{*n}\rho'(g)=$$
$$\int_{\G}{\cdots}\int_{\G}h(\g_1){\ldots}h(\g_n)\rho'(\g_1){\ldots}\rho'(\g_n).$$
One then applies Lemma \ref{lemma:estimateoncompostion} applied to
each product of the form $\rho'(\g_1){\ldots}\rho'(\g_n)$ and
integrates.
\end{proof}

The polynomial $Q$ is computable in a straightforward manner for
any given $k$ and dimension as follows easily from the proof.

\begin{proof}[Proof of Lemma
\ref{lemma:smootherinvariantsections}] For the proof of this
lemma, we let $p$ be such that $\frac{\dim(X)}{p}<1$. By Lemma
\ref{lemma:multiplespaces}, for any $0<C<1$ and $F>0$ and
$\beta>0$, we can choose a neighborhood $U$ of the identity in
$\Diff^k(X)$ such that if $\rho'$ is a $C^{\infty}$ action with
$\rho'(\g)\rho(\g){\inv}{\in}U$ for all $\g{\in}K$, there exists
$h{\in}\mathcal{U}(\G)$ such that $\rho'(h)^ns$ converges to a
$\rho'$ invariant section $s'$ where
$$\|\rho'(h)^ns-s'\|_{p,k}{\leq}\beta$$
\noindent for all $n$ and
$$\|\rho'(h)^{n+1}s-\rho'(h)^ns\|_{p,k}{\leq}C^nF.$$

Proposition \ref{proposition:sobolev} then implies that
\begin{equation}
\label{equation:initialestimate1}
\|\rho'(h)^{n+1}s-\rho'(h)^ns\|_{k-1}{\leq}C^nAF
\end{equation}
\noindent and
$$\|\rho'(h)^ns-s'\|_{k-1}{\leq}A\beta$$
\noindent where $A$ depends only on $(X,g)$.  The last inequality
implies the first conclusion of Lemma
\ref{lemma:smootherinvariantsections} provided we chose
$\beta<\frac{\eta}{A}$.

We will show that, possibly after shrinking $U$, $\rho'(h)^ns$
satisfies
\begin{equation}
\label{equation:keystep1}
\|(\rho'(h)^{n+1}s-\rho'(h)^ns)\|_{l}{\leq}C'^nP(nF_l)F
\end{equation}
where $P$ is a fixed polynomial and $F_l=F_l(l)>0$ and
$0<C'=C'(C,l)<1$. This estimate immediately implies that
$\rho'(h)^{n}s$ converges in $\Sect^l(X)$  so to prove the lemma
it suffices to prove inequality (\ref{equation:keystep1}).

We let $b=2l-k+1$ and define $F_l=\sup_{\supp(h)}\|\rho'(\g)\|_b$.
 We shrink $U$ so that
$\|\rho'(\g)\|^b_1C<1$ for every $\g{\in}\supp(h)$, let
$C_h=\sup_{\supp(h)}\|\rho'(\g)\|_1$ and fix a constant $C'$ with
$\sqrt{C_h^bC}<C'<1$.  Let $f_n=\rho'(h)^{n+1}s-\rho'(h)^ns$. Then
Lemma \ref{lemma:convexityofderivatives} implies that
\begin{equation}
\label{equation:fromconvexity1}
\|f_n\|_l{\leq}B\|f_n\|^{\frac{1}{2}}_{k-1}\|f_n\|^{\frac{1}{2}}_{b}
\end{equation}
for $B$ depending only on $X$ and $b$.  Inequality
(\ref{equation:initialestimate1}) provides a bound on
$\|f_n\|_{k-1}$, so it remains to find a bound on $\|f_n\|_{b}$.
Noting that $f_n=\rho'(h)^n(\rho'(h)s-s)$ Corollary
\ref{corollary:estimateoncompostion} implies that
\begin{equation}
\label{equation:fromcomposition1} \|f_n\|_b{\leq}C_h^{nb}P(nF_l).
\end{equation}
Inequality (\ref{equation:keystep1}) is now immediate from
inequalities (\ref{equation:initialestimate1}),
(\ref{equation:fromconvexity1}) and
(\ref{equation:fromcomposition1}) and the definition of $C'$.
\end{proof}

\noindent{\bf Remark on the choice of $U$:} There are two
constraints on the choice of $U$:
\begin{enumerate}
\item $U$ is small enough so that we can apply Lemma
\ref{lemma:multiplespaces} as described in the first paragraph of
the proof for $\beta<\frac{\eta}{A}$ and some $0<C<1$  and \item
$U$ is small enough so that $\|\rho'(\g)\|_1^bC<1$ for every
$\g{\in}\supp(h)$.
\end{enumerate}
It is easy to see that we can choose $U$ to satisfy these two
conditions. An analogous remark applies to the proof of Lemma
\ref{lemma:smoothinvariantfunctionsonpairs} below.

\section{\bf Foliated results}
\label{section:foliatedproofs}

This section is devoted to the proof of Theorem
\ref{theorem:almostconjugacygen}. Though we can prove some special
cases of Theorem \ref{theorem:almostconjugacygen} by the method of
subsection \ref{subsection:isomrigid}, the general result requires
that we use the method described in subsection
\ref{subsection:isomrigid2} for isometric actions.  We begin by
recalling some facts about foliations and their holonomy
groupoids.

\subsection{\bf Holonomy groupoids and regular atlases.}
\label{subsection:functionsonpairs}

We would like to be able to apply the definitions and results of
section \ref{section:parametrizing} to the ``foliated space"
defined by taking pairs of points on the same leaf of a foliation
$\ff$ of $X$. There is a well-known difficulty in topologizing the
set of pairs of points on the same leaf as a foliated space and it
seems difficult even to make this space a measure space in a
natural way without some additional assumption on the foliation.
For product foliations $X=Y{\times}Z$ foliated by copies of $Z$,
no difficulties occur and the space is simply
$Y{\times}Z{\times}Z$. More generally, one usually considers the
holonomy groupoid or graph of the foliation, which is a, possibly
non-Hausdorff, foliated space.  To avoid technical difficulties,
we have assumed that our foliated spaces have Hausdorff holonomy
groupoids.

We now briefly describe the holonomy groupoid $P$ of the foliated
space $(X,\ff)$ in order to define group actions on $P$ associated
to $\rho$ and $\rho'$. At each point $x$ in $X$, we fix a local
transversal, $T_x$. Given a curve $c$ contained in a leaf $\fL_x$
of $\ff$, with endpoints $x$ and $y$, one can define the holonomy
$h(c)$ of $c$ as the germ of the map from $T_x$ to $T_y$ given by
moving along (parallel copies of) $c$. It is clear that $h(c)$
depends on the homotopy class of $c$. We can define an equivalence
relation on paths $c$ from $x$ to $y$ by saying two paths $c$ and
$c'$ are equivalent if $h(c)=h(c')$. Then $P$ is the set of
equivalence classes of triples $(x,y,c)$ where $c$ is a curve
joining $x$ to $y$ and two triples $(x,y,c)$ and $(x',y',c')$ are
equivalent if $x=x',y=y'$ and $h(c)=h(c')$. There is an obvious
topology on $P$ in which $P$ is a foliated space with leaves of
the form $\fL_x{\times}\tilde \fL_x$ where $\fL_x$ is a leaf of
$\ff$ and $\tilde \fL_x$ is the cover of $\fL_x$ corresponding to
homotopy classes of loops at $x$ with trivial holonomy.  When we
wish to refer explicitly to the structure of $P$ as a foliated
space, we will use the notation $(P, \tilde \ff)$. As mentioned
above, we will always assume that $P$ is Hausdorff in it's natural
topology. There are two natural projections $\pi_1$ and $\pi_2$
from $P$ to $X$ defined by $\pi_1(x,y,c)=x$ and $\pi_2(x,y,c)=y$
both of which are continuous and leafwise smooth. A transverse
invariant measure on $X$ defines one on $P$ and  a leafwise volume
form on $(X,\ff)$ defines one on $(P,\tilde \ff)$. Therefore,
under the hypotheses of Theorem \ref{theorem:almostconjugacygen},
we have a, possibly infinite, measure $\tilde \mu$ on $P$ defined
by integrating the leafwise volume form against the transverse
invariant measure. It is easy to see that $\tilde \mu=\int_X\tilde
\nu_{\ff}d\mu$ where $\tilde \nu_{\ff}$ is the pullback of the
leafwise volume form on leaves of $\ff$ to their holonomy
coverings.  For more detailed discussion, the reader should see
either \cite{CC} or \cite{MS}.

Given an action $\rho$ of $\G$ on $(X,\ff)$ defined by a
homomorphism $\rho:\G{\rightarrow}\Diff^k(X,\ff)$ we can define an
action $\rho_P$ on $P$ as follows.  Take the diagonal action of
$\rho$ on $X{\times}X$.  This defines an action of $\G$ on curves
$c$ as above, which then descends to an action on $P$. It is
immediate that if $\rho$ preserves $\mu$ then $\rho_P$ preserves
$\tilde \mu$. If $\rho'$ is a $C^k$ foliated perturbation of
$\rho$, then we can define an action $\rho'_P$ similarly, provided
$\G$ is compactly presented. We take the action on $X{\times}X$
defined by acting by $\rho$ on the first coordinate and $\rho'$ on
the second.  As long as $\rho(\g)$ is close to $\rho'(g)$, we can
define $\rho'_P(\g)$ on $P$, since there is a canonical choice of
a short, null homotopic, path from $\rho(\g)y$ to $\rho'(\g)y$
given by the length minimizing geodesic segment. Given a path $c$
from $x$ to $y$, we define
$\rho'_P(\g)(x,y,c)=(\rho(\g)x,\rho'(g)y, c')$ where $c'$ is the
concatenation of the path $\rho(g)c$ with the canonical path from
$\rho(g)y$ to $\rho'(g)y$.  Since $\G$ is compactly presented, if
$\rho$ is close enough to $\rho'$ it is easy to verify that
lifting the generating set $K$ to $P$ defines an action $\rho'_P$
of $\G$ on $P$.

It is immediate that $\pi_1:(P,\rho'_P){\rightarrow}(X,\rho)$ and
$\pi_2:(P,\rho'_P){\rightarrow}(X,\rho')$ are equivariant. Note
that compactness of $X$ implies that $\rho_P$ and $\rho'_P$ are
close in the strong topology on $\Diff^k(P, \tilde \ff)$.

When $P$ is Hausdorff, we can define a family of norms on sections
of $J^k(P)$ by $\|f\|^p_{p}=\int_X\int_{\tilde
\fL_x}\|f(x,y)\|^pd\nu_{\ff}(y)d\mu(x)$ where $\tilde
\fL_x=\pi_1{\inv}(x)$. We can complete $\Sect(J^k(P))$ to a Banach
space $L^p(J^k(P))$ of type $L^p_n$. Note that $C^k(P)$ is a
linear subspace of $\Sect(J^k(P))$ and let $L^{p,k}(P,\ff)$ be the
closure of $F^k(P)$ in $L^p(J^k(P))$.

To obtain the required estimates for Theorem
\ref{theorem:almostconjugacygen}, we will also need estimates for
the size of functions with respect to certain other $\G$ invariant
measures. Let $\lambda$ be any $\rho(\Gamma)$ invariant
probability measure on $X$. Define an norm on $\Sect(J^k(P))$ by
$\|f\|^p_{p,\lambda}=\int_X\int_{\fL_x}\|f\|^pd\nu_{\ff}d\lambda$.
We can complete $\Sect(J^k(P))$ and $F^k(P)$ with respect to this
norm to obtain Banach spaces $L^{p,\lambda}(J^k(P))$ and
$L^{p,k,\lambda}(P,\ff)$.

Except for the fact that we consider more general invariant
measures and the corresponding function spaces, the following is a
consequence of Propositions \ref{proposition:isometries} and
\ref{proposition:almostisometric} above.  The proofs of those
propositions can be repeated almost verbatim to prove this one.

\begin{proposition}
\label{proposition:almostisometriespairs} Let $\phi$ be a leafwise
isometry of $(X,\ff,g_{\ff},\mu)$.
\begin{enumerate}

\item  The maps $\phi_P$ on $L^p(J^k(P))$ and $L^{p,k}(P,\ff)$ are
isometric. Furthermore for any $\G$ invariant probability measure
$\lambda$, the maps $\phi_P$ on $L^{p,\lambda}(J^k(P))$ and
$L^{p,k,\lambda}(P,\ff)$ are isometric.

\item For any $\varepsilon>0$ and any $p_0>1$ there exists a
neighborhood $U$ of the identity in $\Diff^k(X,\ff)$ such that for
any $(U,C^k)$-foliated perturbation $\phi'$ of $\phi$, the map
$\phi_P'$ induces $\varepsilon$-almost isometries on
$L^{p,k}(P,\ff)$ and $L^p(J^k(P))$ and on $L^{p,\lambda}(J^k(P))$
and $L^{p,k,\lambda}(P,\ff)$, for any $\Gamma$ invariant
probability measure $\lambda$ on $X$ and any $p{\leq}p_0$.

\item  Let $f$ be a $\phi_P$ invariant compactly supported
function in $C^k(P)$. Then for every $\delta>0$ and every $p_0>1$
there exists a neighborhood $U$ of the identity in
$\Diff^k(X,\ff)$ such that any $(U,C^k)$-foliated perturbation
$\phi'$ of $\phi$, we have that $\|\phi_P'f-f\|_{p,k}{\leq}\delta$
and $\|\phi_P'f-f\|_{p,k,\lambda}{\leq}\delta$ for every $\phi$
invariant probability measure $\lambda$ on $X$ and every
$p{\leq}p_0$.

\item If $H$ is a topological group
 and $\rho$ is a continuous leafwise isometric action of $H$ on
 $X$, then the actions of $H$ induced by $\rho_P$ on $L^{2,k}(P,\ff)$
  and on $L^{2,k,\lambda}(P,\ff)$, for any $\G$
invariant probability measure $\lambda$ on $X$, are continuous.
Furthermore the same is true for any continuous action $\rho'$
which is an $(U,C^k)$-foliated perturbation of $\rho$.

\end{enumerate}
\end{proposition}

We need a proposition concerning covers of foliated spaces by
certain kinds of foliated charts.  This proposition follows from
the proofs that any foliation can be defined by a {\em regular
atlas} but since we require information not usually contained in
the definition of a regular atlas, we sketch the proof here.  For
more discussion of regular atlases, see sections $1.2$ and $11.2$
of \cite{CC}.  We recall that for a foliated space $(X,\ff)$ there
is an associated metric space $Y$, such that there is a basis of
foliation charts $(U,\phi)$ in $X$ of the form
$\phi:U{\rightarrow}V{\times}B(0,r)$ where $U$ is an open in $X$,
$V$ is an open set in $Y$ and $B(0,r)$ is a ball in $\Ra^n$.

\begin{proposition}
\label{proposition:regularatlas} Let $(X,\ff,g_{\ff})$ be a
compact foliated space.  Then there exists a positive number $r>0$
and a finite covering of $X$ by foliated charts $(U_i,\phi_i)$
such that:
\begin{enumerate}
\item  each $\phi_i:U_i{\rightarrow}V_i{\times}B(0,r)$ is a
homeomorphism where $U_i{\subset}X$ and $V_i{\subset}Y$ are open
and $\phi_i{\inv}:\{v_i\}{\times}B(0,r){\rightarrow}\fL{\cap}U_i$
is isometric for all $v_i{\in}V_i$ and all $i$,

\item  each $(U_i,\phi_i)$ is contained in a chart $(\tilde
U_i,\tilde \phi_i)$ such that each ${\tilde \phi_i}:{\tilde
U}_i{\rightarrow}V_i{\times}B(0,2r)$ is a homeomorphism where
$\tilde U_i{\subset}X$ and $V_i{\subset}Y$ are open and
$\phi_i{\inv}:\{v_i\}{\times}B(0,2r){\rightarrow}\fL{\cap}{\tilde
U}_i$ is isometric for all $v_i{\in}V_i$ and all $i$
\end{enumerate}
\end{proposition}

\begin{proof}
Let $\mathcal W$ be any maximal foliated atlas for $(X,\ff)$.
Since $X$ is compact, we can choose a finite cover of $X$ by
$(W_j,\psi_j)_{1{\leq}j{\leq}k}{\subset}{\mathcal W}$.  Let $\eta$
be a Lebesgue number for the cover of $X$ by $W_j$, i.e.
$B(x,\eta)$ is entirely contained in one $W_j$ for every
$x{\in}X$. Let $2d$ be the largest number so that $B(x,\eta)$
contains a chart of the form $(W_i,\psi_i)$ in $\mathcal W$ such
that $\psi_i(W_i)=V_i{\times}B(0,2d)$ and
$\psi_i|_{\{v_i\}{\times}B(0,2d)}$ is isometric for every
$v_i{\in}V_i$. Let $U_i$ be $\phi_i{\inv}(V_i{\times}B(0,d)$ and
let $\phi_i=\psi_i|_{U_i}$. Clearly the charts $(U_i,\phi_i)$
satisfy the conclusions of the proposition.  Since $X$ is compact,
we can pick a finite subset of these charts that cover $X$.
\end{proof}

\subsection{\bf Proof of Theorem \ref{theorem:almostconjugacygen}}
\label{subsection:generalfoliated}

In this subsection we prove Theorem
\ref{theorem:almostconjugacygen}.  The approach is based on the
proof of Theorem \ref{theorem:isomrigid} from subsection
\ref{subsection:isomrigid2}. In place of working on functions on
$X{\times}X$ we work with functions in the spaces $C^k(P)$ and
$L^{p,k}(P,\tilde \ff)$ defined in subsection
\ref{subsection:functionsonpairs}. By Theorem
\ref{theorem:compactlypresented}, we can assume without loss of
generality that $\G$ is compactly presented. Let $\rho$ be the
leafwise isometric action specified in Theorem
\ref{theorem:almostconjugacygen} and $\rho'$ the $C^k$ foliated
perturbation of $\rho$.  Then we have two $\G$ actions $\rho_P$
and $\rho'_P$ on $P$ as defined in the last section. As before we
will start with an invariant function $f$ for $\rho_P$, in this
case a compactly supported function in $C^k(P)$ and therefore
$L^{p,k}(P,\tilde \ff)$, and construct the desired conjugacy from
a $\rho'_P$ invariant function $f'$ close to $f$ in
$L^{p,k}(P,\ff)$.  We construct $f$ in a manner analogous to
Proposition \ref{proposition:invariantfunctions}.  We first define
a subset $\Delta{\subset}P$ which is the set of
$\{(x,x,c_x)|x{\in}X\}$ where $c_x$ is the constant loop at $x$.
Given a point $x{\in}X$, we denote by $\Delta(x)=(x,x,c_x)$. Note
that this defines canonically a point $\tilde x$ in $\tilde
\fL_x$, since the leaf of $\tilde \ff$ in $P$ through $\Delta(x)$
is $\fL_x{\times}\tilde \fL_x$ and $\tilde x$ is the projection of
$\Delta(X){\cap}(\fL_x{\times}{\tilde \fL_x})$ to $\tilde \fL_x$.
Given a point $x$ in $X$, we will refer to $B_{\ff}(x,r)$ as the
ball of radius $r$ about $x$ in the leaf through $\fL_x$ and
$B_{\tilde \ff}(\Delta(x),r)$ for the ball of radius $r$ about
$\Delta(x)$ in $\fL_x{\times}\tilde \fL_x{\subset}P$. Recall that
$N(x)$ is the radius of the largest normal ball containing $x$ in
$\fL_x$, and that $N(x)$ is bounded below by a positive number
since $X$ is compact.  We sometimes write coordinates on $P$ as
$p=(\pi_1(p),y)$ where $y{\in}\pi_1{\inv}(\pi_1(p))$.

\begin{proposition}
\label{proposition:foliatedinvariantfunctions} Let $\rho$ be a
leafwise isometric action of a group $\G$ on a compact foliated
space $(X,\ff,g_{\ff})$, then the action $\rho_P$ leaves invariant
any function on $P$ which is a function of
$d_{\fL_x}(\pi_1(p),\pi_2(p))$ or $d_{\tilde
\fL_x}(\widetilde{\pi_1(p)}, y)$. Let $r>0$ be as in Proposition
\ref{proposition:regularatlas} and such that $N(x)>2r$ for all
$x{\in}X$, then for any $0<\varepsilon{\leq}r$, there exists a
compactly supported invariant function $f{\in}C^{\infty}(P)$ such
that:
\begin{enumerate}
\item $f$ takes the value $1$ on $\Delta$,

\item $f{\geq}0$ and $f(p)=0$ if $\pi(p)=x$ and
$p{\notin}B_{\tilde \ff}(\Delta(x),2\varepsilon)$,

\item $f(p)<1/2$ if $\pi(p)=x$ and $p{\notin}B_{\tilde
\ff}(\Delta(x),\varepsilon)$

\item the Hessian of $f$ restricted to $\pi{\inv}(x)$ is negative
definite on the closure of $B_{\tilde
\ff}(\Delta(x),\varepsilon){\cap}\pi{\inv}(x)$

\end{enumerate}
\end{proposition}

\begin{proof}
The proof is very similar to the proof of Proposition
\ref{proposition:invariantfunctions}.

To define $f$ we first define $f_{x_0}:\tilde
\fL_x{\rightarrow}\Ra$ by $1-\frac{1}{2\varepsilon^2}d(x_0,x)^2$.
Our assumptions on $\varepsilon$ imply that  $B(x_0,2\varepsilon)$
is contained in a normal neighborhood of $x_0$ and so $d(x,x_0)^2$
is a smooth function of $x$ and $x_0$ inside
$B(x_0,2\varepsilon)$, see for example \cite[IV.3.6]{KN}. It is
clear that $f(x_1,x_2)=f_{\tilde x_1}(x_2)$, satisfies $1,3$ and
$4$, but it may not be smooth, fails to satisfy $2$ and may not be
compactly supported. Modifying $f_{x}$ outside $B(x,\varepsilon)$
produces a smooth,positive, compactly supported function
satisfying all the above conditions. This is easily done by
 choosing
any smooth function $g:\Ra{\rightarrow}\Ra$ such that $g$ agrees
with $1-\frac{1}{2\varepsilon^2}x^2$ to all orders for all
$x{\leq}\varepsilon$ and $g<1/2$ for all $x{\geq}\varepsilon$ and
$g=0$ for all $x{\geq}2\varepsilon$. We then let
$f(x_1,x_2)=g(d_{\tilde \fL_x}(\tilde x_1,x_2))$ where
$x_1=\pi_1(p)$ and $x_2$ is the coordinate of $p$ in $\tilde
\fL_x=\pi_1{\inv}(x_1)$.
\end{proof}

{\noindent}{\bf Remark:} We need $f$ to be compactly supported on
$\pi_1{\inv}(x)$ for all $x{\in}X$.  For this reason we choose $f$
with a global maximum along $\Delta(X)$ rather than a minimum as
in Proposition \ref{proposition:invariantfunctions}.

 One can now produce an invariant function $f'$ in
$L^{2,k}(P,\tilde \ff)$ exactly as in the proof of Theorem
\ref{theorem:isomrigid} in subsection \ref{subsection:isomrigid2}.
In order to control the behavior of $f'$ in
$L^{2,k,\lambda}(P,\tilde \ff)$ for any $\G$ invariant probability
measure $\lambda$ on $X$ as well, we need to produce $f'$ using
Theorem \ref{theorem:decreasingdisplacement} rather than Theorem
\ref{theorem:fixedpointsimplegen}.  The following is a sharpening
of Lemma \ref{lemma:multiplespaces}, and as in the proof of that
lemma, we use Theorem \ref{theorem:decreasingdisplacement} in
conjunction with Corollary
\ref{corollary:contractingoperatorbanach} and Lemma
\ref{lemma:convergingpointwiseae} to obtain optimal regularity.

\begin{lemma}
\label{lemma:perturbedinvariantfunctions} Let $\rho$ and $\G$ be
as in Theorem \ref{theorem:almostconjugacygen} and $f$ be a
compactly supported $\rho_P(\G)$ invariant function in
$C^{\infty}(P)$. Given constants $\varsigma_1>0,F>0,0<C<1$ and
$p{\geq}2$ there exists a neighborhood $U$ of the identity in
$\Diff^k(X,\ff)$ and a function $h{\in}{\mathcal{U}(\G)}$ such
that if $\rho'$ is any $(U,C^k)$-foliated perturbation of $\rho$:
\begin{enumerate}
\item $\rho_P'(h)^nf$ converges pointwise almost everywhere to
$\rho'_P$ invariant function $f'$, \item
$\|f-f'\|_{p,k}<\varsigma_1$ and
$\|f-f'\|_{p,k,\lambda}<\varsigma_1$ for every $\rho(\G)$
invariant probability measure $\lambda$ on $X$, \item
$\|\rho'_P(h)^{n+1}f-\rho'_P(h)^nf\|_{p,k,\lambda}<C^nF$ for every
$\rho(\G)$ invariant probability measure $\lambda$ on $X$.
\end{enumerate}
\end{lemma}

\begin{proof}
Given $\varsigma_1>0$ and $p{\geq}2$, we choose $\varepsilon>0$
satisfying the hypotheses of Theorem
\ref{theorem:decreasingdisplacement} and
\ref{corollary:contractingoperatorbanach} and $\delta>0$ depending
on $\varsigma_1$ to be specified below.  Choosing $U$ small
enough, Proposition \ref{proposition:almostisometriespairs}
implies that $\dk(f)<\delta$ for the $\rho'_P$ action on
$L^{p,k,\lambda}(P,\ff)$ and $L^{2,k,\lambda}(P,\ff)$ for every
$\rho(\G)$ invariant probability measure $\lambda$  and also that
$\rho'_P(k)$ is an $\varepsilon$-almost isometry on
$L^{2,k,\lambda}(J^k(P))$ and $L^{p,\lambda}(J^k(P))$ for every
$k{\in}K$ and every $\rho_P$ invariant probability measure
$\lambda$ and also that $\rho'_P$ is a continuous action on all of
these spaces.  Fix a constant $0<C<1$ and a function
$h{\in}\mathcal{U}(\G)$ so that Theorem
\ref{theorem:decreasingdisplacement} and Corollary
\ref{corollary:contractingoperatorbanach} are both satisfied for
$h$ and $C$.  Each theorem yields a constant $M_2$ and $M_p$ and
we let $M=\max(M_2,M_p)$.  Note that $M$ and $h$ depend only on
$\G,K,p$ and $C$ and some function $f{\in}\mathcal{U}(\G)$.
Therefore we can choose $\delta$ such that
$\frac{MC}{1-C}{\delta}<\varsigma_1$ and
$\delta{leq}{\frac{F}{M}}$. Then the sequence $\{\rho'_P(h^n)f\}$
satisfies
$\|\rho'_P(h^n)f-\rho'_P(h^{n-1})f\|_{p,\lambda}{\leq}MC^n\delta{\leq}C^nF$
and
$\|\rho'_P(h^n)f-\rho'_P(h^{n-1})f\|_{2,k,\lambda}{\leq}MC^n\delta{\leq}C^nF$
and $K$-displacement of $\rho'_P(h^n)f$ is less than $C^n\delta$
in $L^{2,k,\lambda}(P,\ff)$ for every $\rho(\G)$ invariant measure
$\lambda$ on $X$. By Lemma \ref{lemma:convergingpointwiseae}, this
implies that $\rho'_P(h)^nf$ converges pointwise $\lambda$ almost
everywhere to a $\rho'_P$ invariant function $f'$ which is in
$L^{p,k,\lambda}(P,\ff)$. Furthermore, our choice of $\delta$
implies that $\|f-f'\|_{p,k,\lambda}<\varsigma_1$ for every $\G$
invariant measure $\lambda$ on $X$.
\end{proof}

To obtain control over $f'$ and the resulting conjugacy on a set
$S$ as described in Theorem \ref{theorem:almostconjugacygen}, we
need to consider a certain class of $\G$ invariant measures on
$X$. Let $\mu=\int_{\mathcal{P}(X)}\mu_ed{\bar \mu}(e)$ be an
ergodic decomposition for $\mu$, where each $\mu_e$ is a $\G$
ergodic measure on $X$ and $\bar \mu$ is a measure on the space
$\mathcal{P}(X)$ of probability measures on $X$ supported on the
$\G$ ergodic measures. Let $\mathcal{P}(X)$ be the space of
regular Borel probability measures on $X$ and define a Markov
operator $M:X{\rightarrow}\mathcal{P}(X)$ be letting
$M_x=\frac{1}{\nu_{\ff}(B_{\ff}(x,r))}\nu_{\ff}|_{B_{\ff}(x,r)}$.
Then $M$ defines an operator on continuous functions on $X$ by
$Mg(x)=\int_Xg(y)dM_x(y)$ for $g$ in $C^0(X)$ and dually an
operator on $\mathcal{P}(X)$ by
$M\nu(f)=\int_XMfd\nu=\int_X\int_Xf(y)dM_x(y)d\nu(x)$ for
$\nu{\in}\mathcal{P}(X)$. Note that $M$ commutes with elements of
$\Diff^k(X,\ff)$ which are leafwise isometric.  This implies that
for any $\rho(\G)$ invariant probability  measure $\nu$ on $X$,
the probability measure $M\nu$ is also $\G$ invariant.  We will be
particularly interested in measures of the form $M\mu_e$ where
$\mu_e$ is an ergodic component of $\mu$.

\begin{lemma}
\label{lemma:closefoliatedfunctions} Let $f{\in}C^k(P)$ be
compactly supported and let $f'$ be a function which is in
$L^{p,k,\lambda}(P,\tilde \ff)$ for every $\rho(\G)$ invariant
probability measure  $\lambda$ on $X$, and such that
$\|f-f'\|_{p,k,\lambda}{\leq}A$ for every $\lambda$. Then for any
$\lambda$ and $\lambda$ almost every $x$, the restriction of $f'$
to $\pi_1{\inv}(x)$ is in $C^{k-\frac{d}{p}}$. Furthermore, there
exist constants $C$ and $r$ depending only on $(X,\ff)$ and a set
$S{\subset}X$ depending on $f'$ such that $\lambda(S)>1-\sqrt{A}$
and
$$\|(f-f')|_{\pi_1{\inv}(B_{\ff}(x,r))}\|_{k-\frac{d}{p}}{\leq}C\sqrt{A}$$
for every $x{\in}S$.
\end{lemma}

\begin{proof}
To see the first claim, we consider the measure $M\lambda$. The
fact that $\|f-f'\|_{p,k,M\lambda}<A$ implies that
$\|f'\|_{p,k,M\lambda}$ is finite.  Applying the definition of
$M$, this means that
$\int_{B_{\ff}(x,r)}\int_{\pi_1{\inv}(x)}\|j^k(f')\|^pd\nu_{\ff}d\nu_{\ff}$
is finite for $\lambda$ almost every $x$.   Then Proposition
\ref{proposition:sobolev} implies that $f'$ is $C^{k-\frac{d}{p}}$
on $\pi_1{\inv}(B_{\ff}(x,r))$.

Let $v_0=\min_{x{\in}X}(\nu_{\ff}(B_{\ff}(x,r)))$ and define $S$
to be the set of $x$ where
\begin{equation}
\label{equation:definingS} \int_{B_{\ff}(
x,d)}\int_{\pi{\inv}(x)}\|j^k(f)-j^k(f')\|^p_kd\nu_{\ff}d\nu_{\ff}{\leq}v_0\sqrt{A}.
\end{equation}
We first verify that $\lambda(S){\geq}1-\sqrt{A}$ for every
$\lambda$. We are assuming that $\|f-f'\|_{2,k,\lambda}{\leq}A$
for every $\rho$ invariant probability measure $\lambda$ on $X$.
By definition of $L^{2,k,M\lambda}(P,\tilde \ff)$ this means that
$$\int_X\int_{\pi_1{\inv}(x)}\|j^k(f)-j^k(f')\|^pd\nu_{\ff}dM(\lambda)=$$
$$\int_X\int_X\int_{\pi_1{\inv}(x)}\|j^k(f)-j^k(f')\|^pd\nu_{\ff}dM(x)d\lambda{\leq}A.$$
This implies that
$$\frac{1}{\nu_{\ff}(B_{\ff}(x,r))}\int_{B_{\ff}(x,r)}\int_{\pi_1{\inv}(x)}\|j^k(f)-j^k(f')\|^pd\nu_{\ff}d\nu_{\ff}<\sqrt{A}$$
or
$$\int_{B_{\ff}(
x,r)}\int_{\pi_1{\inv}(x)}\|j^k(f)-j^k(f')\|^p_kd\nu_{\ff}d\nu_{\ff}{\leq}\sqrt{A}\nu_{\ff}(B_{\ff}(x,r))$$
on a set of $\lambda$ measure at least $1-\sqrt{A}$. Therefore the
set $S$ defined by equation \ref{equation:definingS} has $\lambda$
measure at least $1-\sqrt{A}$ as desired.

Proposition \ref{proposition:sobolev} then implies that
$$\|(f-f')|_{\pi_1{\inv}(B_{\ff}(x,r))}\|_{k-\frac{d}{p}}<C'v_0\sqrt{A}$$
for every $x{\in}S$, where $C$ is a constant depending only on
$X,\ff$ and $g_{\ff}$ and letting $C=C'v_0$ completes the proof.
\end{proof}

Fix $p$ with $\frac{\dim(\fL_x)}{p}<\kappa$ and fix a function $f$
as in the conclusion of Proposition
\ref{proposition:foliatedinvariantfunctions} with
$\varepsilon=r/2$ for the remainder of this section.  We choose a
constant $\varsigma_1$ to be specified below and let $f'$ be the
function produced by Lemma
\ref{lemma:perturbedinvariantfunctions}. Then Lemma
\ref{lemma:closefoliatedfunctions} combined with the definition of
$f$ and $S$ implies that for every $x{\in}S$:
\begin{enumerate}
\item $f'(p)<C\sqrt{\varsigma_1}$ if $\pi_1(p)=x$ and
$p{\notin}B_{\tilde \ff}(\Delta(x),r)$,

\item $f'(p)<1/2+C\sqrt{\varsigma_1}$ if $\pi_1(p)=x$ and
$p{\notin}B_{\tilde \ff}(\Delta(x),\frac{r}{2})$

\item for every $y{\in}B_{\ff}(x,r)$, the Hessian of $f'$
restricted to $\pi_1{\inv}(y)$ is negative definite on
$B(\Delta(y),\frac{r}{2}){\cap}\pi_1{\inv}(y)$
\end{enumerate}

Choosing $\varsigma_1{\leq}\frac{1}{100C^2}$ so
$C\sqrt{\varsigma_1}{\leq}\frac{1}{10}$ this implies that for
$x{\in}X$, the function $f'$ has a maximum on $\pi_1{\inv}(x)$ at
a point $\tilde \phi(x)$ where the value is at least
$\frac{9}{10}$ and that this maximum is the only local maximum
with value greater than $\frac{6}{10}$. Since $f'$ is invariant
under $\rho'_P$ it follows that if $\rho(\g)(x){\in}S$ then
$\tilde \phi(\rho(g)(x))=\rho'(g)\tilde \phi(x)$ since both points
will be the global maxima of $f'$ on $\pi_1{\inv}(\rho(g)x)$.
Furthermore, it follows by $\rho'_P(\G)$ invariance of $f'$ that
for every $x{\in}\G{\cdot}S$, there is a unique global maximum for
$f'|_{\pi_1{\inv}(x)}$.   Therefore we can define the conjugacy
between $\rho$ and $\rho'$ on a set of full measure in $X$ by
letting $\tilde \phi(x)$ be the unique global maximum of $f'$ on
the fiber $\pi_1{\inv}(x)=\tilde \fL_x$ and letting
$\phi(x)=\pi_2(\tilde \phi(x))$.

 We remark that it is possible to show that
$d_{\fL_x}(x,\phi(x))^2{\leq}{C'\sqrt{\varsigma_1}r^2}$ for all
$x{\in}S$ directly from the definition of $f,f'$ and $\phi$, where
$C'=\frac{C}{v_0}$ is as in the proof of Lemma
\ref{lemma:closefoliatedfunctions}.

In order to make the following more readable, we let
$k'=k-\frac{d}{p}-1$.  We now show the map $\phi$ is leafwise
$C^{k'}$ for $x{\in}S$, and therefore, by equivariance, leafwise
$C^{k'}$ almost everywhere. Consider $x{\in}S$. We will show that
$\phi$ is $C^{k'}$ and $C^{k'}$ close to the identity on
$B_{\ff}(x,r)$. Note that $\pi_1{\inv}(B_{\ff}(x,r)))$ is
diffeomorphic to $B_{\ff}(x,r){\times}\tilde \fL_x$.  Let
$D_2f':B_{\ff}(x,r){\times}T(\tilde \fL_x){\rightarrow}\Ra$ be the
derivative of $f'$ in the second variable.  Let
$N(\frac{r}{2},B_{\ff}(x,r))$ be the $\frac{r}{2}$ neighborhood of
$\Delta(B_{\ff}(x,r))$ in $B_{\ff}(x,r){\times}\tilde \fL_x$ and
$TN(\frac{r}{2},B_{\ff}(x,r))$ the restriction of the bundle
$\fL_x{\times}T\tilde \fL_x$ to that set. If $x{\in}S$ then the
set of points $(x,\tilde \phi(x))$ is
$D_2f'{\inv}(0){\cap}TN(\frac{r}{2},B_{\ff}(x,r))$ and $0$ is a
regular value of $D_2f'$, since the Hessian is negative definite
on $\tilde \fL_{\tilde y}{\cap}B(\tilde y,\frac{r}{2})$ for every
$y{\in}B_{\ff}(x,r)$. This implies that the set $(x,\tilde
\phi(x)){\subset}N(\frac{r}{2},B_{\ff}(x,r)){\subset}TN(\frac{r}{2},B_{\ff}(x,r))$
is a $C^{k'}$ submanifold and so $\tilde \phi$ is $C^{k'}$ on
$B_{\ff}(x,r)$. This implies that $\phi$ is $C^{k'}$ on
$B_{\ff}(x,r)$.  Since $f'$ is $C^{k'}$ close to $f$ on
$\pi_1{\inv}(B_{\ff}(x,r))$, the functions $D_2f$ and $D_2f'$ are
$C^{k'}$ close on $TN(\frac{r}{2},B_{\ff}(x,r))$. This implies
that the submanifolds $D_2f{\inv}(0)$ and $D_2f'{\inv}(0)$ are
$C^{k'}$ close, which then implies that $\phi$ is $C^{k'}$ close
to the identity on $B_{\ff}(x,r)$. More precisely by choosing
$\varsigma_1$ small enough, we can assume that the $C^{k'}$ norm
of $\phi-\Id:B_{\ff}(x,r){\rightarrow}B_{\ff}(x,2r)$ is less than
$\varsigma$.  We let $\bar \varsigma_1$ be the value of
$\varsigma_1$ required for this estimate and let
$\varsigma_1=\min(\bar \varsigma_1,
\frac{1}{100C^2v_0^2},\sqrt{\varsigma})$.  We let
$U_0{\subset}\Diff^k(X,\ff)$ be the neighborhood of the identity
such that for any $(U_0,C^k)$-foliated perturbation $\rho'$ of
$\rho$, the function $f'$ produced by Proposition
\ref{proposition:perturbedinvariantfunctions} satisfies
$\|f-f'\|_{p,k,\lambda}{\leq}\varsigma_1$ for every $\rho'(\G)$
invariant measure $\lambda$ on $X$.  Then we have verified
conclusions $(1),(2)$ and $(3)$ of Theorem
\ref{theorem:almostconjugacygen} for any $(U_0,C^k)$-foliated
perturbation $\rho'$ of $\rho$.

To show the final estimate in the statement of Theorem
\ref{theorem:almostconjugacygen}, we need the following
(well-known) quantitative refinement of \cite[III.5.12]{M}.

\begin{lemma}
\label{lemma:spreadingsets} There exists a constant $0<t<1$,
depending only on $\G$ and $K$,  such that for any
$0<\eta<\frac{1}{2}$ and any ergodic action of $\G$ on a finite
measure space $(X,\mu)$ and any set $S$ of measure $1-\eta$, there
is $k$ in $K$ such that $\mu((kS{\cup}S)^c){\leq}t\eta$.
\end{lemma}

\begin{proof}
Assume not.  Then for all $t$ with $0<t<1$, there exists $S$ of
measure $1-\eta$ such that $\mu(kS{\cup}S){\leq}1-t\eta$ for all
$k{\in}K$. We will use this fact to show that the characteristic
function $\chi_S$ has $K$-displacement $(1-t)\eta$ and use this to
produce a $\G$ invariant function which is closer to the
characteristic function of $S$ than any constant function.

Since $\mu(kS)+\mu(S)-\mu(kS{\cap}S)=\mu(kS{\cup}S)$ and
$\mu(S)=\mu(kS)=1-\eta$, we have
$\mu(kS{\cap}S){\geq}(1-\eta)-(1-t)\eta$.  Therefore
$\dk(\chi_S){\leq}\sqrt{(1-t)2\eta}$ in $L^2(X,\mu)$. By the
standard linear analogue of Theorem
\ref{theorem:fixedpointsimplegen} (which is an easy consequence of
Lemma \ref{proposition:decreasingdisplacementrep}) there is a
constant $C$ depending only on $\G$ and $K$ and a $\G$ invariant
function within $C\sqrt{(1-t)2\eta}$ of $\chi_S$.

Since the orthogonal complement of the constant functions are the
functions of integral zero, the distance from $\chi_S$ to the
constant functions is $\sqrt{\eta(1-\eta)}$. Since
$1-\eta>\frac{1}{2}$, we have a contradiction provided
$C\sqrt{(1-t)2\eta}<\frac{\sqrt{\eta}}{\sqrt{2}}$ or
$2C\sqrt{1-t}<1$.  So for $t>1-\frac{1}{4C^2}$ we are done.
\end{proof}

Since $\mu_e(S){\geq}1-\varsigma$ for some $\varsigma>0$ for
almost every ergodic component $\mu_e$ of $\mu$, it follows from
Lemma \ref{lemma:spreadingsets} that the measure of
$S_n=K^n{\cdot}S{\cup}{\cdots}K{\cdot}S{\cup}S$ is at least
$1-t^n\varsigma$.

Choose a neighborhood of the identity $U_1{\subset}\Diff^k(X,\ff)$
such that for every $x{\in}X$ and every $\g{\in}K$ we have
$\|j^{k}(\rho'(\g))(x)\|{\leq}(1+\varsigma)$ for any
$(U_1,C^k)$-foliated perturbation $\rho'$ of $\rho$. Let
$U=U_1{\cap}U_0$ and let $\rho'$ be $(U,C^k)$-foliated
perturbation of $\rho$. Then the fact that
$\|j^{k'}(\phi)(x)\|<1+\varsigma$ for every $x{\in}S$, combined
with the chain rule, the definition of $U_0$ and the fact that
$\mu(S_n^c)<t^n\varsigma$, imply conclusion $(4)$ of the theorem.

The remaining claim of the theorem states that given a positive
integer $l$ then, if $U$ is small enough, $\phi$ is $C^l$.  To
prove this claim, it clearly suffices to see that $f'$ is
$C^{l+1}$ on $\pi_1{\inv}(B_{\ff}(x,r))$ for almost every $x$.
This is exactly the content of Lemma
\ref{lemma:smoothinvariantfunctionsonpairs} in the next
subsection.

\subsection{Improving regularity of $\varphi$ in Theorem
\ref{theorem:almostconjugacygen}}
\label{subsection:smootherinvariantfunctionsonpairs}

We retain all notations and conventions from the previous two
subsections.  As remarked at the end of the last subsection, to
complete the proof of Theorem \ref{theorem:almostconjugacygen}, it
suffices to prove that $f'$ is $C^{l+1}$ on
$\pi_1{\inv}(B_{\ff}(x,r))$ for almost every $x$ in $X$.  This is
exactly the content of the following lemma.

\begin{lemma}
\label{lemma:smoothinvariantfunctionsonpairs} Let $f$ in $C^k(P)$
be a compactly supported, $\rho_P(\G)$ invariant function. For any
$k{\geq}3$, given a positive integer $l{\geq}k$, there exists
$U{\in}\Diff^{k}(X,\ff)$ such that for $\rho'$ a
$(U,C^{k})$-foliated perturbation of $\rho$ defined by a map from
$\G$ to $\Diff^{2l-k+1}(X,\ff)$, the sequence $\{\rho'_P(h)^nf\}$
converges pointwise almost everywhere to a (measurable) function
$f'$ on $P$ such that for almost every $x{\in}X$, the restriction
of $f'$ to $\pi_1{\inv}(B_{\ff}(x,r))$ is $C^{l}$.
\end{lemma}

To simplify the argument, we will use the operator $M$ on
$\mathcal{P}(X)$ defined in subsection
\ref{subsection:generalfoliated} and consider the measure $M\mu$.

We also introduce another technical mechanism to simplify the
proof.  Given $h{\in}\mathcal{U}(\G)$, we can define a measure on
$\G$ by $h{\mu_{\G}}$ where $\mu_{\G}$ is Haar measure on $\G$. We
can then define a probability space $\Omega=\prod_{\Za}\G$ with
measure $\lambda=\prod_{\Za}h{\mu_{\G}}$ and the left shift is an
invertible measure preserving transformation $T$ of
$(\Omega,\lambda)$.  For any measure preserving action $\sigma$ of
$\G$ on a space $Y$, we can define a skew product extension by
$T_{\sigma}(\omega,y)=(T(\omega), \sigma(\omega_0)y)$ and
$T_{\sigma}{\inv}(\omega,y)=(T{\inv}(\omega),
\sigma(\omega_{-1}){\inv}y)$. Identifying functions on $Y$, or
more generally, sections of bundles over $Y$, with their pullbacks
to $\Omega{\times}Y$, it is clear that
$\int_{\Omega}T_{\sigma}fd\lambda=\rho(h)f$ for every function $f$
on $Y$.

Before proving Lemma \ref{lemma:smoothinvariantfunctionsonpairs},
we state the variant of Corollary
\ref{corollary:estimateoncompostion} needed here. Since $(P,\ff)$
is a foliated space, we can use the definitions of norms on
$\Diff^k(X,\ff)$ from section \ref{section:convexityofderivatives}
to define norms on $\Diff^k(P,\ff)$ and the estimates from Lemma
\ref{lemma:estimateoncompostion} clearly hold for maps of $P$ as
well.   If $\psi(x,y)=(\phi_1(x),\phi_2(y))$, then it follows from
the definitions that $\|\psi\|_k=\max(\|\phi_1\|_k,\|\phi_2\|_k)$.
If $\phi_1$ is a $C^k$ leafwise isometry and $\phi_2$ is
$(U,C^k)$-foliated perturbation of $\phi_1$, and we let
$\psi=(\phi_1,\phi_2)$, then $\|\psi\|_k=\|\phi_2\|_k$.  Similarly
for $h{\in}\mathcal{U}(\G)$ we can define the operator norm of
$\rho'_P(h)$ acting on $J^k(P,\ff)$ which we denote by
$\|\rho'_P(h)\|_k$.

\begin{corollary}
\label{corollary:estimateoncompostionfoliated} Under the
assumptions of Lemma \ref{lemma:smoothinvariantfunctionsonpairs},
for any function $h{\in}\mathcal{U}(\G)$ we have the following
estimate:
$$\|\rho'_P(h)^n\|_k{\leq}N_1^{kn}Q(nN_k)$$
where $Q$ is the same polynomial as in Lemma
\ref{lemma:estimateoncompostion} above and
$N_i=\max_{\supp(h)}\|\rho'(\g)\|_i$.
\end{corollary}

\noindent{\bf Remarks:}\begin{enumerate} \item The proof is
identical to the proof of Corollary
\ref{corollary:estimateoncompostion}, so we omit it. \item The
fact that we need only consider $\|\rho'(\g)\|_i$ and not
$\|\rho'_P(\g)\|_i$ in the statement of the corollary follows from
the fact that $\rho'$ is a $(U,C^k)$ foliated perturbation of
$\rho$ which implies $\|\rho'g)\|_i=\|\rho'_P(\g)\|_i$. \item The
need for this estimate is explained following Corollary
\ref{corollary:estimateoncompostion}.
\end{enumerate}

\begin{proof}[Proof of Lemma
\ref{lemma:smoothinvariantfunctionsonpairs}] Fix $p$ such that
$\frac{2\dim(\fL_x)}{p}<1$.   For a choice of $0<C<1$ and any
choice of $\varsigma_1>0$ and $F>0$, choose a neighborhood $U$ in
$\Diff^K(X,\ff)$ and a function $h$ in $\mathcal{U}(\G)$
satisfying Proposition
\ref{proposition:perturbedinvariantfunctions}.  Then for any
$(U,C^k)$-foliated perturbation $\rho'$ of $\rho$, we have that
$\rho_P'(h)^nf$ converges pointwise almost everywhere to $\rho'_P$
invariant function $f'$ with respect to $M\mu$ and that:
$$\|\rho'_P(h)^{n+1}f-\rho_P'(h)^nf\|_{p,k,M\mu}{\leq}C^nF.$$

 Let
$0<D=\sqrt{C}<1$ and applying Lemma
\ref{lemma:closefoliatedfunctions} shows that there exists a set
$S_n$ such that $\mu(S^c_n)<D^{n}\delta$ and for every point
$x{\in}S_n$, we have that
\begin{equation}
\label{equation:initialestimate2}
\|(\rho'_P(h)^{n+1}f-\rho'_P(h)^nf)|_{B_{\ff}(x,r)}\|_{k-1}{\leq}D^nAF
\end{equation}
where $A>0$ is an absolute constant depending only on
$(X,\ff,g_{\ff})$.

We will show that, possibly after shrinking $U$,  $\rho'_P(h)^nf$
satisfies
\begin{equation}
\label{equation:keystep2}
\|(\rho'_P(h)^{n+1}f-\rho'_P(h)^nf)|_{\pi_1{\inv}(B_{\ff}(x,r))}\|_{l}{\leq}D'^nP(nF_l)F
\end{equation}
for $n>j(x)$ where $j$ is an integer valued measurable function on
$X$, where $P$ is a fixed polynomial, and $F_l>0$ and
$0<D'=D'(D,l,h)<1$. This estimate immediately implies that
$\rho'_P(h)^{n}f|_{\pi_1{\inv}(B_{\ff}(x,r))}$ converges in
$C^l(\pi_1{\inv}(B_{\ff}(x,r)))$ which suffices to complete the
proof.

We let $b=2l-k+1$ and can now define
$F_l=\sup_{\supp(h)}\|\rho'(\g)\|_b$.  We shrink $U$ so that
$\|\rho'(\g)\|^b_1D<1$ for every $\g{\in}\supp(h)$, let
$D_h=\sup_{\supp(h)}\|\rho'(\g)\|_1$. We also fix the constant
$D'=D'(l,h,D)$ such that $\sqrt{D_h^bD}<D'<1$.  Letting
$f_n=\rho'_P(h)^{n}(\rho'_P(h)f-f)$ and
$f_n^x=f_n|_{B_{\ff}(x,r)}$, Lemma
\ref{lemma:convexityofderivatives} implies that
\begin{equation}
\label{equation:fromconvexity2}
\|f^x_n\|_l{\leq}B\|f^x_n\|^{\frac{1}{2}}_{k-1}\|f^x_n\|^{\frac{1}{2}}_{b}
\end{equation}
for $B$ depending only on $X,\ff,b$.

We now form the product $\Omega{\times}X$ with measure
$\mu{\times}{\lambda}$ and transformation $T_{\rho'}$ as defined
in the paragraph immediately preceding the proof.  Define subsets
$\tilde S_n=\Omega{\times}S_n$. We now define sets $\bar
S_j=\cap_{i=j+1}^{\infty}T^{-i}\tilde S_i$. This is the set of
$(\omega,x){\in}X$ such that $T_{\rho'}^i(\omega,x){\in}\tilde
S_i$ for all $i{>}j$. The Borel-Cantelli lemma implies that
$\cup\bar S_j$ has full measure in $\Omega{\times}X$. The function
$j$ will be defined so that $j(x)$ is the smallest integer such
that $x{\in}\pi_X(\bar S_j)$ and we will prove inequality
(\ref{equation:keystep2}) by fixing $j$ and assuming
$x{\in}\pi_X(\bar S_j)$.

Applying inequality (\ref{equation:initialestimate2}) to any point
in $\pi_X{\inv}(x){\cap}{\bar S_j}$ implies that
\begin{equation}
\label{equation:initialestimatemodified}
\|f_n^x\|_{k-1}{\leq}D^{n}AF
\end{equation}
for every $x$ with $x{\in}\pi_X(\bar S_j)$ whenever $n>j$. It
remains to find a bound on $\|f^x_n\|_{b}$. Noting that
$f_n=\rho'_P(h)^n(\rho'_P(h)f-f)$ Corollary
\ref{corollary:estimateoncompostion} implies that
\begin{equation}
\label{equation:fromcomposition2} \|f^x_n\|_b{\leq}D_h^{nb}P(nF_l)
\end{equation}
where $P$ is a constant multiple of the polynomial occurring in
Corollary \ref{corollary:estimateoncompostion}. Inequality
(\ref{equation:keystep2}) is now immediate from inequalities
(\ref{equation:fromconvexity2}),
(\ref{equation:initialestimatemodified}) and
(\ref{equation:fromcomposition2}) and the definition of $D'$.
\end{proof}

\appendix
\section{``Good spaces" for continuous limit actions}
\label{appendix:limits}

The purpose of this appendix is to show how to adapt the argument
given in subsection \ref{subsection:Tproofs} to prove the general
cases of Proposition \ref{proposition:continuoussubactionHilbert}
and Proposition
\ref{proposition:convolutionalmostisometrichilbert}.  The proof of
Proposition \ref{proposition:continuoussubactionHilbert} is
completed in the first two subsections and the third subsection
ends with the proof of Proposition
\ref{proposition:convolutionalmostisometrichilbert}. More
generally, this appendix contains a series of remarks concerning
the category of spaces and actions which admit ''good" limit
actions, as well as characterizations of certain of these spaces.

\subsection{Triangles and convexity of continuous subactions}
\label{subsection:appendix}

In this subsection we outline a proof that, under the hypotheses
of Proposition \ref{proposition:continuoussubactionHilbert}
$\rho=\omega$-$\lim\rho_n$ is continuous on an affine subspace. To
see this it suffices to study sequences of triples
$A_n,B_n,C_n{\in}\fh_n$ such that $C_n=tA_n+(1-t)B_n$ and show
that equicontinuity at $A_n$ and $B_n$ implies equicontinuity at
$C_n$. This follows from the fact that almost isometries are
almost affine, i.e. that the image of a convex combination of
points under an almost isometry is close to the same convex
combination of the images of the points. We state this fact
precisely only for globally defined actions, though a more
complicated analogue is clearly true for partially defined
actions.

\begin{lemma}
\label{lemma:almostaffine} For every $\eta>0, t_0>0$ and $R>0$
there exists $\varepsilon>0$ such that if $f$ is an
$\varepsilon$-almost isometry of a Hilbert space $\fh$ and
$A,B{\in}\fh$ with $d(A,B)<R$ and $t<t_0$ then
$$d(f(tA+(1-t)B),tf(A)+(1-t)f(B))<\eta$$
\end{lemma}

As the lemma is easily proved from elementary facts concerning
stability of triples of collinear points in a Euclidean space, we
only indicate what is needed for the proof. Take three collinear
points $A,B,C$ and three arbitrary points $A',B',C'$ such that
$d(A,B)\simeq d(A',B'), d(A,C)\simeq d(A',C')$ and $d(B,C) \simeq
d(B',C')$. Then the triangles $\Delta(ABC)$ and $\Delta(A'B'C')$
are almost congruent.  More precisely, if we move $\Delta(A'B'C')$
by an isometry so that $A=A'$ and so that $B'$ is as close as
possible to $B$, then $C$ will be close to $C'$.  We leave precise
quantification of this fact to the interested reader. As Lemma
\ref{lemma:almostaffine} uses only this fact about triangles, it
is clear that the lemma is true much more generally. For example,
the lemma holds for any $L^p$-type space where $1<p<\infty$, as
well as for $CAT(0)$ spaces. The lemma, and therefore the first
conclusion of Proposition
\ref{proposition:continuoussubactionHilbert}, should hold for any
geodesic metric space which is uniformly convex in any reasonable
sense, see below or \cite{KM} for possible definitions.  If the
space does not admit a linear structure, one needs to interpret
affine subspaces and affine combinations in terms of the geodesic
structure.

\subsection{Barycenters, uniform convexity and almost isometries}
\label{subsection:barycenter} In this subsection we indicate the
proof of the remaining conclusion of Proposition
\ref{proposition:continuoussubactionHilbert}. Given a metric space
$Y$, let $\mathcal{P}(Y)$ be the set of regular, Borel probability
measures on $Y$. Given $\mu{\in}\mathcal{P}(Y)$, we define
$f_{\mu}(x)=\int_Yd(y,x)^2d\mu(y)$.  If $f_{\mu}$ attains a global
minimum at a unique point, we call that point the {\em barycenter}
of the measure, and we denote by $b:\mathcal{P}(Y){\rightarrow}Y$
the map taking a measure to its barycenter (when it exists).  For
any point $x_0$ in a Hilbert space $\fh$, the function
$f_{x_0}=d(x_0,x)^2$ has the property that it's restriction to any
geodesic has second derivative $2$ at every point.  By definition
this property is inherited by $f_{\mu}$ for any measure $\mu$.
This implies that $f_{\mu}$ has at most one minimum and easily
implies that the barycenter is defined at least when the support
of $\mu$ is compact.  The barycenter is not defined for $\mu$ with
non-compact support as can be seen by taking an atomic measure
supported on an infinite sequence of points $\{x_n\}$ which go to
infinity much faster than $\mu(x_n)$ goes to zero.  More
generally, barycenters will exist for measures which decay fast
enough at infinity.  We leave the precise formulation to the
reader.

The relevance of this discussion for subsection
\ref{subsection:Tproofs} follows from the fact that for Hilbert
spaces $b(\mu)=\int_{\fh}vd\mu(v)$. This is easily seen by showing
that $\int_{\fh}vd\mu(v)$ is a critical point for $f_{\mu}$.
Combined with our observation on the second derivative of
$f_{\mu}$ along any geodesic this implies:

\begin{lemma}
\label{lemma:strictconvexity} For every Hilbert space $\fh$ and
every compactly supported $\mu{\in}\mathcal{P}(Y)$, the function
$f_{\mu}$ has a unique global minimum $m_{\mu}$ at a point
$y_{\mu}=\int_{\fh}v d\mu(v)$. Furthermore, for every
$\varepsilon>0$ and any compactly supported probability measure
$\mu$ on any Hilbert space $\fh$ the set of points where
$f_{\mu}(x)<m_{\mu}+\varepsilon$ is contained in
$B(b(\mu),{\sqrt{\varepsilon}})$.
\end{lemma}

It is immediate from the definition that $b$ is $\Isom(\fh)$
equivariant.  We now describe the behavior of $b$ under
$\varepsilon$-almost isometries.

\begin{lemma}
\label{lemma:barycentervarepsilonisometry} For every
$D,\varepsilon>0$ there is an $\eta>0$ such that if $\fh$ is a
Hilbert space, $y_0{\in}\fh$ is a  basepoint and
$\mu{\in}\mathcal{P}(\fh)$ with $m_{\mu}<D$ and
$\supp(\mu){\in}B(y_0,R)$ and $g$ is a $\eta$-almost isometry from
$B(y_0,R)$ to $Y$, then
$$d(g(b(\mu)),b(g_*\mu))<\varepsilon.$$
\end{lemma}

\begin{proof}
Since $g$ is an $\eta$-almost isometry, we know that
$$(1-\eta)d(x,y){\leq}d(g(x),g(y)){\leq}(1+\eta)d(x,y)$$
for every $x,y{\in}B(y_0,R)$.  Since $g_*\mu(S)=\mu(g{\inv}(S))$,
squaring and integrating implies that
$(1-\eta)^2f_{\mu}(y){\leq}f_{g_*\mu}(g(y)){\leq}(1+\eta)^2f_{\mu}(y)$.
In particular
$(1-\eta)^2m_{\mu}{\leq}m_{g_*\mu}{\leq}(1+\eta)^2m_{\mu}$ and
therefore $f_{\mu}(b(g_*\mu)){\leq}(1+\eta)^4m_{\mu}$. Combined
with Lemma \ref{lemma:strictconvexity} this implies that
$$d(g(b(\mu)),b(g_*\mu))<\sqrt{((1+\eta)^4-1)D}.$$
\end{proof}

To complete the proof of Proposition
\ref{proposition:continuoussubactionHilbert} it suffices to show
that $\{y_n\}$ is in $C$ whenever $y_n=\rho_n(f)z_n$ for
$z_n{\in}\tilde X$ and $f{\in}\mathcal{U}(\G)$.  Letting $\mu_n$
be the push-forward of $fd\mu_\G$ under $\rho_n^{z_n}$ we need to
show that for every $\varepsilon>0$, there exists a neighborhood
of the identity $U$ in $\G$ such that  $d_n(\g
y'_n,y'_n)<\varepsilon$ for every $\g{\in}\G$. Note that $d_n(\g
y'_n,y'_n){\leq}d_n(\g y'_n,b(\g\mu_n))+d(b(\g\mu_n) ,b(\mu_n))$.
The first term can be made arbitrary small by Lemma
\ref{lemma:barycentervarepsilonisometry} since
$\omega$-$\lim{\varepsilon_n}=0$.  Bounding the second term
follows as in the proof of the affine case of Proposition
\ref{proposition:convolutionalmostisometrichilbert}.

\noindent {\bf Remarks:}
\begin{enumerate}
\item
 We can define
the {\em lower second derivative} of a function
$f:\Ra{\rightarrow}\Ra$ by
$$\underline{f}''(x)=\liminf_{h{\rightarrow}0}\frac{f(x+h)+f(x-h)-2f(x)}{h^2}.$$
For any $\CAT(0)$ space $Y$ and any point $y_0$, it is easy to
show that the restriction of $f_{y_0}(y)=d(y,y_0)^2$ to any
geodesic satisfies $\underline{f}_{y_0}''{\geq}2$.  Only slightly
more difficult is showing that this property characterizes
$\CAT(0)$ spaces. A similar remark is made in \cite{Gr2}.

\item A harder exercise is to show that if for every point
$y_0{\in}Y$ and every geodesic $c$ in $Y$, we have
$\underline{f}_{y_0}''=2$ on $c$, then $Y$ is Hilbert space.

\item An analog of Lemma \ref{lemma:barycentervarepsilonisometry},
and therefore Proposition
\ref{proposition:continuoussubactionHilbert}, is true for more
general spaces $X$ in place of the Hilbert space $\fh$, provided
we define $\rho(h)x=b(\nu)$ where $\nu$ is push-forward of
$h\mu_{\Gamma}$ under the orbit map $\rho(\G)x{\rightarrow}X$. In
particular, if all spaces acted upon are $\CAT(0)$ spaces or are
$L^{p_n}(Y,\nu)$ for $(Y,\nu)$ a standard measure space and
$1<\omega$-$\lim{p_n}<\infty$.  More generally, this will be true
for any uniformly convex metric space in the sense of say
\cite{KM}.

\end{enumerate}

\subsection{Convolutions and linear structure}
\label{subsection:convolution}

 We will now proceed to prove Proposition
\ref{proposition:convolutionalmostisometrichilbert} from the
following lemma and Lemma
\ref{lemma:barycentervarepsilonisometry}.

\begin{lemma}
\label{lemma:affinecombination} Let $X=\{x_1,{\ldots},x_n\}$ and
$Y=\{y_1,{\ldots},y_n\}$ be finite subsets of a Hilbert space
$\fh$ such that $d(x_i,y_i)<\eta$ for all $i$.  Then for any
coefficients $a_1,{\ldots},a_n{\in}{\mathbb R}$, we have
$d(\sum_{i=1}^n{a_i}{x_i},\sum_{i=1}^n{a_i}{y_i})<(\sum_{i=1}^n|a_i|){\eta}$.
\end{lemma}

\begin{proof}
All statements are easy consequences of the triangle inequality
for the norm on the Hilbert space and the fact that
$d(u,v)=\|u-v\|$.
\end{proof}

\begin{proof}[Proof of Proposition \ref{proposition:convolutionalmostisometrichilbert}]
Since atomic measures with finite support are dense in
$\mathcal{P}(\G)$ we assume that $\mu=\sum_ia_i\delta_{\g_i}$
where the $g_i{\in}G$ and $a_i$ are positive reals. First we note
that $\rho(\mu*\lambda)x=\sum_ia_ib(g_i\rho^x_*\lambda)$ and
$\rho(\lambda)x=b(\rho^x_*\lambda)$ where
$\rho^x:G{\rightarrow}G{\cdot}x$ is the orbit map. Now
$\rho(\mu)\rho(\lambda)x=\sum_i a_ig_ib(\rho^x_*\lambda)$. By
Lemma \ref{lemma:barycentervarepsilonisometry} applied to each
$g_i$ and the measure $\rho^x_*{\lambda}$ we have that
$d(g_ib(\rho^x_*\lambda),b(g_i\rho^x_*\lambda)){\leq}{\eta}$. The
Proposition now follows from Lemma \ref{lemma:affinecombination}.
\end{proof}

Unlike Proposition \ref{proposition:continuoussubactionHilbert},
Proposition \ref{proposition:convolutionalmostisometrichilbert}
holds in much less generality, since it depends on the affine
structure of $X_n$ and the equation $b(\mu)=\int_{X_n}vd{\mu}$. In
fact, to prove more general variants of our results it is probably
best to simply define $\rho(h)x=b(\rho^x_*h\mu_{\G})$ and work
with this averaging operator instead of the linear one.

\section{Estimates on Compositions}
\label{appendix:compositions}

 This appendix contains a proof of
Lemma \ref{lemma:estimateoncompostion}. Given the definitions, it
suffices to prove the Lemma for $\phi{\in}\Diff^k(X,\ff)$.   We
deduce this from some elementary facts about block upper
triangular matrix.

Given a number $N$, we consider $N{\times}N$ matrices which are
{\em block upper triangular}.  By this we mean that there are
number $i_1,{\ldots}i_n$ such that $\sum_{l=1}^ni_l=N$ and the
matrices $M$ have $i_l{\times}i_l$ blocks, which we denote
$A_{l}$, along the diagonal, are zero below these blocks, and have
arbitrary entries above them.  We call such $M$ {\em block upper
triangular of type $i_1, {\cdots}, i_n$}.  We define a norm on
matrices by taking the maximum of the matrix coefficients.  It is
easy to see that this is equivalent to the operator norm.

\begin{lemma}
\label{lemma:blockdiagonal} Let $M_1,{\ldots},M_j$ be block upper
triangular matrices of type $i_1,{\ldots},i_n$.  Assume that
$\|A_l\|<C_1$ for all $l$ and that all other entries of each $M_k$
are bounded by $C_2$.   Then there exists a polynomial $Q$
depending on the type of the $M_l$ such that
$$\|M_1M_2{\cdots}M_j\|{\leq}C_1^jQ(jC_2).$$
\end{lemma}

\begin{proof}
It is easy to see that the diagonal blocks of
$M=M_1M_2{\cdots}M_j$ satisfy this bound, and in fact are less
than $C_1^j$.  For any coefficient of the product outside of the
diagonal blocks, we can write the matrix coefficients of $M$ as:
$$M_{\alpha,\beta}=\sum_{\alpha=\eta_0{\leq}\eta_1{\cdots}{\leq}\eta_{j-1}{\leq}{\beta}=\eta_j}(M_1)_{\alpha,\eta_1}
(M_2)_{\eta_1,\eta_2}{\ldots}(M_j)_{\eta_{j-1},\beta}.$$ It is
easy to see that at most $N$ of the $(M_l)_{\eta_{l-1},\eta_{l}}$
can be outside the diagonal blocks of $M_l$ (or even off the
diagonal), since each entry of this form has $\eta_l>\eta_{l-1}$.
The number of choices of such sequences is $\binom{j}{N}$ which is
a polynomial $Q$ in $j$ of degree $N$.  The norm of
$M_{\alpha\beta}$ is then bounded by $C_1^jQ(jC_2)$ as desired.
\end{proof}

\begin{proof}[Proof of Lemma \ref{lemma:estimateoncompostion}]
We use the fact that
$J^k(X,\ff){\simeq}{\bigoplus}_{j=0}^k(S^j(T{\ff}^*)).$ Given
$\phi{\in}\Diff^k(X,\ff)$ and $x{\in}X$, we can write
$j^k(\phi)(x)$ with respect to bases of $J^k(X,\ff)$ at $x$ and
$\phi(x)$ which respects this splitting.  Then it is clear that
$j^k(\phi)(x)$ is block upper triangular, where the diagonal
blocks are of the form $S^j(D\phi)(x)=S^j(j^1(\phi))(x)$ where
$j=0,{\ldots},k$. Therefore the norm of the blocks is bounded by
$N_1^k$, and Lemma \ref{lemma:estimateoncompostion} is an
immediate consequence of Lemma \ref{lemma:blockdiagonal}.
\end{proof}

\section{Locally compact groups and free topological groups}
\label{appendix:lc}

In this appendix we sketch a proof of the following Proposition.
We believe this Proposition to be well-known, and experts we
consulted all provided proofs more or less along the following
lines, but none could provide a reference.

\begin{proposition}
\label{lemma:lc} Let $\G$ be a locally compact, $\sigma$-compact
topological group and $K$ a compact generating set containing a
neighborhood of the identity.  Then the group $\G'$, generated by
$K$ and satisfying all the relations of $\G$ of the form $xy=z$
where $x,y$ and $z$ are in $K$ can be given a topology as a
locally compact, $\sigma$-compact group.
\end{proposition}

We first note the following lemma.

\begin{lemma}
\label{lemma:leftinvariant} Let $\G$ be a locally compact,
$\sigma$-compact group.  Then $\G$ admits a left invariant metric
$d_L$ which defines the topology.
\end{lemma}

\begin{proof}
It is well-known that $\G$ is completely regular and therefore
metrizable and so we can choose a metric $d$ on $\G$, which is not
necessarily left invariant. To find a left invariant metric we
choose a continuous function $f$ supported on a compact
neighborhood $C$ of the identity such that $f$ is continuous,
$f=1-d(x,y)$ on a smaller neighborhood of the identity and
$f{\equiv}0$ outside $C$. We then define a map from $\G$ to
continuous function on $\G$ by $\g{\rightarrow}\g{\cdot}f$ and
define a function $d_L$ on $\G{\times}\G$ by
$d_L(\g_1,\g_2)=\|\g_1{\cdot}f-\g_2{\cdot}f\|_{C^0}$.  This is
clearly left invariant and it is easy to check that $d_L$ is a
metric and defines the same topology as $d$ on $\G$.
\end{proof}

\noindent It is clear from the construction of $d_L$ that we can
normalize so that the ball of radius one is contained in $K$.

To prove the proposition, we need to define a topology on $\G'$.
The group $\G'$ as a group is the quotient of the free group
$F(K)$ and we call the projection map $\pi$.  We can define a norm
on $F(K)$ by letting
$\|k_1k_2{\ldots}k_n\|_F=\sum_{i=1}^nd_L(e_{\G},k_i)$ and then
define
$$\|\g\|_{\G'}=\inf_{\{w{\in}F|\pi(w)=\g'\}}\|w\|_F.$$
The fact that $\|{\cdot}\|$ is a norm is straightforward and we
define a topology on $\G'$ by taking a system of neighborhoods of
a point $x$ to be sets of the form $$\{y| \hskip 2pt
\|x{\inv}y\|_{\G'}<\frac{1}{n}\}.$$ It remains to check that this
defines a topology on $\G'$ that makes $\G'$ a topological group,
and we will indicate a proof of this below, though it also follows
easily from results in \cite{Ma}. For more details on norms,
topologizing topological groups via norms and a construction of a
topology on $F(K)$ which makes $\pi$ continuous, see \cite{Ma}. We
now note an essentially trivial lemma.  We leave the proof to the
reader.

\begin{lemma}
\label{lemma:onK} For any $\g$ in $\G'$ with $\|\g\|_{\G'}<1$, we
can write $\g=\pi(k)$ for some $k{\in}K$ such that
$\|\g\|_{\G'}=d(e_{\G},k)$.
\end{lemma}

Let $U$ be the set of words in the free group with $\|w\|_F<1$. As
immediate consequences of the lemma we have:
\begin{enumerate}
\item the map from $\G'{\rightarrow}\G$ is a homeomorphism on the
set $\pi(U)$ in $\G'$ \item the set $\pi(U)$ contains a
neighborhood of the identity in $\G'$.
\end{enumerate}

It only remains to check that the topology we have defined on
$\G'$ makes $\G'$ a topological group. To see this one merely
needs to check that the topology is invariant under conjugation in
some neighborhood of the identity.  To check this, it suffices to
check it for conjugation by elements of $K$, but there it is more
or less obvious, as the  action of $K$ by conjugation on a small
enough neighborhood of the identity contained in $\G'$ is now
easily seen to be conjugate by a homeomorphism to the action of
$K$ by conjugation on a small neighborhood of the identity in
$\G$.

\section{Historical remarks, relations to other work, and further
generalizations}
\label{appendix:remarks}

This appendix attempts to clarify the relationship of our work to
the work of others
and also contains some remarks that may be
useful for future generalizations of our results.

\subsection{Proofs of Theorem \ref{theorem:isomrigid} and \ref{theorem:almostconjugacygen} and KAM theory}
\label{appendix:kam} This subsection first discusses the failed
proof mentioned in the introduction to this paper, and then goes
one to compare that failed proof, the current successful one, and
the KAM method.

We recall our original approach to proving Theorem
\ref{theorem:isomrigid}. Given an isometric action $\rho$ of
$\Gamma$ on a compact manifold $X$ and a perturbation $\rho'$ of
$\rho$, a conjugacy is a diffeomorphism $f:X{\rightarrow}X$ such
that $\rho(\gamma){\circ}f=f{\circ}\rho'(\gamma)$ for all $\gamma$
in $\Gamma$. Rearranging, the conjugacy is a fixed point for the
$\Gamma$ action on the group $\Diff^k(X)$ of diffeomorphisms of
$X$ defined by
$f{\rightarrow}{\rho(\gamma)}{\circ}f{\circ}\rho'(\gamma){\inv}$.
Ideally we would parameterize diffeomorphisms  of $X$ locally as a
Hilbert space and then use Theorem \ref{theorem:fixedpointpartial}
below to find a fixed point or conjugacy.

We briefly describe an approach to this parametrization and the
difficulty encountered.  Let $\Vect^k(X)$ be the set of vector
fields on $X$ and $C^k(X,X)$ be the set of $C^k$ maps from $X$ to
$X$. Given a Riemannian metric on $X$ there is a natural
exponential map $\Exp:\Vect^k(X){\rightarrow}C^k(X,X)$ defined by
taking a vector $V$ to the time one map of the geodesic flow along
$V$ and projecting back to $X$, i.e. by
$V{\rightarrow}(x{\rightarrow}\exp_xV_x)$. If $\rho=\rho'$ and we
define $\Exp$ using the $\rho$ invariant metric, we have a natural
action of $\Gamma$ on $\Vect^k(X)$ such that $\Exp$ is
equivariant. As will be shown in section
\ref{section:parametrizing} it is fairly straightforward to
complete $\Vect^k(X)$ with respect to a Sobolev metric in such a
way that the completion $\Vect^{2,k}(X)$ is a Hilbert space on
which the $\Gamma$ action defined by $\rho$ is isometric. We had
hoped to show that if $\rho'$ is close enough to $\rho$ then we
would have a partially defined $\Gamma$ action on $\Vect^{2,k}(X)$
that was by $\varepsilon$-almost isometries.  We would then apply
Theorem \ref{theorem:fixedpointpartial} to this partially defined
action to find a fixed vector field $V$, and $\Exp(V)$ would be
the desired conjugacy. While it is possible to construct a
partially defined action on a ball in $\Vect^k(X)$, we were unable
to show that the action is $\varepsilon$-almost isometric if one
considers a metric on $\Vect^{2,k}(X)$ with $k>1$. This is
important, since to show that $V$ and $\Exp(V)$ are smooth, and
that $\Exp(V)$ is invertible, one needs to use the Sobolev
embedding theorems, which require a loss of derivatives
proportional to the dimension of $X$. This method fails even if we
could use an $L^p$ type Sobolev space $\Vect^{p,k}$ and Corollary
\ref{corollary:contractingoperatorbanach}, since at most these
results will produce a continuous invariant vector field $V$, and
it is not clear that $\Exp(V)$ is even a homeomorphism.  We remark
that we cannot use Corollary
\ref{corollary:contractingoperatorbanach} since the action on
$\Vect^k(X)$ is not linear.

It is worth noting that this is different than the difficulty with
loss of derivatives usually encountered by $KAM$ type methods.
Here the problem is that no matter what topology we assume $\rho$
and $\rho'$ are close in,  we cannot prove that the action we
define on vector fields satisfies {\em any estimate of any kind}
on higher derivatives.   In KAM the typical problem is that
estimates for the solutions to the linearized equation are only
``uniformly good" for low order, but one has some a priori
estimate at higher order.

 This difficulty arises from the fact that the method suggested does not
actually involve linearizing the action on diffeomorphisms.  What
we had hoped to do was to use the parametrization mentioned to
provide a linear structure in which the action constructed was
``almost isometric", and therefore ``close enough" to linear, so
as to be able to apply Theorem
\ref{theorem:decreasingdisplacementgen}. The problem is that we
can only do this in function spaces where this yields no
meaningful results.

To resolve this difficulty, we linearize the problem, and in fact
give two different linearizations. Our linearizations are not very
similar to the $KAM$ linearization, and their utility depends
heavily on our results concerning groups with property $(T)$. More
or less our method takes advantage of the fact that, for groups
with property $(T)$, contracting properties of certain operators
are preserved under small perturbations for actions on a wide
variety of uniformly convex Banach spaces. The disadvantage of our
method is that to obtain such a perturbation, we need to only
consider Banach spaces whose definition involve only finitely many
derivatives.  For a long time, this left a $C^{\infty,\infty}$
result out of reach.  Our proof of the $C^{\infty}$ case was
inspired by a study of the $KAM$ method and particularly of the
paper \cite{DK}, but the only concrete similarities to $KAM$
arguments is the use of an iteration and the types of estimates
used.

From the point of view of $KAM$ theory it is surprising that we
need the estimates from Lemma \ref{lemma:estimateoncompostion}
given the strong contracting properties of the averaging operators
we consider.  The need for these estimates is explained following
Corollary \ref{corollary:estimateoncompostion} and Proposition
\ref{proposition:almostisometric}. It is possible to give a proof
of the $C^{\infty,\infty}$ case of Theorem \ref{theorem:isomrigid}
without using these estimates. This has been done very recently by
the first author using Hamilton's implicit function theorem and an
approach similar to Weil's work on local rigidity of lattices in
Lie groups \cite{F}. This approach has applications to local
rigidity of isometric actions for some groups that do not have
property $(T)$, but is unlikely to yield a result in the
generality of Theorem \ref{theorem:almostconjugacygen}. The proof
uses many facts concerning harmonic analysis on compact manifolds
that are unknown, unlikely to be true, or known to be false in the
context of general compact foliated spaces.

\subsection{Further fixed point properties and relations to the
work of M.Gromov} \label{appendix:gromov}
 Examination of the proof of Theorem
\ref{theorem:decreasingdisplacementgen} shows that one can state
more general variants of the theorems discussed here. The limiting
procedure applied in the proof is quite flexible, and allows one
to limit over almost any set of parameters.  To some extent this
is illustrated in the proofs of the results from subsection
\ref{subsection:TresultsBanach}. In another direction, one can
replace $\fh$ by a non-positively curved space that is
``$\varepsilon$-almost flat". By this one should mean anything
that implies that, given a sequence of ``$\varepsilon$-almost
flat" spaces with $\varepsilon$ tending to zero, the limit space
constructed by the method of subsection
\ref{subsection:limitactions} is a Hilbert space.   To actually
prove this variant, we need to define the operator $\rho(h)$ for
actions (or partially defined actions) of $\G$ on spaces of
non-positive curvature.  A method for doing this is described in
subsection \ref{subsection:barycenter} of the appendix.  Theorem
\ref{theorem:fixedpointpartial} can be generalized even further to
``$\varepsilon$-almost flat" spaces which are not non-positively
curved see the discussion related to Lemma \ref{lemma:key} below.
For discrete groups, these more general assertions are easy
exercises from the proofs in section \ref{section:Tproofs} below.
For non-discrete groups, the issue of finding a continuous
subaction of the limit action constructed in subsection
\ref{subsection:limitactions} can present non-trivial difficulties
or require additional assumptions.

In \cite{Gr2}, Gromov proves that certain ``random" infinite,
discrete groups have a fixed point property that is stronger than
property $(T)$. He proves that these groups have fixed points for
any isometric action on any finite or infinite dimensional
``regular" non-positively curved space. After having completed an
earlier draft of this paper, we discovered that the ideas in
\cite{Gr2} have many points in common with ours. In particular, in
section 3.13B, Gromov outlines a proof of a special case of
Theorem \ref{theorem:fixedpointsimplegen},  for a certain class of
``random" infinite, discrete groups with property $(T)$ and for
affine actions. This is a class of groups whose Cayley graphs
``contain" a family of expander graphs as subsets. By a graph
being contained in the Cayley graph, we mean that the Cayley graph
contains an embedded copy of the graph.  Actually, the Cayley
graphs of Gromov's groups only contain ``most" of the relations
that would arise from containing the collection of expander graphs
in a sense made precise in \cite{Gr2}. By a family of expander
graphs we mean a collection of $(n,k,c)$ expanders with $k$ and
$c>0$ fixed and $n$ going to infinity.   Although one can build a
family of expander graphs of this kind as a series of quotients of
any residually finite group with property $(T)$, it is far from
clear that one can realize a family of expanders as subsets of the
Cayley graph for an arbitrary discrete group with property $(T)$,
even in Gromov's probabilistic sense.

More generally, a central philosophy of \cite{Gr2} is that if a
collection of spaces $\mathcal C$ is ``closed under scaling
limits" then for a group $\G$ to have almost fixed points (i.e.
sequences of points with $K$-displacement converging to zero) for
all isometric actions on spaces in $\mathcal C$ is equivalent to
having fixed points for all such actions.  From this point of view
the emphasis of our results on groups with property $(T)$ is on
extending the fixed point property to (partially defined) actions
that are close enough to being isometric on spaces that are
``close enough" to $\mathcal C$.  Finally, we note that it should
be possible to prove a common generalization, and show that
Gromov's groups have fixed points for partially defined
$\varepsilon$-almost isometric actions on ``regular"
non-positively curved spaces.

A primary technical difference between our work and Gromov's is
the functional used.  Where we use the $K$ displacement, Gromov
uses a $K$ energy.  Despite this, various variants of Proposition
\ref{proposition:decreasingdisplacementstandard}, for discrete
groups with property $(T)$, permeates section $3$ of \cite{Gr2},
see particularly $3.8-3.13$. The precise formulations given there
are somewhat more complicated because they are phrased in terms of
energy rather than displacement. For a gentler presentation of
some of the ideas in \cite{Gr2}, see the commentaries \cite{Si}
and section $6$ of \cite{Gh}.

{\noindent}{\bf On the use of ultrafilters:} It is possible to
construct the limit isometric action of $\rho$ on a Hilbert space
$\fh$ ``by hand" without using ultrafilters, at least when $\G$ is
discrete.  To do this, one chooses an explicit isometric
identification of the orbits $\rho_n(\G)x_n$ with subsets of a
fixed Hilbert space $\fh$, always identifying $x_n$ with $0$.  By
passing to a subsequence were $\{\rho_n(\g)x_n\}$ converges for
every $\g$, we can obtain an isometric action of $\G$ on a
countable set in $\fh$ that extends to an action on a closed
linear subspace of $\fh$. Verifying this and then obtaining the
contradiction between the properties of the $\G$ action on $\fh$
and the $\G$ action on $\fh_n$ is considerably more involved than
the proof above, though the argument does not use much more than
simple linear algebra and geometry. The argument is similar to the
proof of Proposition
\ref{proposition:convolutionalmostisometrichilbert} and uses some
of the same lemmas.  It does not seem possible to carry out
arguments of this type in the generality of subsection
\ref{subsection:limitactions}.

\noindent {\bf Fixed points without iterative method:} If one is
more interested in Theorem \ref{theorem:fixedpointpartial} than
Theorem \ref{theorem:decreasingdisplacementgen}, it is possible to
provide an independent proof of that theorem along the same lines.
This may be useful for generalizations to spaces where the
operators $\rho(h)$ are either not defined or not well-behaved. To
do so one needs to produce a Cauchy sequence of points with
smaller and smaller displacement for all partially defined actions
which are ``close enough" to being isometric actions.  The
following lemma, stated by the second author in Jerusalem in 1997,
suffices:

\begin{lemma}
\label{lemma:key} Let $\Gamma$ be a  group with property $(T)$ and
fix a compact generating set $K$. Given $\delta_0>0$ there exist
$\varepsilon>0,r=r(\delta_0)>0$, and positive integers $s$ and
$M$, such that for any Hilbert space $\fh$, any $\delta<\delta_0$
and any $x{\in}\fh$, and any continuous
$(r,s,\varepsilon,\delta,K)$-action of $\Gamma$ on $\fh$ one can
find a point $y$ such that:
\begin{enumerate}
\item $d(x,y){\leq}M\dk(x)$ and, \item
$\dk(y){\leq}\frac{1}{2}\dk(x)$.
\end{enumerate}
\end{lemma}

\noindent To prove Lemma \ref{lemma:key} one argues by
contradiction as in the proof of Theorem
\ref{theorem:decreasingdisplacementgen}. The contradiction follows
since if $\{x_n\}$ is our sequence of basepoints, then there is a
fixed point $y_{\omega}$ in the limit action with
$d(x_{\omega},y_{\omega})$ less than $M$ times the $K$
displacement of $x_{\omega}$, where $M>0$ is a constant depending
only on $\G$ and $K$.  (This fact for isometric actions on Hilbert
spaces is, for example, an easy corollary of Proposition
\ref{proposition:decreasingdisplacementstandard}.)  To prove
Theorem \ref{theorem:fixedpointpartial} one then argues as in the
proof that Theorem \ref{theorem:decreasingdisplacement} implies
Theorem \ref{theorem:fixedpointsimplegen}.

Lemma \ref{lemma:key} suffices to prove Theorems
\ref{theorem:fixedpointsimplegen} and
\ref{theorem:fixedpointpartial}, but these results do not suffice
for our applications.  In particular, we need the precise
iterative method of finding fixed points:
\begin{enumerate}
\item to obtain optimal regularity by finding estimates in $L^p$
type Sobolev spaces for large $p$, \item to control the
non-uniformities that arise in applying the Sobolev embedding
theorems on foliated spaces and, \item to be able to use the
estimates of Section \ref{section:convexityofderivatives} to
obtain $C^{\infty,\infty}$ results.
\end{enumerate}

\bigskip
\noindent
David Fisher\\
Department of Mathematics and Computer Science\\
Lehman College - CUNY\\
250 Bedford Park Boulevard W\\
Bronx, NY 10468\\
dfisher@lehman.cuny.edu\\
\\
Gregory Margulis\\
Department of Mathematics\\
Yale University\\
P.O. Box 208283\\
New Haven, CT 06520\\
margulis@math.yale.edu\\


\begin{thebibliography}{Palais}



\bibitem[Ab]{Ab} Abels, Herbert {\it Finite presentability of $S$-arithmetic groups.
Compact presentability of solvable groups.} Lecture Notes in
Mathematics, 1261. Springer-Verlag, Berlin, 1987.

\bibitem[BFGM]{BG} U.Bader, A.Furman, T.Gelander, N.Monod, Property $(T)$ and actions on $L^p$ spaces, preprint.

\bibitem[BL]{BL} Y.Benjamini, J.Lindenstrauss, {\it Geometric
Non-linear Functional Analysis: Volume 1}, Colloquium
Publications, AMS, 2000.

\bibitem[Be]{Be}
E.J.Benveniste, Rigidity of isometric lattice actions on compact Riemannian manifolds,
{\it GAFA} 10 (2000) 516-542.




\bibitem[BTG]{BTG}
N. Bourbaki, \'{E}l\'{e}ments de math\'{e}matique. Topologie
g\'{e}n\'{e}rale. Chapitres 1 à 4. Hermann, Paris, 1971.









\bibitem[CC]{CC}
A.Candel and L.Conlon, {\it Foliations I}, American Mathematical
Society Graduate Studies in Mathematics, Volume 23, Providence,
2000.

\bibitem[doC]{doC}
M.P. do Carmo,  {\em Riemannian Geometry}, Birkh\"{a}user, Boston,
1992.

\bibitem[De]{De} Delorme, Patrick $1$-cohomologie des représentations unitaires des
groupes de Lie semi-simples et résolubles. Produits tensoriels
continus de représentations. (French){\it Bull. Soc. Math. France}
105 (1977), no. 3, 281--336.

\bibitem[D]{D}
Dieudonne, {\it Elements d'Analyse}, vol. 4, Gauthiers-Villars, Paris, 1971.

\bibitem[DK]{DK} D. Dolgopyat and R. Krikorian, On simultaneous linearization of diffeomorphisms
of the sphere, preprint.






\bibitem[F]{F} D.Fisher, Harmonic analysis, Hamilton's implicit
function theorem and local rigidity, in preparation.

\bibitem[FM1]{FM1}
D.Fisher and G.A.Margulis, Local rigidity for cocycles,in {\em
Surv.~Diff.~Geom.~Vol VIII}, refereed volume in honor of Calabi,
Lawson, Siu and Uhlenbeck , editor: S.T. Yau, 45 pages, 2003.

\bibitem[FM2]{FM2}
D.Fisher and G.A.Margulis, Local rigidity of affine actions of
higher rank groups and lattices, preprint.

\bibitem[Gh]{Gh} E.Ghys, Groupes al\'{e}atoire [d'apr\'{e}s Misha Gromov], to appear {\em Sem. Bourbaki}.

\bibitem[Gr1]{Gr1}
M. Gromov,  Asymptotic invariants of infinite groups.{\it
Geometric group theory, Vol. 2 (Sussex, 1991)}, 1--295, London
Math. Soc. Lecture Note Ser., 182, Cambridge Univ. Press,
Cambridge, 1993.

\bibitem[Gr2]{Gr2} M.Gromov, Random walk in random groups, {\em GAFA} 13 (2003) vol 1, 73-146.


\bibitem[Gu]{Gu}
Guichardet, Alain Sur la cohomologie des groupes topologiques. II.
(French) {\it Bull. Sci. Math.} (2) 96 (1972), 305--332.

\bibitem[HV]{HV}
P. de la Harpe, A.Valette, {\it La propriete (T) de Kazhdan pour les groupes localement compacts},
 Asterisque 175, Soc. Math. de France, Paris, 1989.

\bibitem[He]{He}
Heinrich, S, Ultraproducts in Banach space theory, {\em J. Reine
Angew. Math.} 313 (1980) 72-104.



\bibitem[H]{H}
Hirsch, Morris W.  {\it Differential topology.}  Graduate Texts in Mathematics, No. 33.
Springer-Verlag, New York-Heidelberg, 1976.

\bibitem[Ho]{Ho}
H\"{o}rmander, Lars The boundary problems of physical geodesy.
{\it Arch. Rational Mech. Anal.} 62 (1976), no. 1, 1--52.









\bibitem[K]{K} D. Kazhdan, On the connection of the dual space of a group with the structure of its closed subgroups.
(Russian) {\it Funct. Anal. Appl.} 1 1967 71--74.

\bibitem[KM]{KM}
A.Karlsson and G.A.Margulis. A multiplicative ergodic theorem and
nonpositively curved spaces. {\it Comm. Math. Phys.} 208 (1999),
no. 1, 107--123.



\bibitem[KN]{KN}
S.Kobayashi and K.Nomizu, {\it Foundations of Differential
Geometry, Vol. I}, Wiley-Interscience Publications, New York,
1963.







\bibitem[M]{M}
G.A. Margulis, {\it Discrete subgroups of semisimple Lie groups},
Springer-Verlag, New York, 1991.

\bibitem[Ma]{Ma}
Markov, A. A. Three papers on topological groups: I. On the
existence of periodic connected topological groups. II. On free
topological groups. III. On unconditionally closed sets. {\it
Amer. Math. Soc. Translation}, (1950). no. 30.

\bibitem[MU]{MU}
S.Mazur, S.Ulam, Sur les transformations isomtriques d'espaces
vectoriels normes, {\it C.R. Acad. Sci. Paris} 194 (1932) 946-948.

\bibitem[MS]{MS}
C.C. Moore and C. Schochet, {\em Global analysis on foliated
spaces.} With appendices by S. Hurder, Moore, Schochet and Robert
J. Zimmer. Mathematical Sciences Research Institute Publications,
9. Springer-Verlag, New York, 1988.

\bibitem[Mo]{Mo}
Mostow, G. D. Equivariant embeddings in Euclidean space. {\it Ann.
of Math.} (2) 65 (1957), 432--446.










\bibitem[P1]{Palais}
Palais, R.S., Differential operators on vector bundles, {\it Seminar on the Atiyah-Singer index theorem},
Annals of Math. Studies 57, Princeton University Press, Princeton, 1965.

\bibitem[P2]{P2} Palais, Richard S. Imbedding of compact, differentiable
transformation groups in orthogonal representations. {\em J. Math.
Mech.} 6 (1957), 673--678.











\bibitem[S]{S}
Y.Shalom, Rigidity of commensurators and irreducible lattices.
{\em Invent. Math.} 141 (2000), no. 1, 1--54.

\bibitem[Si]{Si}
L.Silberman, Addendum by L. Silberman To ''Random walk in random
group" by M. Gromov, {\em GAFA} 13 (2003), no. 1, 147-177.






\bibitem[vDW]{vDW}
van den Dries, L.; Wilkie, A. J. Gromov's theorem on groups of
polynomial growth and elementary logic. {\it J. Algebra} 89
(1984), no. 2, 349--374.





\bibitem[Z1]{Z1}
R.J. Zimmer, Volume preserving actions of lattices in semisimple
groups on compact manifolds, {\it Pub. Math. IHES}, 59 (1984)
5-33.


\bibitem[Z2]{Z1.5}
R.J.Zimmer, Lattices in semisimple groups and distal geometric
structures, {\it Invent. Math.} 80 (1985) 123-137.

\bibitem[Z3]{Z3}
R.J.Zimmer, Lattices in semisimple groups and invariant geometric
structures on compact manifolds, {\it Discrete Groups in Geometry
and Analysis:  Papers in Honor of G.D.Mostow}, ed. Roger Howe,
Birkh\"{a}user, Boston, 1987.





\bibitem[Zk]{Zk} A.Zuk, Property (T) and Kazhdan constants for discrete groups, {\em GAFA} 13 (2003) no. 3 643-670.


\end{thebibliography}
\end{document}